\documentclass[11pt,twoside]{article}

\usepackage{latexsym,amssymb,amsthm,amsmath,amscd}

\theoremstyle{plain}

\newtheorem{theorem}{Theorem}[section]

\newtheorem{proposition}{Proposition}[section]

\newtheorem{corollary}{Corollary}[section]

\newtheorem*{Maximum Principle}{Maximum Principle}
\newtheorem*{WP1}{World Problem 1}
\newtheorem*{WP2}{World Problem 2}

\theoremstyle{remark}
\newtheorem{remark}{Remark} [section]

\theoremstyle{definition}
\def\Cal{\mathcal }
\def\<{\left < }
\def\>{\right >}
\def\({\left ( }
\def\){\right )}
\def\e{\epsilon }

\def\e{\eqref}
\def\de{\delta} \def\wt{\tilde}
\def\cn{\hbox{\rm cn}}

\def\dn{\hbox{\rm dn}}

\def\i{{\rm i}\hskip.01in}
\def\rme{{\rm e}\hskip.01in}

\def\n2{\left[{n\over2}\right]}

\def\r{\hat \theta_k(p) }

\newtheorem{definition}{Definition}[section]
\newtheorem{example}{Example}[section]
\numberwithin{equation}{section}

\begin{document}

\vskip.3in

\thispagestyle{empty}

\centerline{\bf {\LARGE $\delta$-Invariants, Inequalities of Submanifolds}}\medskip

\centerline{\bf {\LARGE and Their Applications}}\medskip

\centerline{\it Dedicated to Prof. Leopold Verstraelen on the occasion of his $60$th birthday}\bigskip

\centerline{\bf {\large Bang-Yen Chen}}\medskip

\centerline{\it Department of Mathematics, Michigan State University}

\centerline{\it East Lansing, Michigan $48824$--$1027$, U.S.A.}

\centerline{{\it E-mail}: {\tt bychen@math.msu.edu}}

\footnote[0]{{\it Topics in Differential Geometry}, 29-155, Ed. Acad. Rom\^{a}ne, Bucharest, 2008 (Edited by A. Mihai, I. Mihai and R. Miron).}

\begin{abstract}

The famous Nash embedding theorem was aimed for in the hope that  if  Riemannian  manifolds  could  be regarded as Riemannian submanifolds, this would then yield the opportunity to use extrinsic help. However, as late as 1985 (see \cite{G}) this hope had not been materialized.   The main reason for this is due to the lack of controls of the extrinsic properties of the submanifolds by the known intrinsic invariants. In order to overcome such difficulties as well as to provide answers to an open question on minimal immersions,  we  introduced in the early 1990's new types of Riemannian invariants, known as the $\delta$-invariants or the so-called Chen invariants,  different  in nature from  the ``classical''  Ricci and scalar curvatures. At the same time we also able to establish general optimal relations between the new intrinsic invariants and the main extrinsic invariants  for Riemannian submanifolds. Since then many  results  concerning these invariants, inequalities,  related subjects, and their applications have been obtained by many geometers. 

The main purpose of  this article is thus to provide an extensive and comprehensive survey of  results over  this very active field of research done during the last fifteen years.  Several related inequalities and their applications are  presented in this survey article as well. 
\end{abstract}\newpage

{\it {\large Contents}}\smallskip

1. Motivation to introduce $\delta$-invariants.

2. Definition of $\delta$-invariants.

3. Relations between $\delta$-invariants and Einstein and conformally flat

\hskip.2in  manifolds.

4. Fundamental inequalities involving $\delta$-invariants.

5. Special cases of  fundamental inequalities.

6. Ideal immersions--best ways of living.

7. Applications of $\delta$-invariants to the estimates of eigenvalues of the

\hskip.2in   Laplacian $\Delta$.

8. Applications of $\delta$-invariants to minimal  immersions.

9. Applications of $\delta$-invariants to Lagrangian and slant  immersions.

10. Applications to rigidity problems.

11. Applications to warped products.

12. Growth estimates of warping functions.

13. A $\delta$-invariant  and its applications to Riemannian submersions.

14. A $\delta$-invariant  and its applications to  Einstein manifolds.

15. A $\delta$-invariant  and its applications to  conformally flat manifolds.

16. A $\delta$-invariant  and its applications to contact manifolds.

17. $\delta(2)$ and Lagrangian submanifolds.

18. $\delta(2)$ and $CR$-submanifolds.

19. $\delta(2)$ and $CMC$ hypersurfaces.

20. $\delta(2)$ and submanifolds of nearly K\"ahler $S^6$.

21. Ricci curvature and its applications.

22. K\"ahlerian $\delta$-invariants and applications to complex geometry.

23. Applications  to affine differential geometry (I): $\delta^{\#}$-invariants.

24. Applications  to affine differential geometry (II): warped products.

25. Applications  to affine differential geometry (III): eigenvalues.

26. $\delta^{\#}(2)$ and affine hypersurfaces.

27. Applications of $\delta$-invariants to general relativity.

28. $k$-Ricci curvature and shape operator.

29. General inequalities (I): $CR$-products.

30. General inequalities (II): $CR$-warped products

31. Inequality involving normal scalar curvature and DDVV conjecture.

32. Inequalities involving scalar curvature for Lagrangian and slant

\hskip.28in   submanifolds.

33. Related articles.

34. Bibliography.\bigskip

\section{Motivation to introduce $\delta$-invariants}

Curvature invariants are the $N^o\,1$  Riemannian invariants  and the most natural ones. Curvature invariants also play key roles in physics. For instance, the magnitude of a force required to move an object at constant speed, according to Newton's laws, is a constant multiple of the curvature of the trajectory. The motion of a body in a gravitational field is determined, according to Einstein's general theory of relativity, by the curvatures of space time. All sorts of shapes, from soap bubbles to red blood cells, seem to be determined by various curvatures (cf. \cite{Os}). Borrowing a term from biology, Riemannian invariants are the DNA of Riemannian manifolds.  Classically, among the Riemannian curvature invariants, people have been studying sectional, scalar and Ricci curvatures in great detail.  

One of the most fundamental problems in the theory of submanifolds is the immersibility (or non-immersibility) of a Riemannian manifold in a Euclidean space (or, more generally, in a space form). According to the 1956 celebrated embedding theorem of J. F. Nash \cite{nash}, every Riemannian manifold can be isometrically embedded in some Euclidean spaces with sufficiently high codimension.  

The Nash embedding theorem was aimed for in the hope that  if  Riemannian  manifolds  could always be regarded as Riemannian submanifolds, this would then yield the opportunity to use extrinsic
help. However, this hope had not been materialized as late as 1985 according to M. Gromov \cite{G} (see also \cite{c-nash}).   The main reason for this is due to the lack of controls of the extrinsic properties of the submanifolds by the known intrinsic invariants.

In view of Nash's theorem,  to study  embedding problems  it is natural to impose some suitable condition(s) on the immersions. For example, if one  imposes the minimality condition, it leads to
\vskip.1in
\noindent {\bf Problem 1.} {\it Given a Riemannian manifold $M$, what are necessary  conditions for $M$ to admit a minimal isometric immersion in a Euclidean $m$-space $\mathbb E^m$}?
\vskip.1in
It is well-known that for a minimal submanifold in $\mathbb E^m$, the Ricci tensor satisfies Ric $\leq 0$.  For many years this was the only known necessary Riemannian condition for a general Riemannian manifold to admit a minimal isometric immersion in a Euclidean space regardless of codimension.
That is why  S. S. Chern asked  in his 1968 monograph to search for further   Riemannian obstructions  for $M$ to admit  an isometric minimal immersion into a Euclidean space. Also, no solutions to Chern's problem were known for many years before the invention of the $\delta$-invariants.

In order to overcome those difficulties,  we need to introduce certain new types of Riemannian invariants, different  in nature from  the ``classical''  invariants. Moreover, we also need to  establish general optimal relationships between the main extrinsic invariants  with  the new  intrinsic invariants  on the submanifolds. These are the author's original motivation in 1990's to introduce his so-called $\delta$-invariants on Riemannian manifolds. 

The $\delta$-invariants are {\it very different in nature} from the ``classical'' scalar and Ricci curvatures; simply due to the fact that both scalar and Ricci curvatures are {\it ``total sum''} of sectional curvatures on a Riemannian manifold. In contrast, all of the non-trivial $\delta$-invariants are obtained  from the scalar curvature by throwing away a certain amount of sectional curvatures. 

\section{Definition of $\delta$-invariants}

Let $M$ be a  Riemannian $n$-manifold. Denote by $K(\pi)$ the sectional curvature of $M$ associated with a plane section $\pi\subset T_pM$, $p\in M$. For any orthonormal basis $e_1,\ldots,e_n$ of the tangent space $T_pM$, the scalar curvature $\tau$ at $p$ is defined to be $$\tau(p)=\sum_{i<j} K(e_i\wedge e_j).\notag$$

Let $L$ be a subspace of $T_pM$  of dimension $r\geq 2$  and $\{e_1,\ldots,e_r\}$ an orthonormal basis of $L$. We define the scalar curvature $\tau(L)$ of the $r$-plane section $L$ by  $$\tau(L)=\sum_{\alpha<\beta} K(e_\alpha\wedge e_\beta),\quad 1\leq
\alpha,\beta\leq r.$$

Given an orthonormal basis $\{e_1,\ldots,e_n\}$ of the tangent space $T_pM$, we simply denote by $\tau_{1\cdots r}$ the scalar curvature of the $r$-plane section spanned by $e_1,\ldots,e_r$. The scalar curvature $\tau(p)$ of $M$ at $p$ is nothing but the scalar curvature of the tangent space of $M$ at $p$; and if $L$ is a 2-plane section,  $\tau(L)$ is nothing but the sectional curvature $K(L)$ of $L$. 

Geometrically, $\tau(L)$ is nothing but the scalar  curvature of the image $\exp_p(L)$ of $L$ at $p$ under the exponential map at $p$.

For an integer $k\geq 0$ denote by ${\mathcal S}(n,k)$ the finite set  consisting of unordered $k$-tuples $(n_1,\ldots,n_k)$ of integers $\geq 2$ satisfying  $n_1< n$ and $n_1+\cdots+n_k\leq n$. Denote by ${\mathcal S}(n)$ the set of unordered $k$-tuples with $k\geq 0$ for a fixed $n$. 

 For each $k$-tuple $(n_1,\ldots,n_k)\in {\mathcal S}(n)$ the  Riemannian invariant $\delta{(n_1,\ldots,n_k)}$ is defined to be
\begin{align}\label{2.1} \delta(n_1,\ldots,n_k)(p)=\tau(p)-\inf\{\tau(L_1)+\cdots+\tau(L_k)\}, \end{align}
where $L_1,\ldots,L_k$ run over all $k$ mutually orthogonal subspaces of $T_pM$ such that  $\dim L_j=n_j,\, j=1,\ldots,k$.

Similarly, we  have also defined $\hat \delta(n_1,\ldots,n_k)(p)$ by
\begin{align}\label{2.2} \hat \delta(n_1,\ldots,n_k)(p)=\tau(p)-\sup\{\tau(L_1)+\cdots+\tau(L_k)\},\end{align}
where $L_1,\ldots,L_k$ run over all $k$ mutually orthogonal subspaces of $T_pM$ such that  $\dim L_j=n_j,\, j=1,\ldots,k$. 

Let $\#{\mathcal S}(n)$ denote the cardinal number of ${\mathcal S}(n)$. Then $\#{\mathcal S}(n)$ increases quite rapidly with $n$. For instance, for
$$\begin{array}{c}n=2,3,4,5,6,7,8,9,10,11,12,\ldots,20,\ldots,\\50,\ldots,100,\ldots,200,\dots,\end{array}$$  
$\#{\mathcal S}(n)$  are given respectively by
$$\begin{array}{c}1,2,4,6,10,14,21,29,41,54,76,\ldots,626,\ldots,\\204225,\ldots, 190569291,\ldots,3972999029387,\ldots\;.\end{array}$$

 In general, the cardinal number $\#{\mathcal S}(n)$ is equal to $p(n)-1$, where $p(n)$ denotes the partition function. The asymptotic behavior of $\#{\mathcal S}(n)$ is
given by
$$\#{\mathcal S}(n)\approx {1\over {4n\sqrt{3}}}\exp\left[\pi\sqrt{\text{\small$ {{2n}\over 3}$}}\, \right]\quad \hbox{as}\;\;n\to\infty.$$
As a consequence, there exist many $\delta$-invariants for Riemannian manifolds of higher dimension. 

In terms of the $\delta$-invariants, the scalar curvature $\tau$ is nothing but $\delta(\emptyset)$ or $\hat \delta(\emptyset)$ (with $k=0$). The simplest non-trivial $\delta$-invariants are $\delta(2)$ and $\hat \delta(2)$. The  scalar curvature and the $\delta$-invariants $\delta(n_1,\ldots,n_k)$ with $ k>0$ differ greatly in nature. 

Obviously, one has $$\delta{(n_1,\ldots,n_k)}\geq \hat\delta{(n_1,\ldots,n_k)}$$ for any $k$-tuple
$(n_1,n_2,\ldots,n_k)\in {\mathcal S}(n)$. A Riemannian $n$-manifold $M$ is  called an {\it $S(n_1,\ldots,n_k)$-space\/} if it satisfies $$\delta{(n_1,\ldots,n_k)}=
\hat \delta{(n_1,\ldots,n_k)}$$ identically for a fixed $k$-tuple $(n_1,\ldots,n_k)\in {\mathcal S}(n)$.

In this article, some other invariants of a similar nature, i.e., those invariants obtained from the scalar curvature by deleting certain amount of sectional curvature,  are also called $\delta$-{\it invariants}. Those invariants have a similar nature as $\delta(n_1,\ldots,n_k)$ or $\hat \delta(n_1,\ldots,n_k)$ with $k>0$.
For instance, we  have the so-called affine $\delta$-invariants, K\"ahlerian $\delta$-invariants, normal $\delta$-invariant, ..., etc.

\section[Relations between $\delta$-invariants]{Relations between $\delta$-invariants and Einstein and conformally flat  manifolds}

The $S(n_1,\ldots,n_k)$-spaces are completely determined by the following two propositions.

\begin{proposition} {\rm \cite{cdvv4}}  Let $M$ be a Riemannian $n$-manifold with $n>2$. Then 
\vskip.04in

$(1)$ For any integer $j$ with $2\leq j\leq n-2$, $M$ is an $S(j)$-space if and only if $M$ is a Riemannian space form.
\vskip.04in

$(2)$ $M$ is an $S(n-1)$-space if and only if $M$ is an Einstein space.\end{proposition}
\vskip.1in

\begin{proposition} {\rm \cite{cdvv4}}  Let $M$ be a Riemannian $n$-manifold such that $n$ is not a prime and $k$ an integer $\geq 2$. Then
\vskip.04in

$(1)$ if $M$ is an $S(n_1,\ldots,n_k)$-space, then $M$ is a Riemannian space form unless $n_1=\ldots=n_k$ and $n_1+\cdots+n_k=n$, and
\vskip.04in

$(2)$  $M$ is an $S(n_1,\ldots,n_k)$-space with $n_1=\ldots=n_k$ and $n_1+\cdots+n_k=n$ if and only if  $M$ is a conformally flat space.\end{proposition}

By using the notion of $\delta$-invariant, we have the following simple characterization of Einstein spaces which generalizes the well-known characterization of Einstein 4-manifolds given by I. M. Singer and J. A. Thorpe \cite{IMS}.

\begin{theorem} {\rm \cite{cdvv4}} \label{T:3.1} Let $M$ be a Riemannian $2r$-manifold. Then $M$ is an Einstein space if and only if we have
\begin{align}\tau(L)=\tau(L^\perp)\end{align} for any $r$-plane section $L\subset T_pM,\, p\in M$.\end{theorem}

Moreover, also by using the notion of the scalar curvature of $r$-plane sections, we have the following simple characterization of conformally flat spaces which generalizes a well-known result  of R. S. Kulkarni \cite{Kul}.

\begin{theorem} {\rm \cite{cdvv4}} Let $M^n$ be a Riemannian manifold with $n\geq 4$, and let $s$ be an integer satisfying $2<2s\leq n$. Then $M$ is a conformally flat manifold if and only if,  for any orthonormal  set $\{e_1,\ldots,e_{2s}\}$ of vectors, one has
$$\tau_{1\cdots s}+\tau_{s+1\cdots 2s}=\tau_{1\cdots s-1 \,s+1}+\tau_{s\,s+2\cdots  2s}.$$
\end{theorem}

 In general, the $\delta$-invariants $\delta(n_1,\ldots,n_k)$ are independent invariants. However, Theorem \ref{T:3.1}
implies that, for a $2r$-dimensional Einstein manifold, we have the
following relations:
\begin{align} 2\delta(r)-\delta(r,r)=2\hat\delta(r)-\hat\delta(r,r).\end{align}

For any $k$-tuple $(n_1,\ldots,n_k)\in{\mathcal S}(n)$, let us put
\begin{align} \label{3.3}\Delta(n_1,\ldots,n_k)= {{\delta(n_1,\ldots,n_k)}\over{c(n_1,\ldots,n_k)}},\end{align}
where $c(n_1,\ldots,n_k)$ is defined by\begin{equation}\begin{aligned} &c(n_1,\ldots,n_k)= {{n^2(n+k-1-\sum n_j)}\over{2(n+k-\sum n_j)}}.\end{aligned}\end{equation}

Since a Riemannian $n$-manifold with $n\geq 3$ satisfies inequality $\Delta(2)>\Delta(\emptyset)=\tau$ if and only if  $\inf K<\tau/(n-1)^2.$ Thus, a Riemannian $n$-manifold $(n\geq 3)$ with vanishing scalar curvature satisfies $$\Delta_0(2)>\Delta_0(\emptyset)$$ automatically, 
unless $M$ is flat. 

For compact homogeneous Einstein K\"ahler manifolds, we also have the following relationship between the $\delta$-invariants and scalar curvature. 

\begin{proposition} {\rm \cite{c12}}  Let $M$ be a compact homogeneous Einstein Kaehler manifold with positive scalar curvature. Then, for each $(n_1,\ldots,n_k)\in {\mathcal S}(n)$, we have
$$\Delta(n_1,\ldots,n_k)\leq  \Big(2-\frac{2}{n}\Big) \Delta(\emptyset),$$ where
$n$ denotes the real dimension of $M$. \end{proposition}

\section{Fundamental inequalities involving $\delta$-invariants}

Let $M$ be an $n$-dimensional submanifold of a
 Riemannian $m$-manifold $\tilde M^m$.   We choose a local field of orthonormal
frame
$$e_1,\ldots,e_n,e_{n+1},\ldots,e_m$$ in $\tilde M^m$ such that, restricted to  $M$, the vectors $e_1,\ldots,e_n$ are tangent to $M$ and hence $e_{n+1},\ldots,e_m$ are normal to $M$.  Let $K(e_i\wedge e_j)$ and $\tilde K(e_i\wedge e_j)$ denote respectively the sectional curvatures of $M$ and $\tilde M^m$ of the plane section spanned by $e_i$ and $e_j$. 

For the submanifold $M$ in  $\tilde M^m$ we denote by $\nabla$ and ${\tilde \nabla}$ the Levi-Civita connections of $M$ and $\tilde M^m$, respectively. The Gauss and Weingarten formulas are given respectively by (see, for instance, \cite{c1})
\begin{align} &{\tilde \nabla}_{X}Y=\nabla_{X} Y + h(X,Y),\\ &{\tilde\nabla}_{X}\xi =
-A_{\xi}X+D_{X}\xi \end{align} for any  vector fields $X,Y$ tangent to $M$ and  vector field $\xi$ normal to $M$, where $h$ denotes the second fundamental form, $D$ the normal connection, and $A$ the shape operator of the submanifold. 

Let $\{h^r_{ij}\}$, $i,j=1,\ldots,n;\,r=n+1,\ldots,m$,  denote the coefficients of the second fundamental form $h$ with respect to $e_1,\ldots,e_n,e_{n+1},\ldots,e_m$. Then we have $$h^r_{ij}=\<h(e_i,e_j),e_r\>=\<A_{e_r}e_i,e_j\>,$$ where $\<\;\,,\;\>$ denotes the inner product.

The mean curvature vector $\overrightarrow{H}$ is defined by
\begin{align}\overrightarrow{H} = {1\over n}\,\hbox{\rm trace}\,h = {1\over n}\sum_{i=1}^{n} h(e_{i},e_{i}), \end{align}
where $\{e_{1},\ldots,e_{n}\}$ is a local  orthonormal frame of the tangent bundle $TM$ of $M$. The squared mean curvature is then given by
$$H^2=\left<\right.\hskip-.02in \overrightarrow{H},\overrightarrow{H}\hskip-.02in\left.\right>.$$ A submanifold $M$  is called  minimal  in $\tilde M^m$ if  its mean curvature vector  vanishes identically. 

Denote by $R$ and $\tilde R$  the Riemann curvature tensors of $M$ and $\tilde M^m$, respectively. Then the {\it equations of Gauss and Codazzi\/} are given respectively by \begin{equation}\begin{aligned} &R(X,Y;Z,W)=\tilde R(X,Y;Z,W) +\<h(X,W),h(Y,Z)\>\\& \hskip1in -\<h(X,Z),h(Y,W)\>,
\\&(\tilde R(X,Y)Z)^\perp= (\bar\nabla_X h)(Y,Z) - (\bar\nabla_Y h)(X,Z),\end{aligned}\end{equation}
where $X,Y,Z,W$ are tangent to $M$ and $\nabla h$ is defined by \begin{align}  \label{2.7}
(\nabla h)(X,Y,Z) = D_X h(Y,Z) - h(\nabla_X Y,Z) - h(Y,\nabla_X Z).\end{align}
 for vectors $X,Y,Z,W$ tangent to $M$. 
 
 A submanifold $M$ is called a {\it parallel submanifold} if we have $\bar \nabla h=0$ identically.

For each $(n_1,\ldots,n_k)\in {\mathcal S}(n)$, let  $b(n_1 \ldots,n_k)$ denote the constant given by
\begin{equation}\begin{aligned} & b(n_1,\ldots,n_k)={1\over2} {{n(n-1)}}- {1\over2}\sum_{j=1}^k {n_j(n_j-1)} .\end{aligned}\end{equation}

For any isometric immersion from a Riemannian submanifold into another Riemannian manifold, we have the following general optimal inequality.

\begin{theorem} \label{T:4.1} {\rm \cite{c39}} Let $\phi:M\to \tilde M$ be an isometric immersion of a Riemannian $n$-manifold into a  Riemannian $m$-manifold. Then, for each point $p\in M$ and each $k$-tuple $(n_1,\ldots,n_k)\in \mathcal S(n)$,  we have the following inequality:
\begin{align}\label{4.6}\delta{(n_1,\ldots,n_k)}(p) \leq  c(n_1,\ldots,n_k)H^2(p)+b(n_1,\ldots,n_k)\max\tilde K(p),\end{align}  where $\max\tilde K(p)$ denotes the maximum of the sectional curvature function of $\tilde M^m$ restricted to $2$-plane sections of the tangent space $T_pM$ of $M$ at $p$.

The equality case of inequality \eqref{4.6} holds at  $p\in M$ if and only if the following  conditions hold:
\vskip.04in

{\rm (a)} There exists an  orthonormal basis  $e_1,\ldots,e_m$ at $p$, such that  the shape operators of $M$ in $\tilde M^m$ at $p$ take the following  form
 :\begin{align}\label{4.7}\font\b=cmr8 scaled \magstep2 \def\bigzerol{\smash{\hbox{ 0}}} \def\bigzerou{\smash{\lower.0ex\hbox{\b 0}}} A_{e_r}=\left( \begin{matrix} A^r_{1} & \hdots & 0 \\ \vdots  & \ddots& \vdots &\bigzerou \\ 0 &\hdots &A^r_k\\  \\&\bigzerou & &\mu_rI \end{matrix} \right),\quad  r=n+1,\ldots,m, \end{align}
where $I$ is an identity matrix and  $A^r_j$ is a symmetric $n_j\times n_j$  submatrix such  that 
\begin{align}\label{4.8}\hbox{\rm trace}\,(A^r_1)=\cdots=\hbox{\rm trace}\,(A^r_k)=\mu_r.\end{align}

{\rm (b)} For any $k$ mutual orthogonal subspaces $L_1,\ldots,L_k$ of $T_pM$ which  satisfy
$$\delta(n_1,\ldots,n_k)=\tau-\sum_{j=1}^k \tau (L_j)$$ at $p$, we have
$ \tilde K(e_{\alpha_i},e_{\alpha_j})=\max \tilde K(p)$
 for any $\alpha_i\in \Gamma_i,\alpha_j\in \Gamma_j$ with $0\leq i\ne j\leq k$, where \begin{equation}\begin{aligned}\notag &\Gamma_0=\{1,\ldots, n_1\},\\&\ldots\ldots\ldots\ldots\ldots\ldots\\&\Gamma_{k-1}=\{n_1+\cdots+n_{k-1}+1,\ldots, n_1+\cdots+n_k\},\\& \Gamma_{k}=\{n_1+\cdots+n_{k}+1,\ldots, n\}.\end{aligned}\end{equation}
\end{theorem}

\section{Special cases of  fundamental inequalities}

\subsection{Submanifolds in real, complex and quaternionic space forms}
The following results are special cases of  Theorem \ref{T:4.1}.

\begin{theorem} \label{T:5.1}{\rm \cite{c5,c16}}   For each $k$-tuple $(n_1,\ldots,n_k)\in{\mathcal S}(n)$ and for each $n$-dimensional submanifold $M$ in a Riemannian space form $R^m(\epsilon)$ of constant sectional curvature $\epsilon$, we have
\begin{align}\label{5.1}\delta{(n_1,\ldots,n_k)} \leq  c(n_1,\ldots,n_k)H^2+b(n_1,\ldots,n_k)
\epsilon.\end{align}

The equality case of  \eqref{5.1} holds at a point $p\in M$ if and only if there exists an  orthonormal basis  $e_1,\ldots,e_m$ at $p$ such that  the shape operators of $M$ at $p$ take the  forms \e{4.7} and \e{4.8}. \end{theorem}

In particular, for any submanifold $M$ of a Euclidean $m$-space, we have  the following general optimal inequality
\begin{theorem} \label{T:5.2}{\rm \cite{c5,c16}} For any $k$-tuple $(n_1,\ldots,n_k)\in{\mathcal S}(n)$ and any $n$-dimensional submanifold $M$  of a Euclidean space $\mathbb E^m$ with arbitrary codimension, we have
\begin{align}\label{5.2}\delta{(n_1,\ldots,n_k)} \leq  c(n_1,\ldots,n_k)H^2.
\end{align}
\end{theorem}

Since the sectional curvatures of a complex projective space $CP^m(4\epsilon)$ (or quaternion projective space $QP^m(4\epsilon)$) satisfies $\epsilon\leq K\leq 4\epsilon$, 
Theorem \ref{T:4.1} implies  

\begin{theorem}  \label{T:5.3}
Let $M$ be an $n$-dimensional submanifold  of the complex projective $m$-space $CP^m(4\epsilon)$ of constant holomorphic sectional curvature $4\epsilon$  $($or the quaternionic projective $m$-space $QP^m(4\epsilon)$ of quaternionic sectional curvature $4\epsilon)$. Then  we have
\begin{align}\delta{(n_1,\ldots,n_k)}(p) \leq  c(n_1,\ldots,n_k)H^2(p)+4b(n_1,\ldots,n_k)\epsilon\end{align} \noindent for any $k$-tuple $(n_1,\ldots,n_k)\in
\mathcal S(n)$. \end{theorem}

Since the sectional curvatures of a complex hyperbolic space $CH^m(4\epsilon)$ (or quaternion hyperbolic space $QH^m(4\epsilon)$) satisfies $$4\epsilon\leq K\leq \epsilon,$$
Theorem \ref{T:4.1} also gives the following.

\begin{theorem} Let $M$ be an $n$-dimensional submanifold  of the complex hyperbolic $m$-space $CH^m(4\epsilon)$ of constant holomorphic sectional curvature $4\epsilon$ $($or the quaternionic hyperbolic $m$-space $QH^m(4\epsilon)$ of quaternionic sectional curvature $4\epsilon)$. Then  we have
\begin{align}\delta{(n_1,\ldots,n_k)}(p) \leq  c(n_1,\ldots,n_k)H^2(p)+b(n_1,\ldots,n_k)\epsilon\end{align} \noindent for any $k$-tuple $(n_1,\ldots,n_k)\in
\mathcal S(n)$.\end{theorem}

\subsection{Submanifolds in Sasakian space forms}
  
 A $(2m+1)$-dimensional manifold is said to be {\it almost contact\/} if it admits a tensor field $\phi$ of type $(1,1)$, a vector field $\xi$ and a 1-form $\eta$ satisfying \begin{equation}\phi^2=-I+\eta\otimes \xi,\;\;\eta(\xi)=1,\end{equation}  where $I$ is the identity endomorphism. It is well-known that $\phi\xi=0,\;\;\eta\circ\phi=0.$

 An almost contact manifold $(\tilde M,\phi,\xi,\eta)$ is called an {\it almost contact metric manifold\/} if it admits a Riemannian metric $g$ such that
\begin{equation}\begin{aligned}\label{5.4}&g(\phi X,\phi Y)=g(X,Y)-\eta (X)\eta (Y)\;\;\end{aligned}\end{equation} for  vector fields $X,Y$ tangent to $\tilde M$. Setting $Y=\xi$ we have $ \eta(X)=g(X,\xi).$

By a {\it contact manifold} we mean a $(2m+1)$-manifold $\tilde M$ together with  a global 1-form $\eta$  satisfying $$\eta\wedge ({\rm d}\eta)^m \ne 0$$ on $M$. If $\eta$ of  an almost contact metric manifold $(\tilde M,\phi,\xi,\eta,g)$ is a contact form and if  $\eta$ satisfies $${\rm d}\eta (X,Y)=g(X,\phi Y)$$ for all vectors  $X,Y$ tangent to $\tilde M$, then $\tilde M$ is called a {\it contact metric manifold}. 

A contact metric manifold is called {\it $K$-contact} if its characteristic vector field $\xi$ is a Killing vector field.   A $K$-contact manifold  is called  {\it  Sasakian} if
we have $$N_\phi+2{\rm d}\eta\otimes\xi=0,$$ where $N_\phi$ is the Nijenhuis tensor associated to $\phi$. A plane section $\sigma$ in $T_p\tilde M^{2m+1}$ of a Sasakian manifold $\tilde M^{2m+1}$ is called $\phi$-section if it is spanned by $X$ and $\phi(X)$, where $X$ is a unit tangent vector orthogonal to $\xi$. The sectional curvature  with respect to a $\phi$-section $\sigma$ is called a $\phi$-sectional curvature. If a Sasakian manifold has constant $\phi$-sectional curvature, it is called a Sasakian space form.
 
 An $n$-dimensional submanifold $M^n$ of a Sasakian space form $\tilde{M}^{2m+1} (c)$ is called a $C$-totally real submanifold of $\tilde{M}^{2m+1}(c)$ if $\xi $
is a normal vector field on $M^n$. A direct consequence of this definition is that $\phi (TM^n) \subset T^{\bot}M^n$, which means that  $M^n$ is an anti-invariant submanifold of $\tilde{M}^{2m+1}(c)$
 
 It is well-known that the  Riemannian curvature tensor of a Sasakian space form $\tilde M^{2m+1}(\epsilon)$ of constant $\phi$-sectional curvature $\epsilon$ is given by \cite{b}:
\begin{equation}\begin{aligned}\label{5.7} \tilde R&(X,Y)Z=\frac{\epsilon+3}{4}(\<Y,Z\>X-\<X,Z\>Y)
\\& +\frac{\epsilon-1}{4}(\eta(X)\eta(Z)Y-\eta(Y)\eta(Z)X+\<X,Z\>\eta (Y)\xi
\\& -\<Y,Z\>\eta(X)\xi+\<\phi Y,Z\>\phi X-\<\phi X,Z\>\phi Y-2\<\phi X,Y\>\phi Z)
\end{aligned}\end{equation} for $X,Y,Z$ tangent to $\tilde M^{2m+1}(\epsilon)$.
Hence if $\epsilon \geq 1$, the sectional curvature function $\tilde K$ of $\tilde M^{2m+1}(\epsilon)$ satisfies
\begin{align}\label{5.8} \frac{\epsilon+3}{4}\leq \tilde K(X,Y)\leq \epsilon
\end{align} for $X,Y\in \ker \eta$; if $\epsilon<1$, the inequalities are reversed.

 From Theorem \ref{T:4.1} and these sectional curvature properties \eqref{5.7} and \eqref{5.8} of Sasakian space forms, we obtain the following results for arbitrary Riemannian submanifolds in Sasakian space forms.

\begin{corollary}  \label{C:5.1} If $M$ is an $n$-dimensional submanifold  of a Sasakian space form $\tilde M(\epsilon)$ of constant  $\phi$-sectional curvature $\epsilon \geq 1$, then, for any $k$-tuple $(n_1,\ldots,n_k)\in \mathcal S(n)$,   we have
\begin{align}\delta{(n_1,\ldots,n_k)}(p) \leq  c(n_1,\ldots,n_k)H^2(p)+b(n_1,\ldots,n_k)\epsilon.\end{align} 
\end{corollary}

\begin{corollary}   \label{C:5.2}If $M$ is an
$n$-dimensional submanifold  of a Sasakian space form $\tilde M(\epsilon)$ of constant  $\phi$-sectional curvature $\epsilon<1$, then, for any $k$-tuple $(n_1,\ldots,n_k)\in \mathcal S(n)$,   we have
\begin{align}\delta{(n_1,\ldots,n_k)}(p)
\leq  c(n_1,\ldots,n_k)H^2(p)+b(n_1,\ldots,n_k).\end{align} 
\end{corollary}

\begin{corollary}  \label{C:5.3} If $M$ is an $n$-dimensional $C$-totally real submanifold  of a Sasakian space form $\tilde M(\epsilon)$ of constant  $\phi$-sectional curvature $\epsilon \leq 1$, then, for any $k$-tuple $(n_1,\ldots,n_k)\in \mathcal S(n)$,   we have
\begin{align}\delta{(n_1,\ldots,n_k)}(p)
\leq  c(n_1,\ldots,n_k)H^2(p)+b(n_1,\ldots,n_k)\frac{\epsilon+3}{4}.\end{align} 
\end{corollary}

Corollary \ref{C:5.3} has been obtained in  \cite{DMV2}.

\subsection{Lagrangian and totally real submanifolds in complex  space forms}

Since the proof of inequality \e{5.1} is based only on the equation of Gauss, the same inequality holds for  Lagrangian submanifolds (or more generally, totally real submanifolds) in a complex space form. In fact, we have the following inequality for totally real submanifolds (see \cite{c19,cdvv3}).

\begin{theorem}  
Let $M$ be a totally real  submanifold  of  a K\"ahler manifold $\tilde M^m(4\epsilon)$ of constant holomorphic sectional curvature $4\epsilon$. Then  we have
\begin{align}\label{05.10} \delta{(n_1,\ldots,n_k)}(p) \leq  c(n_1,\ldots,n_k)H^2(p)+b(n_1,\ldots,n_k)\epsilon\end{align} \noindent for any $k$-tuple $(n_1,\ldots,n_k)\in \mathcal S(n)$. \end{theorem}
 
 When the totally real submanifolds are Lagrangian, we have the following result \cite{c19,cdvv1}.
 
\begin{theorem}  \label{T:5.5}
Let $M$ be a Lagrangian  submanifold  of  a K\"ahler manifold $\tilde M^n(4\epsilon)$ of constant holomorphic sectional curvature $4\epsilon$. Then  we have
\begin{align}\label{5.10} \delta{(n_1,\ldots,n_k)}(p) \leq  c(n_1,\ldots,n_k)H^2(p)+b(n_1,\ldots,n_k)\epsilon\end{align} \noindent for any $k$-tuple $(n_1,\ldots,n_k)\in \mathcal S(n)$. 

If the equality case of  \e{5.10} holds identically on $M$, then $M$ is a minimal Lagrangian submanifold of $\tilde M^n(4\epsilon)$.
\end{theorem}

 A Lagrangian immersion is said to have  {\it full first normal bundle\/} if the first normal space of $M^n$ equals to the normal space at each point $p\in M^n$, i.e. Im$\,h=T^\perp M^n$.

In \cite{c19}, the author has determined ideal Lagrangian submanifolds in complex space forms as follows (see, also \cite{c41}) (see section 6 for the definition of ideal immersions).

\begin{theorem} If $x:M^n\to \hbox{\bf C}^n$ is a Lagrangian immersion of a Riemannian $n$-manifold into the complex Euclidean $n$-space {\bf C}$^n$ with full first normal bundle, then $x$ is an ideal Lagrangian immersion if and
only if $x$ is locally the product of some minimal Lagrangian immersions with full first normal bundle.\end{theorem}

It is known that there exist ample examples of ideal Lagrangian submanifolds in complex projective and complex hyperbolic spaces. On the contrast, we had proved the following two non-existence results in \cite{c19}. 

\begin{theorem} There do not exist  ideal Lagrangian submanifolds in a complex projective space with full first normal bundle.\end{theorem}

\begin{theorem}  There do not exist  ideal Lagrangian submanifolds in a complex hyperbolic space with full first normal bundle.\end{theorem}

A submanifold $M$ in a Riemannian manifold $N$ is called {\it ruled\/} if at each point $p\in M$, $M$ contains a geodesic $\gamma_P$ of $N$ through $p$.

\begin{theorem} Let $M^n$ be a Lagrangian submanifold of {\bf C}$^n$ such that \hbox{\rm Im}$\,h_p\ne T^\perp_pM^n$ at each point $p\in M^n$. If $M^n$ is ideal, then it is  a ruled minimal submanifold.\end{theorem}

\begin{theorem} Let $M^n$ be a Lagrangian submanifold of a complex space form $\tilde M^n(4c)$ with $c\ne 0$. If $M^n$ is ideal, then it is a ruled minimal submanifold.\end{theorem}

\section{Ideal immersions--best ways of living}

The fundamental inequalities \eqref{4.6} and \e{5.1} give prima controls on the most important extrinsic curvature; namely, the squared mean curvature $H^2$,  by the initial  intrinsic curvatures, the $\delta$-invariants $\delta(n_1,\ldots,n_k)$,  of the Riemannian manifold.

\subsection{A maximum principle}
In general there do not exist direct relationship between  $\delta$-invariants $\delta(n_1,\ldots,n_k)$.  On the other hand, we have the following. 

\begin{Maximum Principle} Let $M$ be an $n$-dimensional submanifold of a Euclidean $m$-space $\mathbb E^m$. If it satisfies the equality case of \eqref{5.2}, i.e., it satisfies
\begin{align} \label{6.2} H^2=\Delta(n_1,\ldots,n_k)\end{align} 
for any $k$-tuple $(n_1,\ldots,n_k)\in {\mathcal S}(n)$, then 
\begin{align} \label{6.3}\Delta(n_1,\ldots,n_k)\geq\Delta(m_1,\ldots,m_j).\end{align} 
 \end{Maximum Principle}

For any isometric immersion  $x:M\to  \mathbb E^m$ of a Riemannian $n$-manifold $M$ in $\mathbb E^m$. Theorem \ref{T:5.1} yields
\begin{align}\label{6.4} H^2(p)\geq \hat\Delta_0(p),\end{align} 
where $\hat\Delta_0$ is the Riemannian invariant on $M$ defined by
$$\hat\Delta_0=\max\,\{\Delta(n_1,\ldots,n_k):(n_1,\ldots,n_k)\in {\mathcal S}(n)\}.$$

\subsection{Ideal immersions}
Inequality \e{6.4} enables us to introduce the  notion of ideal immersions.

\begin{definition} An isometric immersion  of a Riemannian $n$-manifold $M$ in $\mathbb E^m$ is called an {\it ideal immersion\/} if it satisfies the equality case of \eqref{6.4} identically. \end{definition}

 The above maximum principle yields the following important fact:
\vskip.1in

\begin{theorem}  {\rm \cite{c12,c16}} If an isometric immersion $x:M\to \mathbb E^m$ of a Riemannian $n$-manifold into $\mathbb E^m$  satisfies equality \e{6.4} for a given $k$-tuple $(n_1,\ldots,n_k)\in {\mathcal S}(n)$, then it is an ideal immersion automatically. \end{theorem}

\begin{remark} ({\it Physical Interpretation of Ideal Immersions}\/) An isometric immersion  $x:M\to  E^m$ is an {\it ideal immersion} means that  $M$  receives the least possible amount of tension  (given by $\hat\Delta_0(p)$) from the surrounding space at each point $p$ on $M$. This is due to \e{6.4} and the well-known fact that the mean curvature vector field is exactly the tension field for an isometric immersion of a Riemannian manifold in another Riemannian manifold; thus  the squared mean curvature at each point on the submanifold simply measures the amount of tension the submanifold is receiving from the surrounding space at that point. 

For this reason, an ideal immersion is also called a {\it best way of living.} \end{remark}

Although a standard $n$-sphere $S^n$ does admit an ideal immersion in $\mathbb E^{n+1}$,  the following results show that  other compact rank one symmetric spaces do not admit ideal immersions in any Euclidean space.

\begin{proposition}  Let $FP^n\, (n>1)$ denote a projective space over real, complex, or quaternion field equipped with a standard  Riemannian metric, where the real dimension of $FP^n$ is $n,2n$ or $4n$, according to
$F=\mbox{\bf R},\mbox{\bf C}$ or $\mbox{\bf H}$. Then $FP^n$  doesn't admit ideal immersions in a Euclidean
space, regardless of codimension. \end{proposition}

\begin{proposition} The Cayley plane  ${\mathcal O}P^2$ with a standard Riemannian metric does not admit ideal immersions in a Euclidean space, regardless of codimension. \end{proposition}

\subsection{Two world problems}

A goal of  research in the direction is to solve the following.

\begin{WP1} Determine those individuals (those Riemannian manifolds)  who admit a best way of living (ideal immersions) in a best world  (in a space form).\end{WP1}

\begin{WP2} Determine the  best ways of living for those individuals   who admit an ideal immersion in a best world.\end{WP2}

Here, by a ``best world'' we mean {\it a space with the highest degree of  homogeneity}. According to work  of Lie, Klein and Killing,
the family of  Riemannian manifolds  with the highest degree of homogeneity consists of Euclidean spaces, Riemannian spheres, real
projective spaces, and real hyperbolic spaces. 

Such spaces have the  highest degree of homogeneity simply because they have the largest groups of isometries. Hence, a best world in the terminology of differential geometry is nothing but a  Riemannian space form.

\subsection{Ideal immersions as stable critical points of total tension  functional} Ideal immersions are closely related with critical point problem in the theory of total mean curvature. This can be seen as follows:

Let $M$ be a compact Riemannian manifold (with or without boundary). Denote by ${\mathcal I}(M,R^m(\epsilon))$
the family of isometric immersions of $M$ into a real space form $R^m(\epsilon)$. For each $\phi\in {\mathcal I}(M,R^m(\epsilon))$, we define its total tension (or the total squared mean curvature, or Willmore)
functional by the formula
$$T(\phi)=\int_M H_\phi ^2{\rm d}V,$$
where $H^2_\phi$ denotes the squared mean curvature of $\phi:M\to R^m(\epsilon)$. 

It follows from  Theorem \ref{T:5.1}  that an ideal immersion of $M$ into $R^m(\epsilon)$ is a critical point of the total tension functional within the class of  ${\mathcal I}(M,R^m(\epsilon))$  automatically. Clearly, every
ideal immersion of $M$ in $R^m(\epsilon)$ is also stable, i.e., the second variation of $T(\phi)$ is nonnegative for each variation of $\phi$ in the class of ${\mathcal I}(M,R^m(\epsilon))$.

\subsection{Size of the smallest ball containing an ideal submanifold}

According to Nash's embedding theorem, every compact Riemannian $n$-manifold can be isometrically embedded in
any small portion of Euclidean space if the codimension is large enough. In contrast the following theorem states that  ideal compact submanifolds cannot be contained in a very small ball of the Euclidean space. In fact, by applying Theorem \ref{T:5.1}, we can estimate the radius of the smallest ball  in the Euclidean space which contains a given compact ideal submanifold in terms of $\delta$-invariants. 

\begin{theorem} {\rm \cite{c12}} Let $\phi :M\to \Bbb E^m$ be an ideal immersion from a compact Riemannian $n$-manifold $M$ into a  Euclidean $m$-space.  Then, regardless of codimension,  the radius $R$ of  the  smallest ball $B(R)$ containing $\phi(M)$ satisfies 
$$R^2\geq {{{\rm vol}(M)}\over{\int_M\hat\Delta_0
{\rm d}V}},$$ 
with the equality sign  holding if and only if $x$ is a $1$-type ideal immersion. \end{theorem}

\section[Applications to estimates of  eigenvalues of $\Delta$]{Applications of $\delta$-invariants to estimates of eigenvalues of the Laplacian}

\subsection{Type number of immersions}
For an isometric immersion $x:M\to \mathbb E^m$ of $M$ in $\mathbb E^m$, let 
$$x=x_0+\sum_{t=p}^q x_t,\quad \Delta x_t=\lambda_t x_t$$ denote the spectral resolution of $x$, where $x_0$ is center of
mass of $M$ in $\mathbb E^m$ and $\Delta$ is the Laplacian of $M$. 
The set $$T(x)=\{t\in \hbox{\bf Z}:x_t\ne \hbox{constant map}\}$$ is called the order of the submanifold. The smallest element $p$ in $T(x)$ is called the {\it lower order} of $x$ and the supremum $q$ of $T(x)$ is called the {\it upper order} of $x$. The immersion is said to be {\it of finite type\/} if the upper order $q$ is finite; and it is said to be of infinite type if the upper order $q$ is infinite. Moreover, the immersion is said to be {\it of $k$-type} if $T(x)$ contains  exactly $k$ elements. 

Clearly, the immersion is of 1-type if and only if $p=q$. In this case, the immersion is called a 1-{\it type immersion  of order}
$\{p\}$ (see \cite{c9} for a comprehensive survey on submanifolds of finite type). 

\subsection{$\lambda_1$ of compact homogeneous spaces and $\delta$-invariants}
By applying the inequalities \eqref{5.2} and the theory of finite  type submanifolds \cite{c2,c6}, we can establish the following new intrinsic results concerning intrinsic spectral properties of homogeneous spaces  via  extrinsic data.  

\begin{theorem} \label{T:7.1} If $M$ is a compact homogeneous Riemannian $n$-mani\-fold with irreducible isotropy action, then the first nonzero eigenvalue $\lambda_1$ of the Laplacian on $M$ satisfies
\begin{align}\label{7.1} \lambda_1\geq n\,\Delta(n_1,\ldots,n_k)\end{align}
for any $k$-tuple $(n_1,\ldots,n_k)\in \mathcal S(n)$.

The equality sign of \eqref{7.1} holds if and only if $M$ admits a 1-type ideal immersion in a Euclidean space.  \end{theorem}

\begin{remark} If $k=0$, inequality \eqref{7.1}  reduces to  the well-known result of T. Nagano on $\lambda_1$ obtained in \cite{nagano};  namely
\begin{align}\lambda_1\geq n\rho,\end{align} where $\rho=
\tau/{n\choose 2}$ is the normalized scalar curvature.  
In general, we have  $$\Delta(n_1,\ldots,n_k)\geq\rho.$$ Moreover, we have   $\Delta(n_1,\ldots,n_k)>\rho$ for  $k>0$ on most Riemannian manifolds. \end{remark}

For $\delta$-invariants on a  compact homogeneous space, we also have the following.

\begin{theorem} \label{T:7.2}  The following statements hold.
\vskip.04in

$(1)$ A compact homogeneous Riemannian $n$-manifold $M$ with irreducible isotropy action admits an ideal immersion in a Euclidean space if and only if it satisfies $\lambda_1=n\,\hat\Delta_0$.
\vskip.04in

$(2)$ Every ideal immersion of a compact homogeneous Riemannian manifold with irreducible isotropy action  in a
Euclidean space is a 1-type immersion of order $\{1\}$. 
\vskip.04in

$(3)$ If a compact homogeneous  Riemannian $n$-manifold with irreducible isotr\-opy action admits an ideal
immersion in a Euclidean space, then  $$\hat\Delta_0=\Delta(n_1,\ldots,n_1)$$ for some $(n_1,\ldots,n_1)$ $\in {\mathcal S}(n,k)$ with $\, kn_1=n$.
\vskip.04in

$(4) \ $If a compact homogeneous  Riemannian $n$-manifold $\,M\, $ with irreducible isotr\-opy action admits an ideal
immersion in a Euclidean $m$-space such that the image of $M$ is contained in a  hypersphere with radius $r$, then we have $$\lambda_1=\frac{n}{r^2},\quad \hat\Delta_0=\frac{1}{r^2}.$$\end{theorem}

\begin{remark}  For a  compact irreducible homogeneous $n$-manifold $M$, Theorem \ref{T:7.2} can be applied to determine whether $M$ admits an ideal immersion. In principle,  $\lambda_1$ (using Freudenthal's
formula for Casimir's operator) and the invariant $\hat\Delta$ are both ``computable'' for every compact irreducible homogeneous Riemannian manifold. 

For many compact irreducible symmetric spaces $M=G/H$ with $G$ being a classical group, $\lambda_1$ of $M$ has been computed by T. Nagano in \cite{nagano}. \end{remark}

\begin{remark}  Besides Riemannian spheres, there do exist other compact homogeneous Riemannian manifolds which admit ideal immersions in the Euclidean space. For instance, the following three compact homogeneous Riemannian manifolds:
$$SU(3)/T^2,\quad Sp(3)/Sp(1)^3,\quad\mbox{and}\quad F_4/\mbox{Spin}(8)$$ admit ideal immersions  in $\mathbb  E^8,$ $ \mathbb E^{14},$ and $\mathbb E^{26}$ of codimension 2 associated with $$(3,3)\in{\mathcal S}(6),\;\;(3,3,3,3)\in {\mathcal S}(12),\;\; {\rm and}\;\;   (12,12)\in{\mathcal S}(24),$$ respectively. 

These ideal immersions of $SU(3)/T^2,\,Sp(3)/Sp(1)^3,$ and $\,F_4/\mbox{Spin}(8)$ in $\mathbb E^8,\,\mathbb E^{14}$ and $\mathbb E^{26}$ arise from their isometric immersions in $S^7,S^{13}$ and $S^{25}$ respectively as minimal isoparametric hypersurfaces. \end{remark}

\subsection{Estimate of eigenvalues of $\Delta$ and $\delta$-invariants}

For $\delta$-invariants on a general compact Riemannian manifold, we have the following general results.

\begin{theorem} \label{T:7.3}   Let $\phi:M \to \mathbb E^m$ be an isometric immersion from a compact Riemannian $n$-manifold into a Euclidean $m$-space. Then 
$$\lambda_q\geq {n\over{{\rm vol}(M)}} \int_M\Delta(n_1, \ldots,n_k){\rm d}V$$ for each $k$-tuple $(n_1,\ldots,n_k)\in {\mathcal S}(n)$, where $q$ is the upper order of the  immersion.

The equality case  holds for a $k$-tuple $(n_1,\ldots,n_k)\in {\mathcal S}(n)$ if and only if $\phi$ is a $1$-type ideal immersion of order $\{q\}$ associated with $(n_1,\ldots,n_k)$, i.e., $\phi$ is a $1$-type immersion of order $\{q\}$  satisfying
$$\Delta(n_1,\ldots,n_k)= H^2$$ identically.  \end{theorem}

\begin{theorem} \label{T:7.4} Let $M$ be a compact Riemannian $n$-manifold. We have

$(1)$ If $M$ admits an ideal immersion in a Euclidean space  associated with a $k$-tuple
$(n_1,\ldots,n_k)$, then \begin{align}\label{7.3}\lambda_1\leq  {n\over{{\rm vol}(M)}}
\int_M\Delta (n_1, \ldots,n_k){\rm d}V.\end{align}

$(2)$ If $M$ satisfies 
$$\lambda_p\leq  {n\over{{\rm vol}(M)}}\int_M\hat\Delta_0 {\rm d}V,$$ then every order $\{p\}$, $1$-type isometric immersion of $M$ in
a Euclidean space is an ideal immersion.

$(3)$ An ideal immersion of $M$ satisfies the equality case of \e{7.3} if and only if the immersion is a $1$-type ideal immersion of order $\{1\}$.\end{theorem}

\begin{theorem} \label{T:7.5} Let $M$ be a compact Riemannian manifold which admits a 1-type isometric immersion of order $\{1\}$ in a
Euclidean space. If $M$ admits an ideal immersion into some Euclidean space, then
\begin{align}\label{7.4}\lambda_1 =\frac{n}{{\rm vol}(M)} \int_M\hat\Delta_0{\rm d}V.\end{align}
In particular, if a compact strongly harmonic $n$-manifold admits an ideal immersion in a Euclidean space, then \e{7.4} holds.\end{theorem}
\vskip.1in
By applying Theorem \ref{T:7.5} we obtain the following simple  obstruction to ideal immersions for compact Riemannian manifolds in terms of the $\delta$-invariant $\hat\Delta_0$.

\begin{theorem} \label{T:7.6}  Let $M$ be a compact Riemannian $n$-manifold. If $M$ satisfies
\begin{align}\label{7.5}\lambda_1>{n\over{{\rm vol}(M)}} \int_M\hat\Delta_0{\rm d}V,\end{align}
then $M$ admits no ideal immersion in a Euclidean space for any codimension.

In particular, every compact Riemannian manifold with non-positive sectional curvatures admits no ideal immersion in a Euclidean space for any codimension.\end{theorem}

For ideal immersions we have the following relationship between the $\delta$-invariant and the first two eigenvalues of the Laplacian $\Delta$ on compact Riemannian manifolds.

\begin{theorem} Let $x:M\to\mathbb E^m$ be  an ideal immersion of a compact Riemannian $n$-manifold in Euclidean $m$-space $\mathbb E^m$. Then
\begin{align}\label{7.6}\int_M \hat\Delta_0{\rm d}V\geq {1\over{n^2}}\left\{n(\lambda_1+\lambda_2)-R^2\lambda_1 \lambda_2\right\}{\rm vol}(M),\end{align}
where $R$ denotes the radius of the smallest ball $B(R)$  in $\mathbb E^m$ which contains the image of $M$.

 The equality sign of \e{7.6} holds if and only if the image of $M$ is contained in the boundary $S^{m-1}$ of the ball $B(R)$ and the
immersion $x$ is a $1$-type ideal immersion of order $\{1\}$, or a $1$-type ideal immersion of order $\{2\}$, or a $2$-type ideal immersion of order 
$\{1,2\}$.

Moreover, if the equality case of \e{7.6} holds, then $M$ is mass-symmetric in $S^{m-1}$ and $\hat\Delta_0$ is a constant on $M$. \end{theorem}

For further applications of $\delta$-invariants to eigenvalues of Laplacian on Riemannian manifolds, see \cite{c12,c16,c18}.

\section[Applications  to minimal  immersions]{Applications of $\delta$-invariants to minimal  immersions}

Since the fundamental inequalities provide us the simplest relationship between the $\delta$-invariants and the squared  mean curvature, 
one important immediate application of the $\delta$-invariants and the fundamental inequalities is to provide many solutions to Problem 1 proposed by S. S. Chern.

\begin{theorem} \label{T:8.1} Let $M$ be a Riemannian $n$-manifold. If there exists  a point $p\in M$  and a $k$-tuple $(n_1,\ldots,n_k)\in {\mathcal S}(n)$  such that 
$$\delta{(n_1,\ldots,n_k)}(p)>\frac 1 2{{n(n-1)}}\epsilon-  \frac 1 2\sum {{n_j(n_j-1)}}\epsilon,$$ then $M$ admits no minimal isometric immersion into a Riemannian space form $R^m(\epsilon)$, regardless of codimension. \end{theorem}

In particular,  we have the following solution to Problem 1.

\begin{theorem} \label{T:8.2} Let $M$ be a Riemannian $n$-manifold. If there exists  a point $p\in M$  and a $k$-tuple $(n_1,\ldots,n_k)\in {\mathcal S}(n)$  such that 
 $$\delta{(n_1,\ldots,n_k)}(p)>0,$$ then $M$ admits no minimal isometric immersion into Euclidean space for any  codimension.\end{theorem}

\begin{remark} There exist many Riemannian manifolds with scalar $\tau<0$ but with some positive $\delta$-invariants (see, for instance \cite{su2}).\end{remark}

\begin{remark} It is important to mention that, for each integer $n\geq 2$ and each  $k$-tuple $(n_1,\ldots,n_k)\in\mathcal S(n)$, the condition on $\de{(n_1,\ldots,n_k)}$ given in Theorem \ref{T:8.2} is sharp. This
can be seen as follows:  

Let $f_j:M^{n_j}_j\to \mathbb E^{m_j},\, j=1,\ldots, k,$ be $k$ minimal submanifolds and  $\iota$ a totally geodesic immersion of a Euclidean $(n-\sum n_j)$-space  into a Euclidean space. Then the  invariant $\de{(n_1,\ldots,n_k)}$ of the Riemannian product $M^{n_1}_1\times\cdots \times M^{n_k}_k\times\mathbb  E^{n-\sum n_j}$  vanishes
identically.  Clearly, the product immersion $f_1\times\cdots f_k\times \iota\;$ is a minimal  immersion. \end{remark}

\begin{corollary} Let $M^{n_1}_1,\ldots,M^{n_k}_k$ be  Riemannian manifolds of dimensions $\geq 2$ with scalar
curvatures $\leq 0$. Then every minimal isometric immersion
$$f:M^{n_1}_1\times\cdots\times M^{n_k}_k\times  \mathbb E^{n- \sum n_j}\to\mathbb E^m$$ of  $M^{n_1}_1\times\cdots\times M^{n_k}_k\times  \mathbb E^{n-\sum n_j}$ in any Euclidean $m$-space is a product immersion $f_1\times\cdots\times f_{k}\times \iota$ of $k$  minimal immersions $f_j:M^{n_j}_j\to \mathbb E^{m_j},\, j=1,\ldots,k,$ and a totally geodesic immersion
$\iota$. \end{corollary} 

For a submanifold $M$ in a real space form $R^m(\epsilon)$, if we consider $\delta(2)$ on $M$, then inequality \e{5.1} reduces to  \begin{align}\label{8.1} \delta(2)\leq {{n^2(n-2)}\over {2(n-1)}}H^2+{1\over 2}(n+1)(n-2)\epsilon.\end{align} This inequality implies that if $M$ admits an isometric minimal immersion into some Euclidean space, one has  \begin{align}K\geq \tau.\end{align}

\section[Applications to Lagrangian and slant immersions]{Applications of $\delta$-invariants to Lagrangian and slant immersions}

An  immersion of a Riemannian $n$-manifold $M$ in a Hermitian manifold $\tilde M$ is called {\it totally real\/} (or {\it isotropic\/} in symplectic geometry) if the almost complex structure $J$ of $\tilde M$ maps each tangent space of $M$ into its corresponding normal space. In particular, a totally real immersion is said to be {\it Lagrangian\/} if $\dim M=\dim_{\hbox{\bf C}}\tilde M$.   

For Lagrangian immersions in  complex Euclidean $n$-space $\hbox{\bf C}^n$, a result of Gromov states that a compact $n$-manifold $M$ admits a Lagrangian immersion (not necessarily isometric) into $\hbox{\bf C}^n$ if and only if the complexification $TM\otimes \hbox{\bf C}$ of the tangent bundle of $M$ is trivial. Gromov's result \cite{Gro} implies that there is  no topological obstruction to Lagrangian immersions for compact 3-manifolds in $\hbox{\bf C}^3$, because the tangent bundle of a 3-manifold is always trivial.

The class of Lagrangian immersions is included in the class of slant immersions which are defined as follows:

Let $M$ be a Riemannian manifold isometrically immersed in almost Hermitian manifold $\tilde M$ with almost complex structure $J$ and almost Hermitian metric $g$.   For any nonzero vector $X$ tangent to $M$ at a point $p\in M$, the
angle $\theta (X)$  between $JX$ and the tangent space $T_pM$ is called the  Wirtinger angle of $X$. 

A submanifold $M$ of $\tilde M$ is called {\it slant} if the Wirtinger angle $\theta (X)$ is a constant (which is independent of the choice of $x \in M$ and of $X \in T_{x}N$). The Wirtinger angle of a slant submanifold is called the  slant angle of the slant submanifold.

Complex submanifolds  and totally real submanifolds are nothing but slant submanifolds with slant angle $\theta = 0$ and $\theta = \pi /2$, respectively. A slant submanifold is called proper if it is neither complex nor totally real. 
There exist abundant examples of slant submanifolds in ${\bf C}^n$ (cf. \cite{cbook3}).
 
The author obtained in \cite{cbook3} a topological obstruction to slant immersions, namely: 

A compact $2k$-manifold $M$ with
$H^{2i}(M;\hbox{\bf R})=0$ for some $1\leq i\leq k$ admits no slant immersion in any K\"ahlerian manifold $\wt M^m$ unless it is totally real (or Lagrangian when $m=2k$). Moreover, Chen and Y. Tazawa proved in \cite{CT2} that there exist no slant immersions of a compact $n$-manifold in $\hbox{\bf C}^m$ unless it is totally real. On the other hand, there do exist compact slant submanifolds in a complex $n$-torus.

From the Riemannian point of view, it is natural to ask the following basic question.
\vskip.1in

\noindent {\bf Problem 2} {\it What are necessary conditions for a compact Riemannian manifold to admit a Lagrangian $($or more generally, slant\/$)$  immersion in $\hbox{\bf C}^n$? }
\vskip.1in
 
Another important application of the $\delta$-invariants is the following solution to Problem 2.

\begin{theorem} \label{T:9.1}  Let $M$ be a compact Riemannian $n$-manifold with null first
Betti number $b_1$ or with finite fundamental group $\pi_1$. If there is a $k$-tuple
$(n_1,\ldots,n_k)\in {\mathcal S}(n)$  such that 
\begin{align}\label{9.1}\delta{(n_1,\ldots,n_k)}>0,\end{align} then $M$
admits no  slant isometric immersion in a complex $n$-torus $ CT^n$ or in the complex Euclidean $n$-space $\mbox{\bf C}^n$. 

In particular, if \e{9.1} holds for some $(n_1,\ldots,n_k)\in {\mathcal S}(n)$, then $M$ admits no  Lagrangian isometric immersion in a complex n-torus or
in complex Euclidean $n$-space. \end{theorem}

\begin{remark}  The assumptions on the finiteness of $\pi_1(M)$ and vanishing of
$b_1(M)$ given in Theorem \ref{T:9.1} are both necessary for $n\geq 3$. This can be seen from the following example:

Let $F:S^1\to\hbox{\bf C}$ be the unit circle in the complex plane given by $F(s)=\rme^{\i s}$ and let $\iota:S^{n-1}\to E^n$ $(n\geq 3)$ be the unit hypersphere in $E^n$ centered at the origin. Denote by $f:S^1\times S^{n-1}\to {\bf C}^n\,$ the complex extensor defined by $$f(s,p)=F(s)\otimes \iota(p),\;\; p\in S^{n-1}.$$Then $f$ is an isometric Lagrangian immersion of $M=:S^1\times S^{n-1}$ into $C^n$ which carries each pair
$\{(u,p),\,(-u,-p)\}$  of points in $S^1\times S^{n-1}$ to a point in ${\bf C}^n$. Clearly, we have
$\pi_1(M)=\hbox{\bf Z}$ and $b_1(M)=1$. Moreover, for each $k$-tuple
$(n_1,\ldots,n_k)\in {\mathcal S}(n)$, the $\delta$-invariant  $\delta(n_1,\ldots,n_k)$ on $M$ is a positive constant. 

This example shows that both the conditions on $\pi_1(M)$ and $b_1(M)$ cannot be removed.
\end{remark}

For Lagrangian immersions in complex space forms,  we also have the following two results.

\begin{theorem}   Let $M$ be a compact  Riemannian $n$-manifold  with finite
fundamental group or with null first Betti number. If 
$$\delta{(n_1,\ldots,n_k)}>\frac 1 2\left({{n(n-1)}}-\sum_{j=1}^k {{n_j(n_j-1)}}\right) $$ holds  for
some  $k$-tuple $(n_1,\ldots,n_k)\in  {\mathcal S}(n)$, then $M$ admits no Lagrangian  isometric immersion  into  the complex projective $n$-space $ CP^n(4)$. \end{theorem}

\begin{theorem}   Let $M$ be a compact Riemannian $n$-manifold either with finite
fundamental group or with null first Betti number. If 
$$\delta{(n_1,\ldots,n_k)}>\frac 1 2 \left(\sum_{j=1}^k {{n_j(n_j-1)}}-{{n(n-1)}}\right) $$ holds  for
some  $k$-tuple $(n_1,\ldots,n_k)\in {\mathcal S}(n)$, then $M$ admits no Lagrangian isometric immersion into the complex hyperbolic $n$-space $ CH^n(-4)$.\end{theorem}

\section{Applications to rigidity problems}

Although the ideal immersions of a given Riemannian manifold in a Euclidean space is not necessarily unique, very often the $\delta$-invariants and Theorem \ref{T:5.1} can  be applied to obtain the rigidity for isometric immersion of arbitrary codimension from a
Riemannian manifold into a Riemannian space form; in particular,  to obtain a rigidity theorem for open portions of a homogeneous Riemannian manifold isometrically immersed in a Euclidean space, regardless of codimension. 

The philosophy of the rigidity comes from the fact that, for a given Riemannian manifold $M$, inequality \e{5.1} provides us a lower bound of the squared mean curvature. When the inequality is actually an equality, the submanifold is an ideal submanifold according to our maximum principle. In this case Theorem \ref{T:5.1} shows that the shape operators of the submanifold must take a special simple form. In many cases, this information on the Riemannian structure of $M$ and on the shape operators is sufficient to conclude the rigidity of the submanifold {\it without any global assumption and regardless of codimension}.

Here we provide three of many such applications.

\begin{proposition} Let $M$ be an open portion of a unit $n$-sphere $S^n(1)$. Then, for any isometric immersion of $M$  into $\mathbb E^m$, we have
\begin{align}\label{10.1} H^2\geq 1\end{align}
regardless of codimension.

The equality case of  \e{10.1} holds identically if and only if $M$ is immersed as an open portion of  an ordinary hypersphere  in an affine $(n+1)$-subspace $\mathbb E^{n+1}$ of $\mathbb E^m$.\end{proposition}

\begin{proposition}   Let $M$ be an open portion of  $S^k(1)\times  \mathbb E^{n-k}, k>1$. Then, for any isometric immersion of $M$  into  $\mathbb E^m$, we have 
\begin{align}\label{10.2} H^2\geq \left(\frac k n\right)^2\end{align} regardless of codimension.

The equality case of \e{10.2} holds identically if and only if  $M$ is immersed as an open portion of  an ordinary spherical hypercylinder  in an affine $(n+1)$-subspace $\mathbb E^{n+1}$ of $ \mathbb E^m$.\end{proposition}

\begin{proposition}  Let $M$ be an open portion of  $S^k(1)\times S^{n-k}(1), 1<k<n-1$. Then for any isometric immersion of $M$  into 
$\mathbb E^m$, we have 
\begin{align}\label{10.3} H^2\geq \left(\frac k n\right)^2+\left(\frac {n-k} n\right)^2\end{align} regardless of codimension.

The equality case of \e{10.3} holds identically if and only if  $M$ is embedded in a standard way in an affine $(n+2)$-subspace  of $\mathbb E^m$.\end{proposition}
\vskip.1in
Theorem \ref{T:5.1} and Moore's lemma can be applied to provide some decomposition results. For instance, we have

\begin{proposition} \label{P:9.4} Let \ $M^{n_1}_1,\ldots,M^{n_k}_k$ $(k\geq 2)$ be $k$ Riemannian manifolds satisfying $n_1+\cdots+n_k\leq n$. Then, regardless of codimension, we have

{\rm (1)} Every isometric immersion of $M^{n_1}_1\times\cdots\times M^{n_k}_k\times H^{n-\sum n_j}(-\epsilon)$ into the hyperbolic $m$-space $H^m(-\epsilon),\, \epsilon>0,$ satisfies $$H^2\geq  {{b(n_1,\ldots,n_k)}\over {c(n_1,\ldots,n_k)}}\epsilon.$$

In particular, if $$H^2= \frac {b(n_1,\ldots,n_k)}{c(n_1,\ldots,n_k)}\epsilon$$
identically, then the immersion is a product immersion.

{\rm (2)} Every minimal isometric immersion from $M^{n_1}_1\times\cdots\times M^{n_k}_k\times  \mathbb E^{n-\sum n_j}$ into a Euclidean space is a product immersion.\end{proposition}

Statement (2) of Proposition \ref{P:9.4} was due to N. Ejiri \cite{Ej}.

\section{Applications to warped products}

 Let $N_1,\cdots ,N_k$ be Riemannian manifolds and let $N=N_1 \times\cdots \times N_ k$ be the Cartesian product of $N_1,\ldots,N_k$. For each $i$, denote by $\pi_i : N \to N_i$  the canonical projection of $N$ onto $N_i$.  When there is no confusion, we identify $N_i$ with a horizontal lift of $N_i$ in $N$ via $\pi_i$.

If $f_2,\cdots ,f_k : N_1\to \hbox{\bf R}_+$ are positive-valued functions, then \begin{align*}
 \left< X,Y \right>:= \left<\pi_{1*}X,\pi_{1*}Y\right> + \sum^k_{i=2} (f_i \circ \pi_1)^2 \left< \pi\sb {i*}X, \pi\sb {i*}Y \right>\end{align*}
defines a Riemannian metric $g$ on $N$, called a multiply warped product metric. The product manifold $N$ endowed with this metric is denoted by $$N_1 \times_{f_2} N_2 \times \cdots \times_{f_k}N_k.$$

For a multiply warped product manifold $N_1 \times_{f_2} N_2 \times \cdots \times_{f_k}N_k$,  let $\mathcal D_i$  denote the distributions obtained from the vectors tangent to  $N_i$ (or more
precisely, vectors tangent to the horizontal lifts of $N_i$). Assume that  $$\phi:N_1 \times_{f_2} N_2 \times \cdots \times_{f_k}N_k\to \tilde M$$ is an isometric immersion of a multiply warped product $N_1 \times_{f_2} N_2 \times \cdots \times_{f_k}N_k$ into a Riemannian manifold $\tilde M$. Denote by $h$ the second fundamental form of $\phi$. Then the immersion $\phi$ is called {\it mixed totally geodesic\/} if $$h(\mathcal D_i,\mathcal D_j)=\{0\}$$ holds for  distinct $i,j\in \{1,\ldots,k\}$.

Let $$\phi:N_1 \times_{f_2} N_2 \times \cdots \times_{f_k} N_k\to \tilde M$$ be an isometric immersion of a multiply warped product $N_1 \times_{f_2} N_2 \times \cdots \times_{f_k} N_k$ into an arbitrary Riemannian manifold $\tilde M$.   Denote by  ${\rm trace}\,h_i$   the trace of $h$ restricted to $N_i$,  that is \begin{align*}{\rm trace}\,h_i=\sum_{\alpha=1}^{n_i} h(e_\alpha,e_\alpha)\end{align*} for some orthonormal frame fields $e_1,\ldots,e_{n_i}$  of $\mathcal D_i$.

By considering a special $\delta$-invariant, we have the following  general optimal result for any isometric immersion of a warped into a real space form for any codimension:

\begin{theorem} {\rm \cite{c26}} \label{T:11.1} Let $\phi:N_1\times_f N_2\to R^m(\epsilon)$ be an isometric immersion of a warped product into a  Riemannian $m$-manifold of constant sectional curvature $\epsilon$. Then we have
\begin{align} \label{r1.2}\frac{\Delta f} {f}\leq {\frac{(n_1+n_2)^2}{4n_2}}H^2+ n_1\epsilon,\end{align} where $n_i=\dim N_i$, $i=1,2$,  $H^2$ is the squared mean curvature of $\phi$, and $\Delta$ is the Laplacian operator of $N_1$.

 The equality sign of \eqref{r1.2} holds identically if and only if $\phi:N_1\times_f N_2\to R^m(\epsilon)$ is a mixed totally geodesic immersion satisfying ${\rm trace}\, h_1={\rm trace }\,h_2$.\end{theorem}

This result was extended to multiply warped products as follows.

\begin{theorem}\label{T:11.2} {\rm \cite{cd,cw}}
Let $\phi:N_1 \times_{f_2} N_2 \times \cdots \times_{f_k} N_k\to
\tilde M^m$ be an isometric immersion of a multiply warped product
$N:=N_1 \times_{f_2} N_2 \times \cdots \times_{f_k} N_k$ into an
arbitrary Riemannian $m$-manifold. Then we have
\begin{align} \label{11.2} & \sum_{j=2}^k n_j\frac{\Delta f_j}{f_j} \leq \frac{n^2}{4} H^2
 + n_1(n-n_1) \max\tilde K,\quad n=\sum_{j=1}^kn_j, \end{align}
where $\max\tilde K(p)$ denotes the maximum of the sectional
curvature function of $\tilde M^m$ restricted to $2$-plane sections of
the tangent space $T_pN$ of $N$ at $p=(p_1,\ldots,p_k)$.

The equality sign of \eqref{11.2} holds identically if and only if the following two statements hold:
\vskip.04in

 {\rm (1)} $\phi$ is a mixed totally geodesic immersion satisfying ${\rm trace}\, h_1=\cdots={\rm trace }\,h_k$;
\vskip.04in

 {\rm (2)} At each point $p\in N$,  the sectional curvature function $\tilde K$
of $\tilde M^m$ satisfies $\tilde K(u,v)=\max \tilde K(p)$ for each
unit vector $u$ in $T_{p_1}(N_1)$ and each unit vector $v$ in $
T_{(p_2,\cdots,p_k)}(N_2\times \cdots\times N_k)$.
 \end{theorem}
 
 The following example shows that inequality \eqref{11.2} is sharp.
 
\begin{example} Let  $M_1\times_{\rho_2} M_2\times \cdots\times_{\rho_k} M_k$ be a {\it multiply warped product representation\/} of a real space form $R^m(\epsilon)$. 
Assume that $\psi_1:N_1\to M_1$ is a minimal immersion of $N_1$ into  $M_1$ and let $f_2,\ldots,f_k$ be the restriction of $\rho_2,\ldots,\rho_k$ on $N_1$. 
Then the following warped product immersion:
$$\psi=(\psi_1,id,\ldots,id):N_1\times_{f_2} M_2\times \cdots\times_{f_k} M_k\to  M_1\times_{\rho_2} M_2\times \cdots\times_{\rho_k} M_k\subset R^m(\epsilon)$$ is a mixed totally geodesic warped product submanifold of $R^m(\epsilon)$ which satisfies the condition:
 \begin{align}&{\rm trace}\, h_1=\cdots={\rm trace }\,h_k =0. \end{align} 
Thus, the immersion $\psi$ satisfies the equality case of \eqref{11.2} according to Theorem \ref{T:11.2}. Therefore, inequality \eqref{11.2} is optimal.

\end{example}

As an immediate consequence of Theorem \ref{T:11.2}, we have the following.

\begin{corollary}\label{C:11.1}  {\rm \cite{cd,cw}}
Let $\phi:N_1 \times_{f_2} N_2 \times \cdots \times_{f_k} N_k\to
R^m(\epsilon)$ be an isometric immersion of a multiply warped product $N_1
\times_{f_2} N_2 \times \cdots \times_{f_k} N_k$ into a Riemannian
$m$-manifold $R^m(\epsilon)$ of constant curvature $\epsilon$. Then we have
 \begin{align} \label{11.3} & \sum_{j=2}^k n_j\frac{\Delta f_j}{f_j} \leq \frac{n^2}{4} H^2
 + n_1(n-n_1)\epsilon,\quad n=\sum_{j=1}^kn_j. \end{align}

  The equality sign of \eqref{11.3} holds identically if and only if $\phi$ is a mixed totally geodesic
immersion satisfying ${\rm trace}\, h_1=\cdots={\rm trace }\,h_k$.
  \end{corollary}

Combining Theorem \ref{T:11.1} and N\"olker's  theorem gives immediately the following.

\begin{corollary}\label{C:11:2} Let $\phi:N_1 \times_{f_2} N_2 \times \cdots \times_{f_k} N_k\to
R^m(\epsilon)$ be an isometric immersion of the multiply warped product $N_1 \times_{f_2} N_2 \times \cdots \times_{f_k} N_k$ into a Riemannian $m$-manifold of constant curvature $\epsilon$.  If we have
 \begin{align*}  & \sum_{j=2}^k n_j\frac{\Delta f_j}{f_j} = \frac{n^2}{4} H^2 + n_1(n-n_1) \epsilon, \end{align*}  then $\phi$ is a warped product immersion.
 \end{corollary}

By applying Theorem \ref{T:11.1} we have  the following.

\begin{corollary} \label{C:11.3} If $N_1$ is a compact Riemannian manifold, then every warped product $N_1\times_f N_2$ does not admit an isometric minimal immersion into a Euclidean space or a hyperbolic space for any codimension.
\end{corollary}

\begin{corollary}\label{C:11.4} If $f_2,\ldots,f_k$ are  harmonic functions on $N_1$ or
eigenfunctions of the Laplacian $\Delta$ on $N_1$ with positive eigenvalues, then the multiply warped product manifold $N_1 \times_{f_2} N_2 \times \cdots \times_{f_k} N_k$  cannot be
isometrically immersed into any Riemannian manifold of negative sectional curvature as a minimal submanifold.

 \end{corollary}

\begin{corollary}\label{C:11.5} If $f_2,\ldots,f_k$ are  eigenfunctions of the Laplacian $\Delta$ on
$N_1$ with nonnegative eigenvalues and at least one of $f_2,\ldots, f_k$ is non-harmonic, then the multiply warped product manifold $N_1 \times_{f_2} N_2 \times \cdots \times_{f_k} N_k$ cannot be isometrically immersed into any Riemannian manifold of non-positive sectional curvature as a minimal submanifold. \end{corollary}

By applying Theorem \ref{T:11.1} we also have the following.

\begin{corollary}\label{C:11.6} If $f_2,\ldots,f_k$ are harmonic functions on $N_1$, then every
isometric minimal immersion of the multiply warped product manifold $N_1 \times_{f_2} N_2 \times \cdots \times_{f_k} N_k$ into a Euclidean space is a warped product immersion. \end{corollary}

Since the proof of Theorem \ref{T:11.1} bases only on the equation of Gauss, the same proof as Theorem \ref{T:11.1} yields the following.

\begin{corollary}\label{C:11.7}
Let $\phi:N_1 \times_{f_2} N_2 \times \cdots \times_{f_k} N_k\to \tilde M^m(4\epsilon)$ be a totally real isometric immersion of the multiply warped product $N_1 \times_{f_2} N_2 \times \cdots \times_{\sigma_k} N_k$ into a complex space form  of constant holomorphic sectional curvature $4\epsilon$
$($or a quaternionic space form  of constant quaternionic sectional curvature $4\epsilon)$. Then we have
 \begin{align*} & \sum_{j=2}^k n_j\frac{\Delta f_j}{f_j} \leq \frac{n^2}{4} H^2  + n_1(n-n_1) \epsilon,\quad n=\sum_{j=0}^kn_j.  \end{align*}\end{corollary}

\begin{remark}  In view of Corollary \ref{C:11.3}, it is interesting to point out that there do exist many isometric minimal immersions from $N_1\times_f N_2$ into Euclidean space with compact $N_2$. For example, a hypercatenoid in
$\mathbb E^{n+1}$ is a minimal hypersurfaces which is isometric to a warped product $ R\times_f S^{n-1}$. Also, for any compact minimal submanifold $N_2$ of $S^{m-1}\subset \mathbb E^m$, the minimal cone $C(N_2)$ is a warped product $
R_+\times_s N_2$ which is also a such example. \end{remark}

\begin{remark} In view of Corollary \ref{C:11.4}, it is interesting to point out that there exist many  minimal submanifolds in Euclidean space which are warped products with harmonic warping function. 

For example, if $N_2$ is a minimal submanifold of the unit $(m-1)$-sphere $S^{m-1}\subset {\mathbb E}^m$,  the
minimal cone $C(N_2)$ over $N_2$ with vertex at the origin of $\mathbb E^m$ is the warped product $ R_+\times_s N_2$ whose warping function $f=s$ is a harmonic function. (Here $s$ is the coordinate function of the positive real line $ R_+$). \end{remark}

\begin{remark}  In views of Corollary \ref{C:11.4}, it is interesting to point out that there exist isometric minimal immersions from warped products $N_1\times_f N_2$  into a hyperbolic space such that the warping function $f$ is an eigenfunction with
negative eigenvalue. For example, $\hbox{\bf R}\times_{\rme^x} {\mathbb E}^{n-1}$ admits an isometric minimal immersion into the hyperbolic space $H^{n+1}(-1)$ of constant sectional curvature $-1$. 
\end{remark}
  
\begin{remark}  In contrast to Euclidean and hyperbolic spaces, the standard $m$-sphere $S^{m}$ admits warped product minimal submanifolds $N_1\times_f N_2$ such that $N_1,N_2$ are both compact. 
The simplest such examples are minimal Clifford tori $M_{k,n-k},\, k=2,\ldots,n-1$, in $S^{n+1}$ defined by
$$M_{k,n-k}=S^k\left(\text{\small$\sqrt{k\over n}$}\,\right)\times S^{n-k}\left(\text{\small$\sqrt{{n-k}\over n}$}\,\right).$$ \end{remark}

\section{Growth estimates of warping functions}

In this section, we provide some growth estimates of the warping function $f$ of the warped product $N_1\times_f N_2$ given in the last section. By applying these, we know in particular  that {\it when the warping function $f$ is an $L^q$ function on a complete noncompact Riemannian manifold $N_1$ for some $q>1$, then for any Riemannian manifold $N_2$ the warped product $N_1\times_f N_2$ does not admit any isometric minimal immersion into any Riemannian manifold with non-positive sectional curvature}. 

To see this, let us assume that $N_1$ is a complete noncompact Riemannian manifold  and  $B(x_0;r)$ denotes the geodesic ball of radius $r$ centered at $x_0 \in N_1$.

We recall some definitions from \cite{cw,WLW}

\begin{definition}A  function $f$  on $N_1$ is said to have \emph{$p$-{finite growth}} $($or, simply, \emph{is $p$-{finite}}$)$ if there exists $x_0 \in N_1$ such that
\begin{equation} \lim_{r \to \infty} \inf\frac{1}{r^p}\int_{B(x_0;r)} |f|^{q}{\rm d}v < \infty ;   \label{12.1} \end{equation}
it has \emph{$p$-{infinite growth}} $($or, simply, \emph{is $p$-infinite}$)$ otherwise. \label{def4.1}\end{definition}

\begin{definition}
A function $f$ has \emph{$p$-mild growth} $($or, simply, \emph{is $p$-mild}$)$ if there exists  $ x_0 \in N_1\, ,$ and a strictly increasing sequence of $\{r_j\}^\infty_0$ going to infinity, such
that for every $l_0>0$, we have
\begin{equation}\sum\limits_{j=\ell_0}^{\infty}\bigg(\frac{(r_{j+1}-r_j)^p}{\int_{B(x_0;r_{j+1})\backslash B(x_0;r_{j})}|f|^q{\rm d}v}\bigg)^{\frac1{p-1}}=\infty \,;   \label{12.2}\end{equation}
and has \emph{$p$-severe growth} $($or, simply, \emph{is $p$-severe}$)$ otherwise.\label{def4.2}\end{definition}

\begin{definition}
A function $f$ has \emph{$p$-obtuse growth} $($or, simply, \emph{is $p$-obtuse}$)$ if there exists $x_0 \in N_1$ such that for every $a>0$, we have
\begin{equation}\int^\infty_a\bigg( \frac{1}{\int_{\partial B(x_0;r)}|f|^q{\rm d}v}\bigg)^\frac{1}{p-1}{\rm d}r = \infty \, ;   \label{12.3}\end{equation}
and has \emph{$p$-acute growth} $($or, simply, \emph{is $p$-acute}$)$ otherwise. \label{def4.3}
\end{definition}

\begin{definition}
A function $f$ has \emph{$p$-moderate growth} $($or, simply, \emph{is $p$-moderate}$)$ if there exist  $ x_0 \in N_1$ and \begin{equation}\notag F(r)\in {\mathcal F} = \left\{F:[a,\infty)\longrightarrow (0,\infty)\, : \,\int^{\infty}_{a}\text{\small$
\frac{{\rm d}r}{rF(r)}$}=+\infty \ \ {rm for \ \ some} \ \ a \ge 0 \right\}\, .
\label{11.8}\end{equation} such that \begin{equation}
\lim \sup _{r \to \infty}\frac {1}{r^p F^{p-1} (r)}\int_{B(x_0;r)} |f|^{q}{\rm d}v < \infty .   \label{12.4}
\end{equation}
And it has \emph{$p$-immoderate growth} $($or, simply, \emph{is $p$-immoderate}$)$ otherwise
 
 {\rm (Notice that the functions in
{$\mathcal F$} are not necessarily monotone.)} \label{def4.4} \end{definition}

\begin{definition}A function $f$ has \emph{$p$-small growth} $($or, simply, \emph{is $p$-small}$)$ if there exists $ x_0 \in N_1\, ,$ such that for every $a>0\, ,$we have
\begin{equation}\int^\infty_a\bigg( \frac{r}{\int_{ B(x_0;r)}|f|^q{\rm d}v}\bigg)^\frac{1}{p-1}{\rm d}r = \infty ;   \label{12.5} \end{equation}
and  has \emph{$p$-large growth} $($or, simply, \emph{is $p$-large}$)$ otherwise.\label{def4.5}
\end{definition}

\noindent The above definitions of   ``\emph{$p$-{finite}}, \emph{$p$-mild},  \emph{$p$-obtuse},  \emph{$p$-moderate}, \emph{$p$-small}" and their counter-parts ``\emph{$p$-infinite}, \emph{$p$-severe},  \emph{$p$-acute}, \emph{$p$-immoderate}, \emph{$p$-large}" growth depend on $q$, and $q$ will be specified in the context in which the definition is used.
\smallskip

From now on, we assume that $N_1$ is a complete noncompact Riemannian $n_1$-manifold and $f$ is a $C^2$-function on $N_1$. Denote by $\tilde M^m_\epsilon$ a Riemannian $m$-manifold with sectional curvatures $\tilde K\leq \epsilon$ for some real number $\epsilon\leq 0$.

In \cite{cw}, the author and S. W. Wei study the growth of warping functions for the warping functions of warped products. They obtain the following results.

\begin{theorem} \label{thm12.1} If $f$ is nonconstant and $2$-finite for some $q > 1$, then for any Riemannian $n_2$-manifold $N_2$ and any isometric immersion $\phi$ of the warped product $N_1\times_f N_2$ into any Riemannian manifold $\tilde M^m_\epsilon$ with $\epsilon\leq 0$, the mean curvature $H$ of $\phi$ satisfies
\begin{align}\label{12.7}H^2 > \frac{4n_1 n_2|\epsilon|}{(n_1 + n_2)^2}\end{align} at some point. \end{theorem}

\begin{corollary} Suppose the squared mean curvature of the isometric immersion $\phi:N_1\times_f N_2\to \tilde M^m_\epsilon$ satisfies
\begin{align}\label{12.8}H^2 \leq \frac{4n_1 n_2|\epsilon|}{(n_1 + n_2)^2}\end{align} everywhere on $N_1\times_f N_2$. Then the warping function $f$ is either a constant or for
every $q > 1\, ,$ $f$ has $2$-infinite growth, i.e., for every $x_0\in N_1$,
\begin{equation} \lim _{r \to \infty}\frac{1}{r^2}\int_{B(x_0;r)} |f|^{q}{\rm d}v = \infty\, .  \label{12.9}
\end{equation} \label{cor12.1} \end{corollary}

\begin{theorem}  \label{thm12.2}  If $f$ is nonconstant and $2$-mild for some $q > 1$, then  for any isometric immersion  of  $N_1\times_f N_2$ into a Riemannian manifold $\tilde M^m_\epsilon$ with  $\epsilon\leq 0$ we have  \eqref{12.7}  at some point.\end{theorem}

\begin{corollary} Suppose that the squared mean curvature of the isometric immersion $\phi:N_1\times_f N_2\to \tilde M^m_\epsilon$  satisfies \eqref{12.8} everywhere on $N_1\times_f N_2$. Then the warping function $f$ is either a constant or has $2$-severe growth for every $q > 1\, .$ \label{cor12.2} \end{corollary}

\begin{theorem}
 \label{thm12.3}  If $f$ is nonconstant and \emph{$2$-obtuse} for some $q > 1\, ,$ then  for  any isometric immersion  of  $N_1\times_f N_2$ into a Riemannian manifold $\tilde M^m_\epsilon$ with  $\epsilon\leq 0$ we have  \eqref{12.7}  at some point.
\end{theorem}

\begin{corollary} Suppose the squared mean curvature of the isometric immersion $\phi:N_1\times_f N_2\to \tilde M^m_\epsilon$  satisfies \eqref{12.8} everywhere on $N_1\times_f N_2$. Then the warping function $f$ is either a constant or has $2$-acute growth for every $q > 1\, .$
 \label{cor12.3}\end{corollary}
\begin{theorem}
 \label{thm4.4}   If $f$ is nonconstant and $2$-moderate for some $q > 1$, then  for any  isometric immersion  of  $N_1\times_f N_2$ into a Riemannian manifold $\tilde M^m_\epsilon$ with  $\epsilon\leq 0$ we have  \eqref{12.7}  at some point.\end{theorem}

\begin{corollary} Suppose the squared mean curvature of the isometric immersion $\phi:N_1\times_f N_2\to \tilde M^m_\epsilon$  satisfies \eqref{12.8} everywhere on $N_1\times_f N_2$. Then the warping function $f$ is either a constant or has $2$-immoderate growth for every $q > 1\, .$   \label{cor4.4}\end{corollary}

\begin{theorem}  \label{thm12.5}  If $f$ is nonconstant and $2$-small for some $q > 1$, then  for any  isometric immersion  of  $N_1\times_f N_2$ into a Riemannian manifold $\tilde M^m_\epsilon$ with  $\epsilon\leq 0$ we have \eqref{12.7}  at some point.
\end{theorem}

\begin{remark}  The assumption on the function $f$ in Theorems \ref{thm12.1}--\ref{thm12.5} above cannot be
dropped.  Otherwise, we have counter-examples that violate \e{12.7}.\end{remark}

\begin{corollary} Suppose that the squared mean curvature of the isometric immersion $\phi:N_1\times_f N_2\to \tilde M^m_\epsilon$  satisfies \eqref{12.8} everywhere on $N_1\times_f N_2$. Then the warping
function $f$ is either a constant or has $2$-large growth for every $q > 1\, .$ \label{cor12.5} \end{corollary}

 \begin{remark} {\rm  Corollaries \ref{cor12.1}--\ref{cor12.5} lead to a dichotomy between
constancy and ``infinity" of the warping functions on complete
noncompact Riemannian manifolds for  isometric immersions of the warped products.}
\end{remark}

In particular, we have the following  result of the author and Wei from \cite{cw}.

\begin{theorem}\label{thm12.6}  Let $f$ be a nonconstant, $L^q$ function on $N_1$ for some $q > 1$, then  for any  isometric immersion
 of  $N_1\times_f N_2$ into a Riemannian manifold $\tilde M^m_\epsilon$ with  $\epsilon\leq 0$ we have  \eqref{12.7}  at some points.
\end{theorem}

From Theorem \ref{thm12.6} we have

\begin{corollary} If $f$ is an $L^q$ function on $N_1$ for some $q > 1$, then for any Riemannian manifold $N_2$ the  warped product  $N_1\times_f N_2$ does not admit any isometric minimal immersion into any  Riemannian manifold with non-positive sectional curvature. \end{corollary}

\begin{theorem}
Suppose $q > 1$ and the warping function $f$ is one of the following: $2$-finite, $2$-mild, $2$-obtuse, $2$-moderate, and $2$-small. If $N_2$ is compact,  then there does not exists an isometric minimal immersion from  $N_1 \times _f N_2 $ into any Euclidean space.
\label{thm6.1} \end{theorem}

\begin{remark}  {\rm In views of the above results,  it is interesting to point out that there do exist isometric minimal immersions from a warped product
$N_1\times_f N_2$  into $\tilde M^m_\epsilon$ with $\epsilon\leq 0$ such that
the warping function $f$ is 2-infinite, 2-severe, 2-acute, 2-immoderate and 2-large for any $q > 1$.

A simple such example is the warped product $\hbox{\bf R}\times_{{\rm e}^x}  {\mathbb E}^{n-1}$ (or $\hbox{\bf R}\times_{\rme^{-x}}  {\mathbb E}^{n-1}$) of constant sectional curvature $-1$ which can be  isometrically immersed in $H^{n+1}(-1)$ as a totally geodesic (hence minimal) submanifold.}\label{rem6.1}\end{remark}

\begin{remark}  {\rm The inequality \eqref{12.7} (resp.  inequality \eqref{12.8})  on  $H^2$  as the assumption for  Theorems  \ref{thm12.1}--\ref{thm12.5} (resp.  assumption of Corollaries \ref{cor12.1}--\ref{cor12.5})  is sharp. This can be seen from the following two examples  (cf. \cite{c1}):

First, let us regard the Euclidean $2k$-space $\mathbb E^{2k}$ as the warped product $\mathbb E^k\times_f \mathbb E^k$ with a
constant warping function $f$. Then $\mathbb E^k\times_f \mathbb E^k$  can be isometrically immersed  in $H^{2k+1}(-1)$ as  a
totally umbilical hypersurface with $H^2=1$. Since $n_1=n_2=k$, the right hand side of \eqref{12.8} is also equal to 1. Thus, this
example satisfies the equality case of \eqref{12.8}.  For the case $k \le 2\, ,$ this example also shows that the nonconstant
assumption in Theorems \ref{thm12.1}--\ref{thm12.5} cannot be dropped.

The second example is the warped product ${\bf R}\times_{\cosh b x}{\bf R}$ of constant negative curvature $-b^2$. This warped
product admits an isometric immersion in $H^3(-1)$ as totally umbilical surface, for $0 < b < 1$. The squared mean curvature of
the immersion satisfies 
$$H^2=1-b^2< \frac{4n_1n_2|\epsilon|}{(n_1+n_2)^2}=1 \hskip.2in {\rm and}\hskip.1in H^2\to 1 \hskip.1in {\rm as}\hskip.1 in b\to 0.$$ The warping function $\cosh bx$ with $0 < b < 1\, $  is a nonconstant and non-harmonic function which is 2-infinite, 2-severe, 2-acute, 2-immoderate and 2-large for any real number $q > 1$. }\end{remark}

 \begin{remark}  {\rm Let $\varphi$ be the function on {\bf R} defined by 
 $$\varphi(x)=\begin{cases} 1 &\text{if $|x|\leq 1$;} \\ x^2 &\text{if $|x|>1$}.\end{cases}$$
 
 Denote by $f$ the smooth out function of $\varphi$ at $\pm 1$. Then $f$ is a subharmonic function on {\bf R} which is  {$2$-finite}, {$2$-mild}, {$2$-obtuse}, {$2$-moderate}, and  {$2$-small}  for any $q\leq 1$; but it is 2-infinite, 2-severe, 2-acute, 2-immoderate and 2-large for $q>1$.

The sectional curvature $K$ of the warped product $N={\bf R}\times_{f}{\mathbb E^{n-1}}$ with this subharmonic warping function $f$ satisfies $K\leq 0$. Let $\tilde M^{n+1}_0={\bf R}\times N$ denote the Riemannian product of the real line and $N$. Clearly, the sectional curvatures of $\tilde M^{n+1}_0$ is bounded above by 0 and the warped product $N$ can be trivially isometrically imbedded in $\tilde M^{n+1}_0$ as a totally geodesic hypersurface. 

This isometric imbedding of $N$ in $\tilde M^{n+1}_0$ satisfies $H^2=\epsilon=0$, which shows that the condition ``$q>1$''  given in Theorems 
\ref{thm12.1}--\ref{thm12.5}  and Corollaries \ref{cor12.1}--\ref{cor12.5} is sharp as well.
 } \end{remark}

\begin{remark}  {\rm The assumption on the warping function $f$ given in Theorem \ref{thm6.1} cannot be dropped, since there do exist  minimal
hypersurfaces in $\mathbb E^{n+1}$ which are isometric to some warped products $N_1\times_f N_2$ with compact $N_2$.
A simple such example is the hypercaternoid in $\mathbb E^{n+1}$.  The hypercaternoid is isometric to a  warped product product
 ${\bf R}\times_{f} S^{n-1}$ with compact $N_2=S^{n-1}$ whose warping function is 2-infinite, 2-severe, 2-acute, 2-immoderate and 2-large for any $q > 1$ according to  Theorem \ref{thm12.1}.} \end{remark}

\section[Applications to Riemannian submersions]{A $\delta$-invariant  and its applications to Riemannian submersions}

\subsection{Riemannian submersions}
 Let $M$ and $B$ denote Riemannian manifolds with $n=\dim M>\dim B=b>0$. A {\it Riemannian submersion} $\pi:M\to B$ is a mapping of $M$ onto $B$ satisfying the following two axioms:
\vskip.04in

\noindent (S1) $\pi$ {\it has maximal rank;}
\vskip.04in

\noindent (S2) {\it the differential  $\pi_*$ preserves lengths of horizontal vectors. }
\vskip.04in

 For each  $p\in B$, $\pi^{-1}(p)$ is an $(m-b)$-dimensional submanifold of $M$.  The submanifolds $\pi^{-1}(p),p\in B,$ are called fibers. A vector  on $M$ is called {\it vertical} if it is tangent to fibers; and it is called {\it horizontal} if it is  orthogonal to fibers. We use corresponding terminology for individual tangent vectors as well. 
 Let  ${\mathcal H }$ and ${\mathcal V}$ denote the  horizontal and vertical distributions.
 
  The simplest type of Riemannian submersions is the projection of a Riemannian product manifold on one of its factors. For such Riemannian submersions, both horizontal  and  vertical distributions are  {\it totally geodesic distributions}, i.e., both distributions are completely  integrable and their leaves are totally geodesic submanifolds. 
 
A Riemannian manifold $M$ is said to admit a {\it non-trivial Riemannian submersion} if there exists a Riemannian submersion $\pi:M\to B$ from $M$ onto another Riemannian manifold $B$ such that the horizontal and vertical distributions of the submersion are not both totally geodesic distribution. 

\subsection{A submersion $\delta$-invariant $\breve \delta_\pi$ and its applications}
  
 For a Riemannian submersion $\pi:M\to B$, we consider a  $\delta$-invariant on $M$ defined by
 \begin{align} \label{1.5} \breve \delta_\pi(p)=\tau(p)-\tau({\mathcal H }_p)-\tau({\mathcal V })_p.
 \end{align} 

By applying this {\it submersion $\delta$-invariant,} we have the following results.
 
 \begin{theorem} {\rm \cite{c37}} \label{T:13.1} If a Riemannian manifold admits a non-trivial Riemannian submersion with totally geodesic fibers, then it cannot be isometrically immersed  in any Riemannian manifold of non-positive sectional curvature as a minimal submanifold. 
\end{theorem}

If $\phi_F:F\to \mathbb E^{m_1}$ and $\phi_B: B\to \mathbb E^{m_2}$ are minimal isometric immersions of Riemannian manifolds $F$ and $B$ into Euclidean spaces, then the product immersion of $\phi_F$ and $\phi_B$ is the immersion: \begin{align}\label{13.2}(\phi_F,\phi_B):F\times B\to \mathbb E^{m_1}\oplus \mathbb E^{m_2}\end{align}  which carries $(q,p)\in F\times B$ to $(\phi_F(q),\phi_B(b))$. 
The product immersion $(\phi_F,\phi_B)$ is  also a  minimal isometric immersion.

\begin{theorem} {\rm \cite{c37}} Let $\pi:M\to B$ be a Riemannian submersion with totally geodesic fibers. If $M$ admits a minimal isometric immersion $\phi$  into a Euclidean space, then  locally $M$ is   the Riemannian product  of a fiber $F$ and the base manifold $B$ and $\phi$ is the product immersion $(\phi_F,\phi_B)$ of some minimal isometric immersions $\phi_F:F\to \mathbb E^{m_1}$ and $\phi_B: B\to \mathbb E^{m_2}$ into some Euclidean spaces. \end{theorem}

The proof of these two theorems is based on the following.

\begin{theorem} {\rm \cite{c37}}  Let  $\pi:M\to B$ be a  Riemannian submersion with totally geodesic fibers. Then, for any isometric immersion of $M$ into a Riemannian $m$-manifold $R^m(\epsilon)$ of constant sectional curvature $\epsilon$, the submersion invariant $ \breve \delta_\pi$ on $M$ satisfies the following inequality:
\begin{align}\label{13.3}  \breve \delta_\pi  \leq \frac{n^2}{4} H^2 +b(n-b)\epsilon. \end{align}
\end{theorem}

If the target manifold in Theorem \ref{T:13.1}  is  negatively  curved, we have

\begin{corollary} {\rm \cite{c37}}   If a Riemannian manifold admits a  Riemannian submersion with totally geodesic fibers, then it cannot be isometrically immersed  in any Riemannian manifold  of negative sectional curvature as a minimal submanifold. 
\end{corollary}

\begin{corollary} {rm \cite{c37}}  Every Riemannian manifold which admits a  non-trivial Riemannian submersion with totally geodesic fibers cannot be isometrically immersed  in any Hermitian symmetric space of non-compact type as a minimal submanifold. \end{corollary}

\begin{remark}{\rm The results obtained above can  be applied to various very large families of Riemannian manifolds, since Riemannian submersions with totally geodesic fibers occur widely in geometry.  For examples, we have:
\vskip.04in 
 
  (i) The well-known Hopf fibrations: \begin{align}\notag &\pi: S^{2n+1}\to CP^n(4)\;\;{\rm and}\; \; \pi:S^{4n+3}\to HP^n(4)\end{align} are Riemannian submersions  with totally geodesic fibers.
\vskip.04in 

(ii) Let $\pi:M\to B$ be a  Riemannian submersion with totally geodesic fibers. If $B'$ is a submanifold of $B$, then the restriction of $\pi$ to $\pi^{-1}(B')$: 
$$ \pi:\pi^{-1}(B')\to B' $$ is a Riemannian submersion with totally geodesic fibers.
For instance,  for any submanifold $N$ of the complex projective $n$-space $CP^n(4)$ of constant holomorphic sectional curvature 4,  $\pi: \pi^{-1}(N) \to N$ is a Riemannian submersion with totally geodesic fibers. 
For this submersion, the invariant $\breve \delta_ \pi$ is given by
\begin{align}\breve \delta_ \pi=||P||^2,\end{align} 
where $P:\mathcal H\to \mathcal H$ is the endomorphism such that $PX$ is  the projection of $\phi X$ onto $\mathcal H$,  $\phi$ being the (1,1)-tensor of the natural Sasakian structure on $S^{2n+1}$.

\vskip.04in 

(iii) If $G$ is a  Lie group equipped with a bi-invariant Riemannian metric and $H$ is a closed subgroup, then the usual Riemannian structure on the homogeneous space $G/H$ is characterized by the fact that the natural mapping  $$\pi:G\to G/H$$ is a Riemannian submersion. 

The fibers of such a submersion are the left cosets of $G$ (mod $H$)  which are totally geodesic.
The invariant $\breve \delta_ \pi$ is given by
\begin{align}\breve \delta_ \pi=\frac{1}{4}\sum_{i,j=1}^b\sum_{b+1}^n \<\, [e_i,e_j],e_s\>^2,\;\; b=\dim H,\end{align} 
where $e_1,\ldots,e_b$  are orthonormal  left-invariant horizontal vector fields  and  $e_{b+1},\ldots,$ $e_n$ an orthonormal basis of the vertical distribution $\mathcal V$.
\vskip.04in 

(iv) On an oriented Riemannian 4-manifold $N$,  there exists an $S^{2}$-bundle $Z$, called the twistor space of $N$, whose fiber over any point $x\in N$ consists of all almost complex structures on $T_{x}N$ that are compatible with the metric and the orientation. It is known that there is one-parameter family of metrics $g^{t}$  on $Z$, making the projection $Z\rightarrow N$ into a Riemannian submersion with totally geodesic fibers.}\end{remark}

  \subsection{Riemannian submersions satisfying the equality case} 
In this subsection, we provide  Riemannian submersions $\pi:M\to B$ with totally geodesic fibers and the isometric immersions of $M$ in $S^N$ which satisfy the equality case of inequality \eqref{13.3} identically. In order to do so, we recall briefly the definition of Hopf's fibration.
 
  Consider  $S^{2n+1}$ as the unit hypersphere in $\mathbb C^{n+1}$ centered at the origin and let $z$ be its unit outward normal. Let $\tilde J$ be the natural almost complex structure on $\mathbb C^{n+1}$.  Then $\tilde Jz$ defines an integrable distribution on $S^{2n+1}$ with totally geodesic leaves. Identifying the leaves as points we obtain the complex projective $n$-space $CP^n$.   By taking as the horizontal distribution, the orthogonal complements to $\tilde J z$ in $TS^{2n+1}$, one  can make this into a Riemannian submersion, known as the Hopf fibration:  $$\pi_C:S^{2n+1}\to CP^n$$ with  great circles as fibers.
 
Similarly,  consider  $S^{4k+3}$ as the unit hypersphere in $\mathbb Q^{k+1}$ and let $z$ its unit outward normal. Let $J_1,J_2,J_3$ be the natural almost complex structures on $\mathbb Q^{k+1}$ with \begin{align}J_1J_2=J_3, \;\; J_2J_3=J_1, \;\; J_3J_1=J_2.\end{align} Then $J_1z, J_2z,J_3 z$ define an integrable distribution on $S^{4k+3}$  with totally geodesic leaves. Identifying the leaves as points we obtain the quaternionic projective $k$-space $QP^k$ which can be made into a Riemannian submersion: 
\begin{align} \label{13.5}\pi_Q:S^{4k+3}\to QP^k \end{align}
 by taking as the horizontal distribution, the orthogonal complements to $J_iz,i=1,2,3,$ in $TS^{4k+3}$. Fibers of $\pi_Q$ are totally geodesic 3-spheres in $S^{4k+3}$. The projection \eqref{13.5} is also known as the Hopf fibration.

The following result shows that there exist many Riemannian manifolds which admit Riemannian  submersions with totally geodesic fibers and which also admit isometric immersions satisfying the equality case of \eqref{13.3} identically  into some unit spheres.
 
\begin{theorem}  \label{T:13.4}  {\rm \cite{c40}} We have:
\vskip.04in

{\rm (a)} Let $B$ be a K\"ahler submanifold of $CP^n$ and let $\pi_C^B:\pi_C^{-1}(B)\to B$ be the restriction of Hopf's  fibration $\pi_C$ to $\pi_C^{-1}(B)\subset S^{2n+1}$. Then the inclusion map $\iota_C: \pi_C^{-1}(B)\to S^{2n+1}$ is an isometric immersion such that
\vskip.04in

 {\rm (a.1)} the fibers of $\pi_C^B$ are fibers of the Hopf fibration  $\pi_C$ and 
\vskip.04in

 {\rm (a.2)} $\iota_C$ satisfies the equality case of \eqref{13.3}  identically with $\epsilon=1$ and  $m=1+b; b=\dim_{\mathbb R} B$.
\vskip.04in

{\rm (b)} Let $B$ be an  open portion of a  totally geodesic $QP^\ell \subset QP^k$ and let $\pi_Q^B:\pi_Q^{-1}(B)\to B$ be  the restriction of Hopf's  fibration $\pi_Q$ to $\pi_Q^{-1}(B)\subset S^{4k+3}$. Then the inclusion map $\iota_Q:\pi_Q^{-1}(B)\to S^{4k+3}$ is an isometric immersion such that
\vskip.04in

 {\rm (b.1)} the fibers of $\pi_Q^B$ are  fibers of Hopf's fibration  $\pi_Q$ and 
\vskip.04in

 {\rm (b.2)} $\iota_Q$ satisfies the equality case of \eqref{13.3}  identically with $\epsilon=1$ and $m=3+b; b=\dim_{\mathbb R} B$.
\end{theorem}

 Let $\pi:M\to B$ and $\pi':M'\to B'$ be two Riemannian submersions with totally geodesic fibers. Then $\pi$ and $\pi'$ are said to be {\it equivalent} provided that there exists an isometry $f:M\to M'$ which induces an isometry $f_B:B\to B'$ so that the following diagram commutes:
\begin{equation}\notag  \begin{CD} M@>f >> M'\\  @V{\pi}VV @VV{\pi'}V\\
B@>f_B >> B'\; . \end{CD} \end{equation}
\vskip.1in

As a converse to Theorem \ref{T:13.4}, we have the following.

\begin{theorem} \label{T:13.5} {\rm \cite{c40}} Let $\pi:M\to B$ be a  Riemannian submersion with totally geodesic fibers.  Then we have:

{\rm (i)}  If  $M$ admits an isometric imbedding $\phi:M\to S^{2n+1}$ which carries fibers of $\pi$ to fibers of  $\pi_C:S^{2n+1}\to CP^n$ and $\phi$ satisfies the equality sign of \eqref{13.3}, then there exists a  K\"ahler submanifold $B_1\subset CP^n$ such that $\pi:M\to B$  is equivalent to $\pi_C:\pi_C^{-1}(B_1)\to B_1$ and  $\phi$  is congruent to the inclusion map: $$\iota_C:\pi_C^{-1}(B_1)\to S^{2n+1}.$$
\vskip.04in

 {\rm (ii)}  If  $M$ admits an isometric imbedding $\phi:M\to S^{4k+3}$ which carries fibers of $\pi$ to  fibers of  $\pi_Q:S^{4k+3}\to QP^k$ and $\phi$ satisfies the equality sign of \eqref{13.3}, then there exists an open portion $B_2$ of some totally geodesic  $QP^\ell\subset QP^k$ such that $\pi:M\to B$ is equivalent to $\pi_Q:\pi_Q^{-1}(B_2)\to B_2$ and  $\phi$  is congruent to the inclusion map: $$\iota_Q:\pi_Q^{-1}(B_2)\to S^{4k+3}.$$ 
  \end{theorem}

Let $\pi:M\to B$ be a  Riemannian submersion. An immersion $\phi:M\to S^N$ is called  {\it mixed-totally geodesic} if its second fundamental form $h$ satisfies \begin{align}h(X,V)=0\end{align} for any horizontal vector $X$ and vertical vector $V$ on $M$.

In view of Theorem \ref{T:13.5}, we mention the following.

\begin{theorem} \label{T:13.6} {\rm \cite{c40}} Let $\pi:M\to B$ be a  Riemannian submersion with totally geodesic fibers. If an isometric immersion $\phi:M\to S^{N}$ carries the fibers of $\pi$ to totally geodesic submanifolds of $S^{N}$ and satisfies the equality case of \eqref{13.3} with $\epsilon=1$, then we have
\vskip.04in

{\rm (1)} $\dim M\geq 3$ and 
\vskip.04in

{\rm (2)} $M$ is immersed as a minimal mixed-totally geodesic submanifold of $S^N$.  \end{theorem}

\subsection{Characterization of Cartan hypersurfaces in terms of $\breve \delta_\pi$}
 The Cartan hypersurface in $S^4\subset
\mathbb E^5$ is defined by the equation:
\begin{equation}\begin{aligned}\notag 2x_5^3&+3(x_1^2+x_2^2)x_5-6(x^2_3+x^2_4)x_5+3\sqrt{3}(x_1^2-x_2^2)x_4 \\&\hskip.7in +6\sqrt{3} x_1x_2x_3=2.\end{aligned}\end{equation} 

\'E. Cartan proved that this hypersurface  is the homogeneous Riemannian manifold $SU(3)/SO(3)$ (equipped with a suitable metric) and its principal curvatures in $S^4$ are given by $$0,\sqrt{3},-\sqrt{3}.$$   
It is also known that the Cartan hypersurface is  a tubular hypersurface
about the Veronese surface with radius $r=\frac{\pi}{ 2}$.

The next theorem classifies hypersurfaces in $S^4$ satisfying the hypothesis of Theorem \ref{T:13.6}.  This result provides us a new characterization of the Cartan hypersurface in terms of the submersion invariant $\breve \delta_\pi$.

\begin{theorem} \label{T:13.7} {\rm \cite{c40}} Let $\pi:M^3\to B$ be a  Riemannian submersion with totally geodesic fibers.     If  $\phi:M^3\to S^{4}$ is a non-totally geodesic isometric immersion which carries fibers of the submersion $\pi$  to totally geodesic submanifolds of $S^{4}$, then $\phi$ satisfies the equality case of \eqref{13.3} (with $\epsilon=1$) if and only if  $\phi$ is congruent to the Cartan hypersurface. 
 \end{theorem}

\subsection{A canonical  cohomology class for Riemannian submersions}
There is a canonical  cohomology class, denoted by $c_\pi(M)$, associated with each Riemannian submersion $\pi:M\to B$ with   orientable base manifold $B$ as follows.

Let $b=\dim B, n=\dim M$, and let $e_1,\ldots,e_n$ be a local orthonormal frame on $M$ which satisfies the following two conditions: 

(i)  $e_{b+1},\ldots,e_n$ are vertical vector fields and 

(ii) $e_1,\ldots,e_b$ are basic horizontal vector fields such that $(e_1)_*,\ldots,(e_b)_*$ gives rise to the positive orientation of $B$. 

Let $\omega^1,\ldots,\omega^n$  be the dual frame of $e_1,\ldots,e_n$. 
Consider  the $b$-form $\omega$ on $M$ defined by
\begin{align} \label{13.9} \omega=\omega^1\wedge \cdots\wedge \omega^b.\end{align}
Then we have $d\omega=0$, since $\omega$ is the pull back of the volume form of $B$.  Thus $\omega$ defines  a cohomology class: $$c_\pi(M)=[\omega]\in H^b(M;{\mathbb R}).$$

By applying this cohomology class, the author proves the following result.

\begin{theorem} {\rm \cite{c37}} \label{T:13.8} Let $b=\dim B$ and $\pi:M\to B$ be a Riemannian submersion with minimal  fibers and orientable base manifold $B$.  If $M$ is a closed manifold with $H^b(M;{\mathbb R})=0$, then the horizontal distribution $\mathcal H$ of the Riemannian submersion is never integrable. Thus  the submersion is always non-trivial.
\end{theorem}

Since each nonzero harmonic form represents a non-trivial cohomology class, Theorem \ref{T:13.8}  follows  from the following.

\begin{theorem} {\rm \cite{c37}}  Let $\pi:M\to B$ be a Riemannian submersion from a closed manifold $M$ onto an orientable  base manifold $B$.  Then the pull back of the volume element of $B$ is  harmonic if and only if the horizontal distribution
 $\mathcal H$ is integrable and fibers are minimal.\end{theorem}

\section[Applications to  Einstein manifolds]{A $\delta$-invariants and its applications to  Einstein manifolds}
Let $p$ be a point in a Riemannian $n$-manifold $M$ and $q$ a  natural number $\leq n/2$.  For a given point $p\in M$, let $\pi_1,\ldots,\pi_q$ be $q$  mutually orthogonal plane sections in $T_pM$. Define the invariant $K_q^{\inf}(p)$ to be the infinimum of the average of the sectional
curvatures $K(\pi_1),\ldots, K(\pi_q)$, i.e., 
\begin{align}\label{14.1}K_q^{\inf}(p)=\inf_{\pi_1\perp\cdots \perp \pi_q} \frac{K(\pi_1)+\cdots +K(\pi_q)}{q},\end{align}
where $\pi_1,\ldots,\pi_q$ run over all mutually orthogonal  $q$ plane sections in $T_pM$.

For a natural number $q\leq n/2$, we define a special $\delta$-invariant  $\hat \delta_q^{{\rm Ric}}$ by
\begin{align} \label{14.2} \delta_q^{{\rm Ric}}(p)=\sup_{X\in T^1_pM} {\rm Ric}(X,X)-\frac{2q}{n}K_q^{\inf}(p),\end{align}
 where $n=\dim M$ and $X$ runs over vectors in   $T_p^1M:=\{X\in T_pM:|X|=1\}$.

Recall that a submanifold $M$ of a Riemannian manifold is  called {\it pseudo-umbilical\/} if its mean curvature vector ${\overrightarrow H}$ is nonzero and its shape operator $A_{\overrightarrow H}$ at the mean curvature vector is proportional to the identity map (cf. \cite{c1}).

For the family of  Einstein manifolds, we have  the following results   from \cite{c34}.

\begin{theorem} \label{T:14.1} For any integer $k\geq 2$ and any  isometric immersion of an Einstein $2k$-manifold $M$  into   $R^m(\epsilon)$  with arbitrary codimension, we have
\begin{align}\label{14.3} \delta_k^{{\rm Ric}}\leq 2(k-1)\(H^2+\epsilon\).\end{align}

The equality sign of \eqref{14.3} holds identically if and only  if  one of the following two cases occurs:
\vskip.04in 

$(1)$ $M$ is a minimal Einstein submanifold such that, with  respect to  some suitable orthonormal frame 
$\{e_1,\ldots,e_{2k},e_{2k+1},\ldots,e_m\}$, we have \begin{align}\notag\font\b=cmr8 scaled \magstep2 \def\bigzerol{\smash{\hbox{ 0}}}
\def\bigzerou{\smash{\lower.0ex\hbox{\b 0}}}
A_r= \begin{pmatrix} A^r_{1} &  & \bigzerou \\   & \ddots&   \\ \bigzerou & &A^r_k
 \end{pmatrix},\;\; r=2k+1,\ldots,m,\end{align}  where $A^r_j,  j=1,\ldots,k,$ are symmetric $2\times 2$   submatrices satisfying  $${\rm trace}\,(A^r_1)=\cdots={\rm trace}\,(A^r_k)=0.$$
\vskip.04in

$(2)$ $M$ is a pseudo-umbilical Einstein submanifold such that,  with respect to some suitable orthonormal frame  $\{e_1,\ldots,e_{2k},e_{2k+1},\ldots,e_m\}$, we have
\begin{align}\notag\font\b=cmr8 scaled \magstep2 \def\bigzerol{\smash{\hbox{ 0}}}
\def\bigzerou{\smash{\lower.0ex\hbox{\b 0}}} A_r= \begin{pmatrix} A^r_{1} &  & &\bigzerou
\\   & \ddots&   \\   & &A^r_k\\ \bigzerou & & & \bigzerou
 \end{pmatrix},\; r=2k+2,\ldots,m,\end{align} 
where $A^r_j,\, j=1,\ldots,k,$ are symmetric $2\times 2$  submatrices  satisfying 
$${\rm trace}\,(A^r_1)=\cdots={\rm trace}\,(A^r_k)=0.$$  \end{theorem} 

\begin{remark} \label{R:14.1} The following example illustrates that  inequality \eqref{14.3} does not hold for arbitrary  submanifolds in general.    \end{remark}

\begin{example} \label{E:14.1} Consider the following spherical  hypercylinder:
$$M:= S^2(1)\times \mathbb E^{2k-2}\subset \mathbb E^{2k+1}.$$
We have $\delta_k^{{\rm Ric}}=1$ and $H^2=1/k^2$ on $M$ which  imply that $$\delta_k^{{\rm Ric}}=1>\frac {2(k-1)}{k^2}=2(k-1)H^2$$ for $k\geq 2.$
\end{example}

\begin{theorem} \label{T:14.2} Let $\phi:M\to  R^m(\epsilon)$ be an  isometric immersion of an Einstein $n$-manifold $M$  into $R^m(\epsilon)$.  Then, for every natural number $q<n/2$, we have
\begin{align}\label{14.4} \delta_q^{{\rm Ric}}\leq \frac{n(n-q-1)}{n-q} H^2+\(n-1-\frac{2q}{n}\)\epsilon.\end{align}

The equality sign of \eqref{14.4} holds identically if and only if   $M$ is totally geodesic.
\end{theorem} 

\begin{remark} \label{R:14.2} The next example shows that  inequality \eqref{14.4} does not hold for arbitrary  submanifolds in general as well.    \end{remark}

\begin{example} \label{E:14.2} For the spherical hypercylinder: $S^{n-q}(1)\times \mathbb E^{q}\subset \mathbb  E^{n+1}$, we have $$\delta_q^{{\rm Ric}}=n-q-1,\;\; H^2=\frac{(n-q)^2}{n^2}$$   
for $q<n/2$ which imply that $$\delta_q^{{\rm Ric}}> \frac{n(n-q-1)}{n-q}H^2.$$ \end{example}

Some consequences of Theorem \ref{T:14.1} and   Theorem \ref{T:14.2} are the following.

\begin{corollary}\label{C:14.1} If  a Riemannian manifold  $M$ admits  an isometric immersion  into a Euclidean space which satisfies \begin{align}\delta_q^{{\rm Ric}}>\frac{n(n-q-1)}{n-q}H^2,\;\; n=\dim M,\end{align} for some natural number $q\leq n/2$ at some points, then $M$ is not Einstein.
\end{corollary}

This corollary  applies to a large family of Riemannian manifolds.  For instance, Example \ref{E:14.1} and Corollary \ref{C:14.1} imply immediate that $S^2\times \mathbb E^{2k-2}$ is not Einstein.

Theorem \ref{T:14.1} and  Theorem \ref{T:14.2}  also imply the following.

\begin{corollary}\label{C:14.2} If  an Einstein $n$-manifold   satisfies \begin{align}\delta_q^{{\rm Ric}}>\(n-1-\frac{2q}{n}\)\epsilon \end{align} for some  natural number $q\leq n/2$ at some points, then  it admits no minimal isometric immersion into $R^m(\epsilon)$ regardless of codimension.\end{corollary}

\begin{corollary}\label{C:14.3} Let $M$ be a compact Einstein $n$-manifold with finite fundamental group $\pi_1$ or with null first Betti number, i.e., $b_1=0$. If there is a  natural number $q\leq n/2$  such that $\delta_q^{{\rm Ric}}>0$, then $M$ admits no Lagrangian isometric
 immersion into any complex n-torus or complex Euclidean $n$-space. \end{corollary}

\section[Applications to  conformally flat manifolds]{A $\delta$-invariant  and its applications to  conformally flat manifolds}

  Let $Ric$ denote the Ricci tensor of a Riemannian $n$-manifold $M$. For each $\ell$-subspace
$L$ of $T_p(M),\, p\in M$, we define the {\it Ricci curvature}  $S(L)$ of $L$ as the trace of the
restriction of the Ricci tensor Ric on $L$, i.e., 
\begin{align}\label{15.1} S(L)={\rm Ric}(e_1,e_1)+\ldots +{\rm Ric}(e_\ell,e_\ell)\end{align}
for an orthonormal basis $\{e_1,\ldots,e_\ell\}$ of $L$. 

For each $k$-tuple $(n_1,\ldots,n_k)$ in  ${\mathcal S}(n)$, we consider in \cite{c36} a special $\delta$-invariant  $\sigma(n_1,\ldots,n_k)$ for conformally flat $n$-manifolds which is defined by
\begin{align}\label{15.2}\sigma(n_1,\ldots,n_k)=\tau- \inf \left\{ \frac {(n-1)\sum_{j=1}^k (n_j-1)}
{(n-1)(n-2)+\sum_{j=1}^k n_j(n_j-1)}  S(L_j)\right\},\end{align}
where $L_1,\ldots,L_k$ run over all $k$ mutually orthogonal subspaces of $T_pM$ such that  $\dim L_j=n_j,\, j=1,\ldots,k$. 
\vskip.1in

The following  general optimal inequality   for  the family of conformally flat  submanifolds was obtained by the author in  \cite{c36}.

\begin{theorem}\label{T:15.1}  Let $M$ be a conformally flat $n$-manifold isometrically immersed in a Riemannian $m$-manifold $R^m(\epsilon)$.  Then, for each $k$-tuple $(n_1,\ldots,n_k)$ in ${\mathcal S}(n)$, we have
\begin{align}\label{15.3} \sigma(n_1,\ldots,n_k)\leq \alpha(n_1,\ldots,n_k) H^2+\beta(n_1,\ldots,n_k)\epsilon,\end{align}
where 
\begin{equation}\begin{aligned}\notag 
&\alpha(n_1,\ldots,n_k)=\text { $\frac{n^2 (n-1)(n-2)\(n+k-1-\sum_{j=1}^k n_j\)}{2\((n-1)(n-2)+\sum_{j=1}^k n_j(n_j-1)\)\(n+k-\sum_{j=1}^k n_j\)}$},
\\& \beta(n_1,\ldots,n_k)=\frac{(n-1)(n-2)\big(n(n-1)-\sum_{j=1}^k n_j(n_j-1)\big)} {2(n-1)(n-2)+2\sum_{j=1}^k n_j(n_j-1)}.
\end{aligned}\end{equation}

The equality case of inequality  \eqref{15.3}  holds at a point $p\in M$ if and only if, there exists an  orthonormal basis 
$e_1,\ldots,e_m$ at $p$, such that  the shape operators of $M$ in $R^m(\epsilon)$ at $p$ take the following
 forms:
\begin{align}\label{15.4}  \font\b=cmr8 scaled \magstep2 \def\bigzerol{\smash{\hbox{ 0}}} \def\bigzerou{\smash{\lower.0ex\hbox{\b 0}}} A_r=\begin{pmatrix} A^r_{1} & \hdots & 0 \\ \vdots  & \ddots& \vdots &\bigzerou \\ 0 &\hdots &A^r_k& \\
 \\&\bigzerou & &\mu_rI \end{pmatrix} ,\quad  r=n+1,\ldots,m, \end{align}
where  $I$ is an identity matrix and  $A^r_j$'s are symmetric $n_j\times n_j$  submatrices satisfying 
\begin{align}\label{15.5}\hbox{\rm trace}\,(A^r_1)=\cdots=\hbox{\rm trace}\,(A^r_k)=\mu_r.\end{align} 
\end{theorem}

 \begin{remark}
Inequality \eqref{15.3} does not hold for arbitrary submanifolds   in general. This can be seen from the following example. \end{remark}
 
 \begin{example}\label{E:2}
 Let $$M=:S^{n-2}(1)\times \mathbb E^2\subset  \mathbb E^{n+1}=\mathbb E^{n-1}\times \mathbb E^2$$ denote the standard isometrically embedding of  $S^{n-2}(1)\times \mathbb E^2$ in $\mathbb E^{n+1}$ as a spherical hypercylinder over the unit $(n-2)$-sphere. If we choose  $k=1, n_1=2, $ then we have $\inf S(\pi)=0$,  where $\pi$ runs over all   2-planes of $T_qM$ at a given point $q\in M$. Hence, by \eqref{15.2},  we have
\begin{align} \sigma(2)=\tau-\frac{n-1}{n^2-3n+4} \inf  S(\pi)=\frac{(n-2)(n-3)}{2}.\end{align}

On the other hand, we have 
\begin{align} \alpha(2)=\frac{n^2(n-2)^2}{2(2n^2-3n+4)},\;\;\;  H^2=\frac{(n-2)^2}{n^2}.\end{align}
Hence, we obtain
 $$\frac{(n-2)(n-3)}{2}=\sigma(2)>\alpha(2)H^2=
\frac{(n-2)^4}{2(2n^2-3n+4)}$$ for $n>4$, which shows that
the  inequality \eqref{15.3} does not hold  for an arbitrary submanifold in a Euclidean space in general. \end{example}  

An immediate application of Theorem \ref{T:15.1} is the following.

\begin{corollary}\label{C:15.1} If  a Riemannian $n$-manifold $M$  admits  an isometric immersion  into a Euclidean space whose $\sigma$-invariant satisfies $$\sigma(n_1,\ldots,n_k)>\alpha(n_1,\ldots,n_k)H^2$$ at some points in $M$ for some  $(n_1,\ldots,n_k)$ in ${\mathcal S}(n)$, then $M$ is not conformally flat. \end{corollary}

For instance, by applying this corollary, we conclude from Example \ref{E:2} that, for any $n\geq 4$,  the Riemannian product $S^{n-2}\times \mathbb E^2$  is not conformally flat.  On the other hand, it is well-known that  $S^{n-1}\times {\bf R}$ is a conformally flat space.

Two other  consequences of Theorem \ref{T:15.1} are the  following obstructions  to minimal and Lagrangian immersions.

\begin{corollary}\label{C:1.2}  Let $M$ be  a conformally  flat $n$-manifold $M$. If there exist  a $k$-tuple $(n_1,\ldots,n_k)$ in ${\mathcal S}(n)$ such that $$\sigma(n_1,\ldots,n_k)>0$$ at some points in $M$, then $M$ does not admit any minimal isometric immersion into a Euclidean space.
\end{corollary}

\begin{corollary}\label{C:15.3}  Suppose that $M$ is a  compact conformally flat $n$-manifold either with finite fundamental group $\pi_1$ or with null first betti number $b_1$. If there exists  a $k$-tuple  $(n_1,\ldots,n_k)$ in ${\mathcal S}(n)$ such that the $\sigma$-invariant
$$\sigma(n_1,\ldots,n_k)>0$$ at some points in $M$, then $M$ admits  no  Lagrangian isometric
 immersion into any complex n-torus or into the complex Euclidean $n$-space.
\end{corollary}
 
 \begin{remark} The condition on  $\sigma{(n_1,\ldots,n_k)}$ given 
in Corollary \ref{C:15.3} is  sharp. This is illustrated by the following
example. \end{remark}

\begin{example} Consider  the Whitney  immersion
$w_a:S^n\to {\mathbb C}^n$ defined by
\begin{align}\label{15.8}w_a(y_0,y_1,\ldots,y_n)= \frac{a(1+\text{i}\,y_0)}{1+y_0^2} (y_1,\ldots,y_n),\;\;\; a>0,\end{align} with $y_0^2+y_1^2+\cdots+y_n^2=1$. 

The Whitney $n$-sphere $W^n_a$ is the topological  $n$-sphere $S^n$ endowed with the induced metric via \eqref{15.8}. The Whitney $n$-sphere is a conformally flat space and the Whitney immersion  is a Lagrangian immersion which has a unique self-intersection point at
$w_a(-1,0,\ldots,0)=w_a(1,0,\ldots,0)$. 

For any $k$-tuple $(n_1,\ldots,n_k)$ in ${\mathcal S}(n)$, we have $$\sigma{(n_1,\ldots,n_k)}\geq 0$$ on $W^n_a$ with respect to the induced metric. Moreover, $\sigma{(n_1,\ldots,n_k)}=0$ holds only at the unique point of self-intersection. From these one may conclude  that the condition on the $\sigma$-invariant  in Corollary \ref{C:15.3} is  sharp. \end{example}

\begin{remark} Let $F:S^1\to\hbox{\bf C}$ be the unit circle  in the complex plane defined by
$F(s)=\rme^{\i s}$ and let $\iota:S^{n-1}\to E^n$  $(n\geq 3)$ be the unit hypersphere in $E^n$
centered at the origin. Denote by $$f:S^1\times S^{n-1}\to \hbox{\bf C}^n\,$$ the complex extensor given by $$f(s,p)=F(s)\otimes \iota(p),\quad p\in S^{n-1}.$$Then $f$ defines an isometric Lagrangian immersion of the conformally flat space $M=:S^1\times S^{n-1}$ into $\hbox{\bf C}^n$ which carries each pair $\{(u,p),\,(-u,-p)\}$  of points in
$S^1\times S^{n-1}$ to a point in $\hbox{\bf C}^n$. Clearly, we have $\pi_1(M)=\hbox{\bf Z}$ and $b_1(M)=1$. Moreover, for each $k$-tuple $(n_1,\ldots,n_k)\in {\mathcal S}(n)$,   $\sigma(n_1,\ldots,n_k)$ is a positive constant  on $M$. This example illustrates that both conditions on $\pi_1(M)$ and $b_1(M)$ in Corollary \ref{C:15.3} are necessary.
\end{remark}

 A Riemannian submanifold $M$  is  {\it minimal}   if its second fundamental form $h$ satisfies $\sum_{i=1}^n h(e_i ,e_i)=0$ with respect to some orthonormal  frame $e_1,\ldots,e_n$.  On the other hand, the notion of coordinate-minimality is introduced in \cite{CG2} as follows:
A Riemannian submanifold  is said to be  {\it coordinate-minimal} with respect to a coordinate system $\{x_1,\ldots,x_n\}$ (or with respect to 
$g=\sum_{i,j=1}^n g_{ij}{\rm d}x_i\otimes {\rm d}x_j$) if its second fundamental form  satisfies $$\sum_{i=1}^n h(\partial_{x_i},\partial_{x_i})=0,$$ where $\partial_{x_i}$ denotes the coordinate vector $\partial/\partial x_i$.

Although coordinate-minimal surfaces are not necessary minimal, they  do share some important properties with minimal surfaces in real space forms.

\begin{proposition} {\rm \cite{CG2}} Let $\varphi:N\to R^m(\epsilon)$ be a coordinate-minimal surface with respect to a coordinate system $\{s,t\}$. Then we have:
\vskip.06in

{\rm (a)} At each point $p\in N$, the second fundamental form  at $p$ takes the form:
\begin{align} h(\partial_s,\partial_s)=-h(\partial_t,\partial_t)=\lambda \hat e_3+\phi \hat e_4,\;\; h(\partial_s,\partial_t)=\mu \hat e_3\end{align}
for some real numbers $\lambda,\mu,\phi$ and  orthonormal normal vectors $\hat e_3,\hat e_4$.
\vskip.06in

{\rm (b)} The Gauss curvature $G$ of $N$  satisfies $G\leq \epsilon$. 
\vskip.06in

{\rm (c)}  $G=\epsilon$ holds identically if and only if $\varphi:N\to R^m(\epsilon)$ is totally geodesic.
\end{proposition}

The following result of the author and O. J. Garay \cite{CG2}  completely classifies  conformally flat $n$-manifolds $(n\geq 4)$ in a Euclidean space which satisfy the equality case of \eqref{15.3} with $k=1$ and $n_1=2$.
  
 \begin{theorem}\label{T:15.2} Let $\phi:M\to {\mathbb E}^m$ be an isometric immersion of a conformally flat $n$-manifold with $n\geq 4$ into Euclidean $m$-space $\mathbb E^{m}$ with arbitrary codimension.  Then we have
 \begin{align}\label{15.9} \sigma(2)\leq \frac{n^2(n-2)^2}{2(n^2-3n+4)} H^2.\end{align} 
 The equality case of \eqref{15.9}  holds identically if and only if  one of the following holds:
\vskip.04in
 
{\rm (1)}  $M$ is an open part of a totally geodesic $n$-plane.
\vskip.04in

{\rm (2)}  $M$ is an open part of a spherical hypercylinder $S^{n-1}\times {\bf R}$ in an affine $(n+1)$-subspace of $\mathbb E^m$.
\vskip.04in

{\rm (3)} $M$ is an open part of a round hypercone in an affine $(n+1)$-subspace of $\mathbb E^m$.
\vskip.04in

{\rm (4)}  $m\geq n+3$ and $M$ is the loci of $(n-2)$-spheres  defined by
\begin{align}\label{15.10} \phi=\( \Psi(s,t),\Big(\text{ $\frac{1}{4c^2}$}-c^2(s^2+t^2)\Big)F\),\end{align}
where $c$ is a positive number, $ F:S^{n-2}\to \mathbb E^{n-1}$ is a unit hypersphere in $\mathbb E^{n-1}$, and $\Psi: N_1\to \mathbb E^{m-n+1}$ is 
a  coordinate-minimal isometric immersion with respect to  
\begin{align} \label{15.11}  g_{N_1}=(1-4c^4 s^2){\rm d}s^2-8c^4st {\rm d}s {\rm d}t+(1- 4c^4 t^2){\rm d}t^2.\end{align} 

{\rm (5)}  $m\geq n+3$ and $M$ is the loci of $(n-2)$-spheres  defined  by
$$\phi=(P(s,t),f(s,t)F),$$ where $F$ is a unit hypersphere in $\mathbb E^{n-1}$ and $P:N_2\to \mathbb E^{m-n+1}$ is a coordinate-minimal surface with respect to \begin{align} \label{15.12} g_P=(f\Delta  f+f_t^2){\rm d}s^2-2f_sf_t {\rm d}s {\rm d}t+(f\Delta  f+f_s^2){\rm d}t^2, \end{align}
where $\Delta=-(\partial_s^2+\partial_t^2)$ and $f$ is a positive  solution of  the  system:
\begin{equation}\begin{aligned}\label{15.13}  &(\Delta f)K_s=(\Delta f)_s-f_s,\;\;  (\Delta f)K_t=(\Delta f)_t+f_t,\;\; \Delta f>0 \end{aligned}  \end{equation}
with $K=\ln (f^2\Delta \ln f)$.\end{theorem} 

\begin{remark} It is proved in \cite{CG2} that there exist many coordinate-minimal isometric immersion $\Psi: N_1\to \mathbb E^{m-n+1}$ with respect to the metric \e{15.11}. Moreover, it is also proved in \cite{CG2} that there exist coordinate-minimal isometric immersion $P N_2\to \mathbb E^{m-n+1}$ with respect to the metric \e{15.12} which satisfies system \e{15.13}.
\end{remark}

\section[A $\delta$-invariant for contact manifolds]{A $\delta$-invariant for contact manifolds and its applications}

Let $(M,\phi,\xi,\eta,g)$ be  an almost contact metric  $(2n+1)$-manifold.  For any natural number $k\in [2,2n]$, we define a {\it contact $\delta$-invariant} $\delta^c(k)$ by (cf. \cite{cmihai})
\begin{align} \label{4.5} \delta^c(k)(p)=\tau (p)-\inf_{L^k_\xi} \tau (L^k_\xi),\end{align}
where $L^k_\xi$ runs over all linear $k$-subspace of $T_xM$ 
containing  $\xi$.

As application of the invariant $\delta^c(k)$ we have the following results of Chen and I. Mihai from \cite{cmihai}.

\begin{theorem} \label{T:16.1} Let $(M,\phi,\xi,\eta,g)$ be an  almost contact metric $(2n+1)$-manifold for which $\eta$ is a contact structure. Then, for any integer $k\in [2,2n]$ and any isometric immersion of $M$ into a real space form $R^m(\epsilon)$, we have:
\begin{align} \label{16.2} \delta^c(k)\leq \frac{(2n+1)^2(2n-k+1)}{2(2n-k+2)}\|H\|^2+\frac{1}{2}\big
(2n(2n+1)-k(k-1)\big)\epsilon. \end{align}

Moreover, the equality case holds identically if and only if we have:
 \begin{enumerate}
  \item[{\rm (1)}] $n+1\leq k\leq 2n$;
 
  \item[{\rm (2)}]  With respect to some suitable orthonormal basis $\{e_1,...,e_{2n+1},e_{2n+2},...,e_m\}$ with $\,\xi=e_1$, the shape operator of $M$ takes the following form:
\begin{align} \label{16.3} A_r=\left ( \begin{matrix} A^r_{k+1} & 0\\ 0 & \mu_r I \end{matrix}\right ),\quad r\in\{2n+2,...,m\}, \end{align}
where $A^r_{k+1}$ are symmetric $(k+1)\times (k+1)$ submatrices satisfying $\,\mbox{\rm trace}\;A^r_{k+1}=\mu_r$ for $r=2n+2,\ldots,m$.
 \end{enumerate} \end{theorem}

\begin{theorem}  \label{T:16.2} If $f:M\to  R^m(c)$ is  an isometric immersion of a $K$-contact
$(2n+1)$-manifold $M$ into a  real space form $R^m(\epsilon)$ which satisfies the equality case of \e{16.2} for some integer $k\in [n+1,2n]$, then we have $\epsilon\geq 1$. 

In particular, when $\epsilon=1$, the $K$-contact structure on $M$ is Sasakian. \end{theorem}

\begin{theorem} \label{T:16.3} If a $K$-contact  $(2n+1)$-manifold $M$ admits an isometric immersion into an $m$-sphere $S^m(\epsilon)$ of constant curvature $\epsilon$ which satisfies the equality case of \eqref{16.2} with $\,k=2$ identically, then we have

$(1)$ $\,\epsilon=n=1$. 

$(2)$ $M$ is a Sasakian manifold of constant curvature one. 

$(3)$ The immersion is totally geodesic. \end{theorem}

\begin{theorem}\label{T:16.4}  If a $K$-contact $(2n+1)$-manifold $M$ admits an isometric immersion into  $S^m(1)$ which  satisfies the equality case of \eqref{16.2} with $k=n+1$, then    $M$ is a minimal submanifold of $S^m(1)$.\end{theorem}

\begin{theorem} \label{T:16.5} If a $K$-contact  $(2n+1)$-manifold $M$ admits an isometric immersion into $S^m(1)$ which satisfies the equality case of \eqref{16.2} with $\,k=3$, then  $\,n=2$ and we have either
\vskip.04in

$(a)$ $M$ is a Sasakian manifold of constant curvature one isometrically immersed in $S^m(1)$ as a totally geodesic submanifold, or
\vskip.04in

$(b)$ $M$ is a Sasakian $5$-manifold foliated by Sasakian $3$-manifolds of constant curvature one and $M$ is isometrically immersed in $S^m(1)$ as a minimal submanifold of codimension at least two. Moreover, leaves of the foliation are immersed as totally geodesic submanifolds of $S^m(1)$. \end{theorem}

\begin{remark} Let  $\hat f:N\to CP^m(4)$ be a K\"ahler  immersion of a K\"ahler surface  $N$ into $CP^m(4)$ with relative nullity two at each point. Then the pre-image $M^5:=\pi^{-1}(N)$ of $N$ via the Hopf fibration $\pi$ is a non-totally geodesic Sasakian 5-manifold  in $S^{2m+1}(1)$ which satisfies the equality case of \eqref{16.2} for $k=3$. \end{remark}

\begin{remark}
It follows from Theorem \ref{T:16.2} that there do not exist $K$-contact manifolds in Euclidean spaces which satisfies the equality case of \eqref{16.2} for any integer $k\in [n+1,2n]$. 

In contrast, the following example shows that there do exist  almost contact metric manifolds in Euclidean spaces which satisfy the equality case of \eqref{16.2} identically.  \end{remark}

\begin{example} Consider the cylindrical hypersurface $$f:M:={\mathbb R}\times S^2(1)\to{\mathbb E}^4$$ defined by
\begin{align}&\label{16.4} f(t,\theta,\varphi)=\big(t,\cos\theta\cos\varphi,\sin\theta\cos\varphi,\sin\varphi\big),\end{align} where $\mathbb E^4$ is the Euclidean 4-space endowed with the flat Riemannian metric: 
\begin{align}&\label{16.5} g={\rm d}t^2+{\rm d}\varphi^2+\cos^2\varphi {\rm d}\theta^2.\end{align}

Define an almost contact metric structure $(\phi,\xi,\eta,g)$ on $M$ by
\begin{equation}\begin{aligned} \notag &\eta =\cos\theta {\rm d}t +\sin\theta {\rm d}\varphi,\;\;
\\&\xi=\cos\theta\frac{\partial}{\partial t}+
\sin\theta\frac{\partial}{\partial\varphi},\\& \phi \left(\frac{\partial}{\partial t}\right)=-\tan\theta
\frac{\partial}{\partial\theta},\;\;
\end{aligned}\end{equation} \begin{equation}\begin{aligned}\notag  & \phi \left(\frac{\partial}{\partial
\varphi}\right)=\frac{\partial} {\partial\theta},\\&
\phi \left(\frac{\partial}{\partial
\theta}\right)=\cos\varphi \left(\sin\theta
\frac{\partial}{\partial t}-
\cos\theta \frac{\partial}{\partial\varphi}\right).\end{aligned}\end{equation}

Consider the orthonormal frame $\{e_1,e_2,e_3\}$ on $M$  given by
$$e_1=\xi,\quad e_2=-\sin\theta\frac{\partial}{\partial t}+\cos\theta\frac{\partial}{\partial\varphi},\quad e_3=\sec\varphi\frac{\partial}{\partial \theta}.$$
With respect to this frame, we have
$$\phi e_1=0,\quad \phi e_2=e_3,\quad\phi e_3=-e_2.$$
It is easy to verify that $g$ satisfies $g(\phi X,\phi Y)=g(X,Y)-\eta (X)\eta (Y).$ Moreover, we also have $\eta\wedge {\rm d}\eta={\rm d}t \wedge
{\rm d}\theta\wedge {\rm d}\varphi\ne 0$ on $M$.
So, $(M,\phi,\xi,\eta,g)$ is an almost contact metric manifold. It is easy to verify that
\begin{equation}\begin{aligned}&\nabla_{\partial_t}\partial_t=\nabla_{\partial_t}\partial_\varphi= \nabla_{\partial_t}\partial_\theta=\nabla_{\partial_\varphi}\partial_\varphi=0,\\& \nabla_{\partial_\varphi}\partial_\theta=-\tan\varphi\partial_\theta,\;\; \\& \nabla_{\partial_\theta}\partial_\theta=\sin\varphi \cos\varphi\partial_\varphi, \end{aligned}\end{equation}
which implies $\nabla_{\xi}{\xi}=0.$ 

This almost contact metric hypersurface $(M,\phi,\xi,\eta,g)$ in $\mathbb E^4$ is  non-$K$-contact and it satisfies the equality case of \eqref{16.2} with $k=2$ identically . \end{example}

\section{$\delta(2)$ and  Lagrangian submanifolds}

If $k=1$ and $n_1=2$, inequality \e{05.10} reduces to
\begin{align}\label{17.1} \delta(2)\leq {{n^2(n-2)}\over {2(n-1)}}H^2+{1\over 2}(n+1)(n-2)\epsilon.\end{align}

Recall that if a Lagrangian submanifold in complex space form $\tilde M^n(4\epsilon)$ satisfies the equality case of \e{17.1}, then it is minimal (cf. Theorem \ref{T:5.5}). Hence, we have the following equality:
\begin{align}\label{17.2} \delta(2)={1\over 2}(n+1)(n-2)\epsilon.\end{align}

For an $n$-dimensional Lagrangian submanifold in $CP^n(4)$ satisfying \e{17.2},  We define for every $p\in M$ the kernel of
the second fundamental form by
$$\mathcal D (p) = \{X \in T_p M\, \vert\, \forall\, Y \in T_pM : h(X,Y) = 0\}. $$
If the dimension of $\mathcal D(p)$ is constant, then it follows from \cite{c5} that
either $M$ is totally geodesic or that the distribution $\mathcal D$ is an
$(n-2)$-dimensional completely integrable distribution.   

Denote by $\mathcal D^\perp$ the complementary orthogonal distribution of $\mathcal D$ in $TM$. Contrast to $\mathcal D$, the distribution $\mathcal D^\perp$ is not necessary integrable in general.

\subsection{Lagrangian submanifolds with constant scalar curvature satisfying   \e{17.2}}

To state the next result, we provide an example of minimal Lagrangian submanifold $M$ with constant scalar curvature  in $CP^3(4)$ which satisfies  equality  \e{17.2} with $n=3$ as follows.

\begin{example} Define two complex structures on ${\bf C}^4$ by
$$ \aligned I(v_1,v_2,v_3,v_4)&=(iv_1,iv_2,iv_3,iv_4)\\
J(v_1,v_2,v_3,v_4)&=(-\bar v_4, \bar v_3,-\bar v_2,\bar v_1). \endaligned $$ Clearly $I$ is the standard complex structure.
The corresponding Sasakian structures on $S^7(1)$ have characteristic vector fields $\xi_1=-I(x)$ and $\xi_2=-J(x)$. Since we consider two complex structures on ${\bf C}^4$, we can consider two different Hopf fibrations $\pi_j:S^7(1)\to C P^3(4)$. The vector field $\xi_j$ is vertical for $\pi_j$.

Now we consider the Calabi curve $\mathcal C_3$ of $ C P^1$ into $C P^3(4)$ of constant Gauss curvature $4/3$, given by
$$ \mathcal C_3(z) =\big[1,\sqrt3 z ,\sqrt 3 z^2,z^3\big] $$
Since $\mathcal C_3$ is holomorphic with respect to  $I$, there is a circle bundle $\pi : M^3 \to  CP^1$ over $CP^1$ and an isometric minimal immersion $\mathcal I : M^3\to S^7(1)$ such that $\pi_1(\mathcal I)=\mathcal C_3(\pi)$.  It is known in \cite{cdvv3} that $M$ has constant scalar curvature and $\mathcal I$ is horizontal with respect to  $\pi_2$ such that the immersion $$\mathcal J:M^3\to  CP^3(4),$$ defined by $\mathcal J = \pi_2(\mathcal I)$ is a minimal Lagrangian
immersion which satisfies the equality  \e{17.2} with $n=3$.
\end{example}

For Lagrangian submanifolds with constant scalar curvature, we have

\begin{theorem} \label{T:17.1} {\rm \cite{cdvv3}} Let $\phi:M^n \to\tilde M^n(4\epsilon)$, $\epsilon \in \{-1,0,1\}$ and $n \ge 3$ be a Lagrangian minimal immersion with constant scalar curvature.  Then $M^n$ satisfies
$$ \delta_M = \tfrac12(n+1)(n-2) \epsilon ,\leqno({1.2}) $$
if and only if either
\vskip.05in

{\rm (1)}  $M^n$ is a totally geodesic immersion, or
\vskip.05in

{\rm (2)} $n=3$, $\epsilon=1$ and $\phi $ is locally congruent to the immersion $\mathcal J:M^3\to CP^3(4)$ defined above.
\end{theorem}

The next example shows that there exists a $3$-dimensional totally real submanifold in $\mathbb CP^3$ with non-constant scalar curvature which satisfies the equality \e{17.2}.  From that example, it follows that the condition of constant scalar curvature in Theorem \ref{T:17.1} cannot be omitted.

\begin{example} Consider a totally geodesic $S^5(1)$ in $S^7(1)$. Let $N$ be a unit vector orthogonal to the linear subspace containing $S^5(1)$.  Let $N^2$ be any minimal surface in $S^5(1)$, immersed by $f$. We define an isometric immersion from the warped product manifold $M^3=(-\pi/2,\pi/2)\times_{\cos t} N^2$ into $S^7(1)$ by
$$x: M^3\to S^7(1) : (t,p)\mapsto (\sin t) N +(\cos t)f(p). $$
This immersion is minimal and satisfies equality \e{17.2}. In order for $x$ to be Legendrian, we have to assume that $S^5(1)$ is contained in a complex hyperplane ${\bf C}^3$ of ${\bf C}^4\supset S^7(1)$ and that $f$ is Legendrian.  It is easy to check that $M^3$ does not have constant scalar curvature if $N^2$ is not totally geodesic.
\end{example}

\subsection{Lagrangian submanifolds  of $CP^n$ with integrable $\mathcal D^\perp$ satisfying  \e{17.2}}

Lagrangian submanifolds  of $CP^n$ with integrable $\mathcal D^\perp$ satisfying  \e{17.2} are completely determined by the author, Dillen, Verstraelen and Vrancken in \cite{cdvv1}.

\begin{theorem} Let $\phi: M^n \longrightarrow CP^n(4)$ be a Lagrangian immersion satisfying   equality \e{17.2} and

{\rm (1)} the dimension of $\mathcal D$ is constant,

{\rm (2)} $\mathcal D^\perp$ is an integrable distribution.

\noindent Then either $\phi$ is totally geodesic or $\phi$ has no totally geodesic points and, up to holomorphic transformations, $\phi(M)$ is contained in the image under the Hopf fibration $\pi : S^{2n+1}(1)\to CP^n(4)$ of the image of one of the immersions described in the  next proposition.  \end{theorem}

\begin{proposition} Let  $S^{2n+1}(1)$ be the unit hypersphere of ${\bf C}^{n+1}$ and consider the orthogonal decomposition
${\bf C}^{n+1} = {\bf C}^3\oplus J(\mathbb E^{n-2})\oplus {\mathbb E}^{n-2}$. Let $f:M^2\to S^5(1)\subset {\bf C}^3 $ be a minimal Legendrian  immersion and consider the hypersphere $S^{n-3}(1)$ in $\mathbb E^{n-2}$.
Then
$$ F:(0,\pi/2)\times_{\cos t} M^2\times_{\sin t}S^{n-3}(1) \to S^{2n+1}(1) : (t,p,q)\mapsto \cos t f(p) + \sin t q $$
is a minimal Legendrian immersion satisfying equality \e{17.2}. Moreover, if $f$ has no totally geodesic points,
then  the dimension of $\mathcal D$ is exactly $n-2$. 

Finally, if we extend $F$ to a map $$\wt F:
(-\pi/2,\pi/2)\times M^2\times S^{n-3}(1) \to S^{2n+1}(1) : (t,p,q)\mapsto \cos t f(p) + \sin t q.$$ Then $\wt F$ fails to be immersive at $t=0$, but the image of $\wt F$ is an immersed minimal Legendrian submanifold. If $f$ is not totally
geodesic, then this image can not be extended further.
\end{proposition}

\subsection{Improved inequality for Lagrangian submanifolds}
For Lagrangian submanifolds in $\tilde M^n(4\epsilon)$, T.  Oprea \cite{Op3} improves inequality \e{17.1} to the following.
\begin{align}\label{17.3} \delta(2)\leq {{n^2(2n-3)}\over {2n+3}}H^2+{1\over 2}(n+1)(n-2)\epsilon.\end{align}

The improved inequality \e{17.3} was proved in \cite{Op3} by using a method of optimizations. Recently, a purely algebraic proof of \e{17.3} is obtained in \cite{BDFV}. Moreover, it is also proved in \cite{BDFV,BMV} that an $n$-dimensional Lagrangian submanifold in a complex space form $\tilde M^n(4\epsilon)$ with $n\geq 4$ is minimal if  the equality  case of \e{17.3} is attained at all points. Notice that in the minimal case, inequality \e{17.3} reduces to \e{17.1}.

Three-dimensional non-minimal Lagrangian submanifolds in $CP^3(4)$ which attain the equality case of \e{17.3} have been constructed in \cite{BV2}. More precisely, J. Bolton and L. Vrancken prove the following.

\begin{theorem} Let $M$ be a $3$-dimensional non-minimal Lagrangian submanifold of $CP^3(4)$  which attains equality at every point in \e{17.3}. Then there is a minimal Lagrangian surface  $\tilde W$ in $CP^2(4)$ such that $M$ can be locally written as $[E_0]$,  where
$$ E_0 =\frac{ \rme^{\i t/3}}{\sqrt{1+b_1^2+\lambda_2^2}}(0,W)+ \frac{(-b_1+{\rm i}\lambda_2)}{\sqrt{1+b_1^2+\lambda_2^2}}(\rme^{{\rm i} t}, 0, 0, 0),$$
where $b_1$ and  $\lambda_2$ are solutions of the following system of ordinary differential equations:
$$\frac{{\rm d}b_1}{{\rm d}t}=-\frac{1+3\lambda_2^2+b_1^2}{3\lambda_2},\hskip.2in \frac{{\rm d}\lambda_2}{{\rm d}t}=\frac{2}{3}b_1.$$

Conversely, any $3$-dimensional Lagrangian submanifold $M$ obtained in this way attains equality at each point in \e{17.3}. \end{theorem}

\section{$\delta(2)$ and $CR$-submanifolds}

For $CR$-submanifolds in complex space forms, there exists a sharp relationship between the
invariant $\delta (2)$ and the squared mean curvature $H^2$:

\begin{theorem} Let $M$ be an $n$-dimensional   $CR$-submanifold in a complex space form $\tilde M^m(4\epsilon)$. Then we have
\begin{align}\label{18.1}\delta(2) \leq \begin{cases}  \dfrac{n^2(n-2)}{2(n-1)}H^2 +\left\{ \dfrac{1}{2}(n+1)(n-2)+3h\right\} \epsilon,&\text{if\ $\epsilon>0;$} \\  \dfrac{n^2(n-2)}{2(n-1)}H^2,&\text{if $ \epsilon=0;$} \\ \dfrac{n^2(n-2)}{2(n-1)}H^2 + \dfrac{1}{2}(n+1)(n-2)\epsilon, &\text{if $ \epsilon<0$},\end{cases}\end{align} where $h$ is the complex dimension of the holomorphic distribution.
\end{theorem}

There exist many $CR$-submanifolds in complex space forms which satisfy the equality cases of
the above inequalities.

For instance, for Hopf hypersurfaces in $CP^m$, we have the following classification theorem.

\begin{theorem} {\rm \cite{c96}} \label{T:18.2}Let $M$ be a Hopf hypersurface of $CP^m(4)\;(m\geq2)$. Then $M$ satisfies
 $$\delta(2) =  \frac{(2m-1)^2(2m-3)}{4(m-1)}H^2 +2m^2-3$$ if and only if 

{\rm (1)} $M$ is an open part of a geodesic sphere of radius ${\pi\over 4}$ in $CP^m(4)$ or 

{\rm (2)} $m=2$ and $M$ is an open part of a tube over a complex quadric curve
$Q_1$ with  radius   $\,r_1=\arctan\left({1\over2}\left(1+\sqrt{5}-\sqrt{2+2\sqrt{5}}\,\right)\right).$
  \end{theorem}
By a Hopf hypersurface in $CP^m$ we mean a real hypersurface such that $J\xi$ is an eigenvector of the shape operator $A_\xi$, where $\xi$ is a unit normal vector field.
  
 Up to rigid motions of $CH^m(-4)$, a horosphere in $CH^m(-4)$ is a real hypersurface defined by the equation
\begin{align}|z_1-z_0|=1.\end{align}

  For hypersurfaces in $CH^m(-4)$, we obtain the following
 characterization of horospheres in $CH^2(-4)$ in terms of $\delta(2)$.

\begin{theorem} {\rm \cite{c96}}  \label{T:18.3} Let $M$ be a real hypersurface of $ CH^m(-4)\;(m\geq2)$. Then 
\begin{align}\label{96}\delta(2)= \frac{(2m-1)^2(2m-3)}{4(m-1)} H^2 -(2m^2-6) \end{align}  holds identically if
and only if $m=2$ and $M$ is an open part of a horosphere in $\mathbb CH^2(-4)$.\end{theorem}

Proper $CR$-submanifolds of complex hyperbolic spaces satisfying the equality case of \e{96} were completely determined by the author and Vrancken \cite{CV2}. More precisely, they obtained the following:

\begin{theorem} \label{T:18.4} Let $U$ be a domain of $\hbox{\bf C}$ and $\Psi:U\to {\bf C}^{m-1}$ be a nonconstant holomorphic curve in ${\bf C}^{m-1}$. Define $z:E^2\times U\to {\bf C}^{m+1}_1$  by
\begin{align} z(u,t,w)=\rme^{\i t} \Bigg(\i u-1-{1\over 2}\Psi(w)\bar\Psi(w),\i u-{1\over 2}\Psi(w)\bar\Psi(w),\Psi(w)\Bigg).\end{align}
Then $\left<z,z\right>=-1$ and the image $z(E^2\times U)$ in $H^{2m+1}_1$ is invariant under the group action of $H^1_1$. Moreover, away from points where $\Psi'(w)=0$, the image $\pi(E^2\times U)$, under the projection
$$\pi:H^{2m+1}_1(-1)\to CH^m(-4),$$ is a proper $CR$-submanifold of $ CH^m(-4)$ which
satisfies \begin{equation} \label{18.3}\delta_M =  {{n^2(n-2)}\over{2(n-1)}}H^2 + {1\over 2}(n+1)(n-2)\epsilon.\end{equation}

Conversely, up to rigid motions of $ CH^m(-4)$, every proper $CR$-submanifold of $CH^m(-4)$ satisfying the equality is obtained in such way.
\end{theorem}

Theorem \ref{T:18.4}  can be regarded as a natural extension of  Theorem \ref{T:18.3}.

\section{$\delta(2)$ and $CMC$ hypersurfaces}

A hypersurface in the unit round sphere $S^{n+1}(1)$ is called {\it isoparametric} if it has constant principal curvatures. It is known  that an isoparametric hypersurface in $S^{n+1}(1)$ is either an open portion of a 3-sphere or an open portion of the product of a circle and a 2-sphere, or an open portion of a tube of constant radius over the Veronese embedding. 

Because every isoparametric hypersurface in $S^{n+1}(1)$ has constant mean curvature ($CMC$) and constant scalar curvature, it is an interesting problem to determine all $CMC$ hypersurfaces with constant scalar curvature. Many results in this directions have been obtained (see, for instance, \cite{chang,cran}).

On the other hand, every isoparametric hypersurface in $S^{n+1}(1)$ or in $\mathbb E^{n+1}$ has constant mean
curvature and constant $\delta(2)$-invariant. Thus, it is also a very natural problem to study $CMC$ hypersurfaces
 with  constant $\delta(2)$-invariant in a real space form $R^{n+1}(\epsilon)$. 

In this respect, we mention the following results of the author and O. J. Garay proved in \cite{CG1}.

\begin{theorem} \label{T:19.1} A CMC hypersurface in the Euclidean $4$-space $\mathbb E^4$ has constant
$\delta(2)$-invariant if and only if it is one of the following:
\vskip.04in

{\rm (1)}   An isoparametric hypersurface;
\vskip.04in

{\rm (2)}   A minimal hypersurface with relative nullity greater than or equal to 1;
\vskip.04in

{\rm (3)}   An open portion of a hypercylinder $N\times {\bf R}$  over a surface $N$ in $\mathbb E^3$ with
constant mean curvature and nonpositive Gauss curvature.
\end{theorem}

\begin{theorem} A CMC hypersurface M in the unit $4$-sphere $S^4(1)$  has constant $\delta(2)$-invariant if and only if one of the following two statements holds:
\vskip.04in

{\rm (1)}  M is an isoparametric hypersurface;
\vskip.04in

{\rm (2)}  There is an open dense subset $U$ of $M$  and a non-totally geodesic isometric minimal immersion $\phi:B^2\to S^4(1)$  from a surface $B^2$ into $S^4(1)$ such that $U$ is
an open subset of $NB^2\subset S^4(1)$, where $NB^2$ is defined by
$$N_pB^2=\Big\{\xi\in T_{\phi(p)} S^4(1)\,:\, \<\xi,\xi\>=1,\; \<\xi,\phi_*(T_pB^2)\>=0\Big\}.   $$
\end{theorem}

As an immediate application of Theorem \ref{T:19.1}  we have the following corollary for complete $CMC$ hypersurfaces.

\begin{corollary} Let $M$ be a non-minimal, complete, CMC hypersurface of the Euclidean 4-space $\mathbb E^4$. Then $M$ has constant $\delta(2)$-invariant if and only if $M$ is one of
the following hypersurfaces:
\vskip.04in

{\rm (1)}  An ordinary hypersphere;
\vskip.04in

{\rm (2)}  A spherical hypercylinder: ${\bf R}\times S^2$;
\vskip.04in

{\rm (3)} A hypercylinder over a circle: $\mathbb E^2\times S^1$. \end{corollary}

\section{$\delta(2)$ and submanifolds of nearly K\"ahler $S^6$}

\subsection{General results of nearly K\"ahler $S^6$} It was proved by E. Calabi \cite{calabi} in 1958  that any oriented  submanifold $M^6$ of the hyperplane  $\hbox{Im}\,\mathcal O$ of the imaginary octonions  carries a $U(3)$-structure (that is, an almost Hermitian structure). For instance, let $S^6\subset \hbox{Im}\,\mathcal O$ be the sphere of unit imaginary vectors; then the right multiplication by $u\in S^6$ induces a linear transformation \begin{align}J_u: {\Cal O} \to\mathcal O,\end{align} which is orthogonal and satisfies $(J_u)^2= - I$. The operator $J_u$ preserves the 2-plane spanned by 1 and $u$ and therefore preserves its orthogonal 6-plane which may be identified with $T_uS^6$.
Thus $J_u$ induces an almost complex structure on $T_uS^6$ which is compatible with the inner product induced by the inner product of $\mathcal O$ and $S^6$ has an almost complex structure.

The almost complex structure $J$ on $S^6$ is a  nearly K\"ahler structure in the sense that the (2,1)-tensor field $G$ on $S^6$, defined by $$G(X,Y) = (\widetilde \nabla_XJ)(Y),$$ is skew-symmetric, where $\widetilde \nabla$ denotes the Riemannian connection on $S^6$. 

The group of automorphisms of this nearly K\"ahler structure is the exceptional simple Lie group $G_2$ which acts transitively on $S^6$ as a group of isometries.

 A. Gray \cite{gray} proved the following:
\vskip.04in 

(1)  {\it Every almost complex submanifold of the nearly K\"ahler $S^6$ is a minimal submanifold, and} 
\vskip.04in 

(2) {\it The nearly K\"ahler $S^6$ has no $4$-dimensional almost complex submanifolds. }
\vskip.04in 

A 3-dimensional submanifold $M$ of the nearly K\"ahler $S^6$ is called Lagrangian if the almost complex structure $J$ on the nearly K\"ahler 6-sphere carries each tangent space $T_xM,\,x\in M$ onto the corresponding normal space $T^\perp_xM$.

N.  Ejiri proved in \cite{Ejiri}  that a Lagrangian submanifold $M$ in  the nearly K\"ahler  $S^6$ is always minimal and orientable.  He also proved that if $M$ has constant sectional curvature, then $M$ is either totally geodesic or has
constant curvature $1/16$. The first  nonhomogeneous examples of Lagrangian submanifolds in the nearly K\"ahler 6-sphere were  described in  \cite{Ejiri86}.

\subsection{$G_2$-equivariant Lagrangian submanifolds and their characterizations}  K. Mashimo \cite{Mas} classified the $G_2$-equivariant Lagrangian submanifolds $M$ of the nearly K\"ahler 6-sphere. It turns out that there are five models, and every equivariant  Lagrangian submanifold in the nearly K\"ahler 6-sphere is $G_2$-congruent to one of the five models. 

These five models can be distinguished by the following curvature properties:
\vskip.1in 

(1)  $M^3$ is totally geodesic ($\delta(2)=2$),
\vskip.04in 

(2) $M^3$ has constant curvature $1/16$ ($\delta(2)=1/8$),
\vskip.04in 

(3) the curvature of $M^3$ satisfies  $1/16\leq K \leq 21/16$ ($\delta(2)=11/8$),
\vskip.04in 
 
(4)  the curvature of $M^3$ satisfies $-7/3\leq K \leq 1$ ($\delta(2)=2$), 
\vskip.04in 
 
(5)  the curvature of $M^3$ satisfies $-1\leq K \leq 1$ ($\delta(2)=2$).
\vskip.1in

F. Dillen, L. Verstraelen and L. Vrancken  \cite{DVV}  characterized models (1), (2) and (3) as the only compact Lagrangian submanifolds in $S^6$ whose sectional curvatures satisfy $K\geq 1/16$. They also obtained an explicit expression for the Lagrangian submanifold of constant curvature $1/16$ in terms of harmonic homogeneous polynomials of degree 6. Using these formulas, it follows that the immersion
has degree 24. Further, they also obtained an explicit expression for model (3).

It follows from the fundamental inequality in Theorem \ref{T:5.1} and Ejiri's result that the invariant
$$\delta(2)=\tau-\inf K$$
always satisfies \begin{align}\delta(2)\leq 2,\end{align} for every Lagrangian submanifold of the nearly K\"ahler $S^6$. 

From above, we know that the models (1), (4) and (5) satisfy the basic equality: $\delta(2)=2$ identically.  On the other hand, 
Chen, F. Dillen, L. Verstraelen and L. Vrancken   proved in \cite{cdvv95} that {\it  models}  (1), (2) {\it and} (3) {\it are the only Lagrangian submanifolds of  the nearly K\"ahler $S^6$ with constant scalar  curvature that satisfy the equality $\delta(2)=2$. }

Many further examples of Lagrangian submanifolds in the nearly K\"ahler $S^6$ satisfying the equality $\delta(2)= 2$ have been constructed  in \cite{cdvv2,cdvv95}.

Recall that a Riemannian $n$-manifold $M$ whose Ricci tensor has an eigenvalue of multiplicity at least $n-1$ is called quasi-Einstein. 

R.  Deszcz, F. Dillen, L. Verstraelen  and L. Vrancken \cite{dedvv} proved that {\it Lagrangian submanifolds of the nearly K\"ahler $6$-sphere satisfying $\delta(2)=2$ are always quasi-Einstein.}

\subsection{Classification of Lagrangian submanifolds in the nearly K\"ahler $S^6$}
The complete classification of Lagrangian submanifolds in the nearly K\"ahler 6-sphere satisfying  $\delta(2)= 2$ was established by Dillen and Vrancken in \cite{DV1}. More precisely, they proved the following:
\vskip.04in

(1) {\it Let $\phi: N_1 \to CP^2(4)$ be a holomorphic curve in $CP^2(4)$, $PN_1$  the circle bundle over $N_1$ induced by the Hopf fibration $\pi: S^5(1) \to CP^2(4)$, and $\psi$ the isometric immersion such that the following diagram commutes:
$$\CD PN_1 @>\psi>> S^5\\
@V VV @VV\pi V\\ N_1 @>\phi>> CP^2(4).
\endCD$$
Then, there exists a totally geodesic embedding $i$ of $S^5$ into the nearly K\"ahler $6$-sphere such that the immersion $i \circ \psi : PN_1 \rightarrow S^6$ is a $3$-dimensional Lagrangian immersion in $S^6$
satisfying equality $\delta(2)=2$.}
\vskip.04in

(2) {\it Let $\bar\phi: N_2 \to S^6$ be an almost complex curve (with second fundamental form $h$)  without totally geodesic points. Denote by $UN_2$ the unit tangent bundle over $N_2$ and define a map
\begin{align}\bar\psi : UN_2 \rightarrow S^6 : v \mapsto \bar\phi_\star(v) \times {{h(v,v)}\over{\Vert h(v,v) \Vert}}.\end{align}  
Then $\bar\psi$ is a (possibly branched) Lagrangian immersion into $S^6$ satisfying equality $\delta(2)=2$. Moreover, the immersion is linearly full in $S^6$.}
\vskip.04in

(3) {\it Let $\bar\phi: N_2\to S^6$ be a (branched) almost complex immersion.  Then, $SN_2$ is a $3$-dimensional (possibly branched) Lagrangian submanifold of $S^6$ satisfying equality  $\delta(2)=2$.}
\vskip.04in

(4) {\it Let $f : M\to S^6$ be a Lagrangian immersion which is not linearly full in $S^6$.  Then $M$ automatically satisfies equality $\delta(2)=2$ and there exists a totally geodesic $S^5$, and a holomorphic immersion $\phi:N_1\to  CP^2(4)$ such that $f$ is congruent to $\psi$, which is obtained from $\phi$ as in} (1).
\vskip.04in

(5) {\it Let $f: M\rightarrow S^6$ be a linearly full Lagrangian immersion of a $3$-dimensional manifold satisfying equality  $\delta(2)=2$.  Let $p$ be a non totally geodesic point of $M$. Then there exists a (possibly branched) almost complex curve $\bar\phi : N_2 \rightarrow S^6$ such that $f$ is locally around $p$ congruent to $\bar \psi$, which is obtained from $\bar \phi$ as in} (3).
\vskip.04in

Let $f : S \rightarrow S^6$ be an almost complex curve without totally geodesic points.  Define
\begin{align}\label{F} F : T_1S \rightarrow S^6(1) : v \mapsto \frac{h(v,v)}{|| h(v,v) ||},\end{align} where $T_1S$ denotes the unit tangent bundle of $S$.

Also, Vrancken \cite{v} showed that the following:
\vskip.04in 

(i) $F$ {\it given by \e{F} defines a Lagrangian immersion if and only if $f$ is superminimal, and}
\vskip.04in 

(ii)  {\it If $\psi : M \rightarrow S^6(1)$ be a Lagrangian immersion which admits a unit length Killing vector field whose integral curves are great circles.  Then there exist an open dense subset $U$ of $M$ such that each point $p$ of $U$ has a neighborhood $V$ such that $\psi : V \rightarrow S^6$
satisfies $\delta(2)=2$, or $\psi: V \rightarrow S^6$ is obtained as in} (i).

\subsection{$CR$-submanifolds in the nearly K\"ahler $S^6$} 
M. Djori\'c and Vrancken  \cite{DjV} study 3-dimensional $CR$-submanifolds in $S^6$ and proved  the following.

\begin{theorem} Let $M$ be a $3$-dimensional minimal CR-submanifold in the nearly K\"ahler $S^6$ satisfying $\delta(2)=2$. Then $M$ is a totally real submanifold or locally M is congruent with the immersion:
\begin{align}\label{djv} &f (t, u, v) = \Big(\cos t \cos u \cos v, \sin t,  \cos t \sin u \cos v,\\&  \hskip.6in \notag \cos t \cos u \sin v, 0, -\cos t \sin u \sin v, 0\Big).\end{align}
\end{theorem}

Notice that \e{djv} can also be described algebraically by the equations:
$$x_5=x_7=0,\quad x_1^2+x_2^2+x_3^2+x_4^2+x_6^2=1,\quad x_3x_4+x_1x_6=0,$$
from which we see that the $CR$-submanifold is a hypersurface lying in a totally geodesic $S^4(1)$.

Very recently, Djori\'c and Vrancken proved in \cite{DjV1} that  3-dimensional $CR$-submanifolds in the nearly K\"ahler $S^6$ satisfying $\delta(2)=2$ must be minimal.

For a 4-dimensional submanifold $M$ of  $S^6(1)$, the fundamental inequality \e{5.1} in Theorem \ref{T:5.1} yields the following inequality for $\delta(2)$:

$$\delta(2)\leq \frac{16}{3}H^2+\frac{15}{2}.$$

Recently, Anti\'c, Djori\'c and Vrancken  \cite{ADV} study 4-dimensional $CR$-submanifolds in $S^6$ and proved  the following.

\begin{theorem}
Let M be a $4$-dimensional minimal $CR$-submanifold in the nearly K\"ahler $S^6$ satisfying $\delta(2)=15/2$. Then M
is locally congruent with the following immersion:
\begin{align}\notag  f (x_1, x_2, &x_3, x_4) = \Big(\cos x_4 \cos x_1 \cos x_2 \cos x_3,\\&  \notag\hskip-.2in    \sin x_4 \sin x_1 \cos x_2 \cos x_3, \sin 2x_4 \sin x_3 \cos x_2 +\cos 2x_4 \sin x_2,\\&  \notag  0,  \sin x_4 \cos x_1 \cos x_2 \cos x_3, \cos x_4 \sin x_1 \cos x_2 \cos x_3,\\& \notag  \hskip.4in  \cos 2x_4 \sin x_3 \cos x_2 - \sin 2x_4 \sin x_2\Big).\end{align}
\end{theorem}

\section{Ricci curvature and its applications}

For a Riemannian $n$-manifold $M$ and for each unit vector $X\in T_pM,\, p\in M$, the Ricci curvature Ric$(X)$ of $X$ can be regarded as a 
$\delta$-invariant:
\begin{align}\label{21.1} {\rm Ric}(X)=\tau(p)-\tau(L^\perp_X),\end{align}
where $L^\perp_X$ is the hyperplane of $T_pM$ with $X$ as its hyperplane normal.

For the Ricci curvature we have the following result from \cite{c13}.

\begin{theorem}  Let $x:M^n\to R^m(\epsilon)$ be an isometric immersion of a Riemannian $n$-manifold $M^n$ into a Riemannian space form $R^m(\epsilon)$. Then 
\vskip.04in

{\rm (1)}  For each unit  tangent vector $X\in T_p M^n$, we have
\begin{align} \label{21.2} {\rm Ric}(X)\leq \frac{n^2}{4}H^2+ (n-1)\epsilon. \end{align} 
\vskip.04in

{\rm (2)} If $H(p)=0$, then a unit tangent vector $X$ at $p$ satisfies the equality case of \e{21.2} if
and only if $X$ lies in the relative null space $N_p$ at $p$.
\vskip.04in

{\rm (3)} The equality case of \e{21.2} holds identically for all unit tangent vectors at $p$ if and only if  either $p$ is a totally geodesic point or $n=2$ and $p$ is a totally umbilical point.
\end{theorem}

In particular, we have the following.

\begin{corollary} Let $x:M^n\to \mathbb E^m$ be an isometric immersion of a Riemannian
$n$-manifold $M^n$ in a Euclidean $m$-space with arbitrary codimension. Then we have
\begin{align}\label{21.3}H^2(p)\geq\({{4}\over{n^2}}\)\max_{X} {\rm Ric}(X),\end{align}
where $X$ runs over all unit tangent vectors at $p$.\end{corollary}

\begin{remark} There exist many examples of submanifolds in a Euclidean $m$-space which satisfy the equality case of \e{21.3} identically. Two simple  examples are spherical hypercylinder $S^{2}(r)\times \mathbb R$ and  round hypercone in $\mathbb E^{4}$.\end{remark}

For  a  Riemannian $n$-manifold $M$, denote by $K(\pi)$ the sectional curvature  of a 2-plane section $\pi\subset T_pM$, $p\in M$. Suppose $L^k$ is a $k$-plane section of $T_pM$ and $X$ a unit vector in $L^k$. We choose  an orthonormal basis  $\{e_1,\ldots,e_k\}$ of $L^k$ such that $e_1=X$. Define the Ricci curvature $Ric_{L^k}$ of ${L^k}$ at $X$ by
$${\rm Ric}_{L^k}(X)=K_{12}+\cdots+K_{1k},$$
where $K_{ij}$ is the sectional curvature of the $2$-plane section spanned by $e_i, e_j$.   We simply called such a curvature a
{\it $k$-Ricci curvature. } 

For each integer $k,\,2\leq k\leq n$, we have defined in \cite{c13} a Riemannian invariant  on a Riemannian $n$-manifold $M^n$ by
$$\theta_k(p)=\({1\over {k-1}}\)\inf_{L^k,X}{\rm Ric}_{L^k}(X),\quad p\in M^n,$$ where $L^k$
runs over all $k$-plane sections at $p$ and $X$ runs over all unit vectors in $L^k$.

The following result was also proved in \cite{c13}.

\begin{theorem}   Let $x:M^n\to R^m(\epsilon)$ be an isometric immersion of a Riemannian $n$-manifold $M^n$ in a real space form $R^m(\epsilon)$ of constant sectional curvature $c$. Then, for any integer $k,\,2\leq k\leq n$, we have
\begin{align} \label{21.4}H^2(p)\geq {{4(n-1)}\over{n^2}}\({{\theta_k(p)}\over{k-1}} -\epsilon\).
\end{align} 
 
The equality case of \e{21.4} holds identically for all unit tangent vectors at $p$ if and only if either $p$ is a totally geodesic point or $k=n=2$ and $p$ is a totally umbilical point. \end{theorem}

\begin{remark} In general,  given an integer $k,\, 2\leq k\leq n-1$, there {\it does not} exist a positive constant, say $C(n,k)$, such that 
\begin{align} \label{21.5}H^2(p)\geq C(n,k)\max_{L^k,X}{\rm Ric}_{L^k}(X), \end{align}
where $L^k$ runs over all $k$-plane sections in  $T_pM^n$ and $X$ runs over all unit tangent vectors in $L^k$. This fact can be seen from the following example:

Let $x:M^3\to \mathbb E^4$ be a minimal hypersurface whose shape operator is non-singular at some point $p\in M^3$.  Then by the minimality there exist two principal directions at $p$, say
$e_1,e_2$, such that their corresponding principal curvatures $\kappa_1,\kappa_2$ are of the same sign. This implies that the sectional curvature $K_{12}$ at $p$ is positive. Now, consider the minimal hypersurface in $\mathbb E^{n+1}$ which is given by the product of $x:M^3\to \mathbb E^4$ and the identity map $\iota:\mathbb E^{n-3}\to \mathbb E^{n-3}$. 

It is clear that, for any integer $k,\,2\leq k\leq n-1$, the maximum value of the $k$-th Ricci curvatures of $M^n:=M^3\times \mathbb E^{n-3}$ at a point $(p,q),\,  q\in \mathbb E^{n-3}$ is given by $K_{12}=\kappa_1\kappa_2$ which is positive. Since $H=0$, this shows that there does not exist any positive constant $C(n,k)$ which satisfies \e{21.5}.\end{remark}

\section[K\"ahlerian $\delta$-invariants]{K\"ahlerian $\delta$-invariants and  applications to complex geometry}

Let $M^n$ be a K\"ahler manifold of complex dimension $n$. Denote by $J$ the  complex structure on the K\"ahler manifold. For each plane section $\pi\subset T_xM,\, x\in M$, we denote by $K(\pi)$ the sectional curvature of the plane section $\pi$ as before. A plane section $\pi\subset T_xM$ is called {\it totally real\/} if $J\pi$ is perpendicular to $\pi$. 

For a  (real) $2n$-dimensional K\"ahler submanifold $M$ of a K\"ahlerian space form $\tilde
M^m(4\epsilon)$ with constant holomorphic sectional curvature $4\epsilon$,  the scalar curvature $\tau$ of $M$ satisfies \begin{align}\label{22.1}\tau\leq2n(n+1)\epsilon,\end{align}
with equality holding if and only if $M$ is a totally geodesic K\"ahler submanifold.

\subsection{K\"ahlerian  $\delta$-invariants $\delta^c(2n_1,\ldots,2n_k)$}
Let $M$ be a real $2n$-dimensional K\"ahler manifold.  The {\it K\"ahlerian  $\delta$-invariants\/} $\delta^c(2n_1,\ldots,2n_k)$ is defined in \cite{c12} as: 
$$\delta^c(2n_1,\ldots,2n_k)= \tau-\inf\{ \tau(L^c_1)+\cdots+\tau(L^c_k)\} $$ for  each $k$-tuple
$(2n_1,\ldots, 2n_k)$ $\in {\mathcal S}(2n)$,   where $L_1^c,\ldots,L_k^c$ run over all $k$ mutually orthogonal complex subspaces of $T_pM,$ $p\in M$, with real dimensions $2n_1,\ldots,2n_k$, respectively.

For  a K\"ahler submanifold in a complex space form, we have the following general result from \cite{c12}.

\begin{proposition}  Let $M$ be a (real) $2n$-dimensional K\"ahler submanifold of a
complex space form $\tilde M^m(4\epsilon)$. Then, for each $k$-tuple
$(2n_1,\ldots, 2n_k)\in {\mathcal S}(2n)$, the complex $\delta$-invariant $\delta^c(2n_1,\ldots,2n_k)$ satisfies
\begin{align}\label{22.2}\delta^c(2n_1,\ldots,2n_k)\leq 2\Bigg(n(n+1)-\sum_{j=1}^k
n_j(n_j+1)\Bigg)\epsilon.\end{align}

The equality case of inequality \e{22.2} holds at a point $p\in M$ if and only if, there exists an  orthonormal basis 
\begin{align}\notag \begin{array}{c} e_1,\ldots,e_{n_1},Je_1,\ldots,Je_{n_1},\ldots,
e_{2(n_1+\cdots+n_{k-1})+1},\ldots, e_{2(n_1+\cdots+n_{k_1})+n_k},\\ 
Je_{2(n_1+\cdots+n_{k-1})+1},\ldots, Je_{2(n_1+\cdots+n_{k-1})+n_k},e_{2n+1}, \ldots,e_{2m}\end{array}\end{align}
 at $p$, such that  the shape operators of $M$ in $\tilde M^m(4\epsilon)$ at $p$ take the following form:
\begin{align}\label{06.3}A_r=\left( \begin{array}{llllll} A^r_{1} & \cdots & 0&0 &\cdots & 0
\\\vdots  & \ddots& \vdots&\vdots  & \ddots&\vdots \\ 0 &\cdots &A^r_k&0 &\cdots & 0 
\\ 0&\cdots&0&0 &\cdots&0\\ \vdots  & \ddots & \vdots &\vdots &\ddots&\vdots \\
0  &\cdots& 0&0 &\cdots & 0 \end{array}
\right),
\end{align}
$$ r=2n+1,\ldots,2m,$$
where each $A^r_j$ is a symmetric $(2n_j)\times (2n_j)$  submatrix with zero trace.\end{proposition}

\begin{remark} Contrast to inequality \e{22.1},  besides totally geodesic  K\"ahler  submanifolds there do exist K\"ahler submanifolds of $\tilde M^m(4\epsilon)$ which satisfy the equality case of \e{22.2}  identically. 

For instance, let $Q_n$ denote the complex hyperquadric in $ CP^{n+1}(4)$ defined by
\begin{align}\label{22.5}\{(z_0,z_1,\ldots,z_{n+1})\in CP^{n+1}(4):
z_0^2+z_1^2+\cdots+z_{n+1}^2=0\}.\end{align}
For $Q_n$ we have
\begin{align}\label{22.6}\tau=2n^2,\quad \delta^c_k:=\delta^c(2,\ldots,2)=2n(n-1),\end{align} where $2$ in $\delta^c(2,\ldots,2)$ repeats $n$ times. 

The complex quadric $Q_n$ in $CP^{n+1}(4\epsilon)$ satisfies the equality case of \e{22.2} identically for the $n$-tuple $(2,\ldots,2)\in {\mathcal S}(2n)$.

Also, direct products of  K\"ahler submanifolds of complex Euclidean spaces are K\"ahler submanifolds of complex Euclidean spaces which satisfy the equality case of \e{22.2}  for some suitable $k$-tuples.

In views of the above facts, it is an interesting problem to classify all K\"ahler submanifolds of
K\"ahlerian space forms which satisfy either the equality case of inequality \e{22.2}) or the
equality case of inequality \e{22.4}.
\end{remark}

\subsection{Totally real $\delta$-invariants $\delta^r(n_1,\ldots,n_k)$}
For  $(n_1,\ldots, n_k)$ in ${\mathcal S}(2n)$ we also introduced in \cite{c12}  the {\it  totally real $\delta$-invariants\/} $\delta^r(n_1,\ldots,n_k)$ by
\begin{align}\notag&\delta^r(n_1,\ldots,n_k)=\tau-\inf\{ \tau(L^r_1)+\cdots+\tau(L^r_k)\},\end{align}
where $L_1^r,\ldots,L_k^r$ run over all $k$ mutually orthogonal totally real subspaces of $T_pM$, $p\in M$, with dimensions $n_1,\ldots,n_k$, respectively.

For totally real $\delta$-invariants $\delta^r(n_1,\ldots,n_k)$ of a K\"ahler submanifold in a complex space form, we have the following.

\begin{proposition}   Let $M$ be a (real) $2n$-dimensional K\"ahler submanifold of a complex space form $\tilde M^m(4\epsilon)$. Then, for each $k$-tuple $(n_1,\ldots, n_k)\in {\mathcal S}(2n)$, the totally real $\delta$-invariant $\delta^r(n_1,\ldots,n_k)$
satisfies
\begin{align}\label{22.4}\delta^r(n_1,\ldots,n_k)\leq \Bigg(2n(n+1)-\frac 1 2\sum_{j=1}^k
n_j(n_j+1)\Bigg)\epsilon.\end{align}
The equality case of inequality \e{22.4} holds at a point $p\in M$ if and only if, there exists an  orthonormal basis  $e_1,\ldots,e_{2n},e_{2n+1},\ldots,e_{2m}$ at $p$, such that $${\rm Span} \{e_1,\ldots,e_{n_1}\},\ldots,{\rm Span}\{e_{n_1+\ldots+n_{k-1}+1},\ldots, e_{n_1+\ldots+n_{k}}\}$$ are totally real subspaces of $T_pM$ and the shape operators of $M$ in $\tilde M^m(4\epsilon)$ at $p$ take the following form:
$$A_r=\left( \begin{array}{llllll} A^r_{1} & \cdots & 0&0 &\cdots & 0
\\\vdots  & \ddots& \vdots&\vdots  &\ddots &\vdots \\0 &\cdots &A^r_k&0 &\cdots & 0 
\\ 0&\cdots&0&0&\cdots&0\\\vdots  & \ddots & \vdots &\vdots &\ddots&\vdots \\
0  &\cdots& 0&0 &\cdots & 0 \end{array}\right),$$
$$ r=n+1,\ldots,m,$$ where each $A^r_j$ is a symmetric $n_j\times n_j$  submatrix with zero trace. \end{proposition}

\subsection{K\"ahlerian $\delta$-invariant $\,\delta^r_k$ and strongly minimal K\"ahler submanifolds}
For each real number $k$ we may also define a K\"ahlerian $\delta$-invariant $\,\delta^r_k$ by \cite{c24}:
\begin{align}\delta^r_k(p)=\tau(p)-k\,\inf\, K^r(p),\quad p\in M,\end{align}
 where $$\inf K^r(p)=\inf_{\pi^r}\{K(\pi^r)\}$$ and $\pi^r$ runs over all totally real plane sections in $ T_pM$. 

Let $M$ be a  K\"ahler  submanifold of a K\"ahler manifold $\tilde M^{n+p}$.  Just like in the real case,  we denote by $\,h\,$ and $\,A\,$ the second fundamental form and the shape operator of $M^n$ in  $\tilde M^{n+p}$, respectively, 
 
For the  K\"ahler  submanifold we consider an orthonormal frame
$$e_1,\ldots,e_n,e_{1^*}=Je_1,\ldots, e_{n^*} =Je_n$$ of the tangent bundle and an orthonormal frame
$$\xi_1,\ldots,\xi_p,\,\xi_{1^*}=J\xi_1,\ldots, \xi_{p^*}=J\xi_p$$
of the normal bundle.  With respect to such an orthonormal frame, the complex structure $J$ on $M$ is given by 
\begin{align}\label{22.8} J=\begin{pmatrix} 0 &-I_n\\ \\ I_n &0\end{pmatrix}, \end{align}
where $I_n$ denotes an identity matrix of degree $n$.

For a K\"ahler submanifold $M^n$ in  $\tilde M^{n+p}$ the shape operator of $M^n$ satisfies (see, for instance, \cite{O1974})
\begin{align}\label{22.9} A_{J\xi_r}=JA_r,\quad JA_r=-A_r J,\quad\hbox{\rm for }\;
r=1,\ldots,n,1^*,\ldots,p^*,  \end{align}  where $A_r=A_{\xi_r}$. From these it  follows that the shape operator of $M^n$ takes the  form:
\begin{align}  \label{22.10}A_\alpha=\begin{pmatrix}  A_\alpha' &A_\alpha^{''}\\ \\
A_\alpha^{''} &-A_\alpha^{'}\end{pmatrix},\quad  A_{\alpha^*}=\begin{pmatrix} 
-A_\alpha'' &A_\alpha^{'}\\ \\ A_\alpha^{'} &A_\alpha^{''}\end{pmatrix},\quad
\alpha=1,\ldots,p,\end{align}  where $A_\alpha^{'}$ and $A_\alpha^{''}$ are $n\times n$ matrices. 
This condition implies that every K\"ahler submanifold $M^n$ is minimal, i.e., 
$$\hbox{\rm trace}\,A_\alpha=  \hbox{\rm trace}\,A_{\alpha^*}=0,\quad \alpha=1,\ldots,p.$$

The notion of strongly minimal K\"ahler submanifolds was first  introduced in \cite{c24}.

\begin{definition} A  K\"ahler  submanifold $M^n$ of a K\"ahler manifold  $\tilde M^{n+p}$  is called {\it strongly minimal\/} if it satisfies
\begin{align} \hbox{\rm trace}\;A'_\alpha=\,\hbox{\rm  trace}\;A''_{\alpha}=0,\quad \hbox{\rm for }
\;\alpha=1,\ldots,p,\end{align}  
 with respect to some orthonormal frame: $$e_1,\ldots,e_n,e_{1^*}=Je_1,\ldots,\;e_{n^*}=Je_n,\; \xi_1,\ldots,\xi_p,\,\xi_{1^*}=J\xi_1,\ldots, \xi_{p^*}=J\xi_p.$$ \end{definition}

For  K\"ahler  submanifolds,  we have the following sharp general results.
 
\begin{theorem} \label{T:22.1} {\rm \cite{c24}} For any  K\"ahler  submanifold $M^n$ of complex dimension
$n\geq 2$ in a complex space form $\tilde M^{n+p}(4\epsilon)$, the following statements hold.
\vskip.04in
 
{\rm (1)}  For each $k\in (-\infty,4\,],\; \delta^r_k$ satisfies
\begin{align}\label{22.12} \delta^r_k\leq (2n^2+2n-k)\epsilon.\end{align}
 
\vskip.04in
 
{\rm (2)}  Inequality \e{22.12}  fails for every $\,k>4$.

\vskip.04in
 
{\rm (3)}   $\delta^r_k= (2n^2+2n-k)\epsilon$ holds identically for some $k\in (-\infty,4)$ if and only if $M^n$ is  a totally geodesic  K\"ahler  submanifold of $\tilde M^{n+p}(4\epsilon)$.

\vskip.04in
 
{\rm (4)}  The  K\"ahler  submanifold $M^n$ satisfies $\delta^r_4= (2n^2+2n-4)\epsilon$ at a point $x\in M^n$  if and only if there exists an orthonormal basis 
\begin{align}\notag  e_1,\ldots,e_n,e_{1^*}=Je_1,\ldots, e_{n^*}=Je_n, \xi_1,\ldots,\xi_p,\xi_{1^*}=J\xi_1 ,\ldots,\xi_{p^*}=J\xi_p\end{align}
of $\,T_x\tilde M^{n+p}(4\epsilon)$ such that, with respect to this basis, the shape operator of $M^n$ takes the following form:
\begin{align} & A_\alpha=\begin{pmatrix}  A_\alpha' &A_\alpha^{''}\\ \\ A_\alpha^{''} &-A_\alpha^{'}\end{pmatrix},\quad   A_{\alpha^*}=\begin{pmatrix}  -A_\alpha'' &A_\alpha^{'}\\ \\ A_\alpha^{'} &A_\alpha^{''}\end{pmatrix},\\ \notag\\ & \font\b=cmr8 scaled \magstep2
\def\bigzerol{\smash{\hbox{ 0}}} \def\bigzerou{\smash{\lower.0ex\hbox{\b 0}}}
A_\alpha^{'}= \begin{pmatrix}\begin{matrix} a_\alpha & b_\alpha
 \\ 
b_\alpha  &-a_\alpha \end{matrix}& \bigzerou
\\ \bigzerou & \bigzerou \end{pmatrix} ,\quad \font\b=cmr8 scaled \magstep2
\def\bigzerol{\smash{\hbox{ 0}}}
\def\bigzerou{\smash{\lower.0ex\hbox{\b 0}}}
A_\alpha^{''}=\begin{pmatrix}\begin{matrix}a^*_{\alpha} & b^*_{\alpha}
 \\  b^*_{\alpha}  &-a^*_{\alpha} \end{matrix}& \bigzerou
\\  \bigzerou & \bigzerou \end{pmatrix}\end{align} for
some $n\times n$ matrices $A_{\alpha}^{'}, A_{\alpha}^{''},\,\alpha=1,\ldots,p$.
 \end{theorem}
 
\begin{theorem} \label{T:22.2} {\rm \cite{c24}} A complete  K\"ahler  submanifold $M^n$ $(n\geq 2)$
in  $CP^{n+p}(4\epsilon)$ satisfies
\begin{align}\label{22.16}\delta^r_4= 2(n^2+n-2)\epsilon\end{align}   identically if and only if 
  \begin{enumerate} 
  
  \item $M^n$ is a totally geodesic K\"ahler submanifold, or
 
   \item $n=2$ and $M^2$ is a strongly minimal K\"ahler surface in $CP^{2+p}(4\epsilon)$.
\end{enumerate}\end{theorem}

\begin{theorem} \label{T:22.3} A complete K\"ahler submanifold $M^n$ $(n\geq 2)$ of  $\,{\bf C}^{n+p}$ satisfies $\delta^r_4=0$ identically if and only if 
 \vskip.04in
 
{\rm (1)} $M^n$ is  a complex $n$-plane of $\,{\bf C}^{n+p}$, or  
 \vskip.04in
 
 {\rm (2)} $M^n$ is a complex cylinder over  a strongly minimal K\"ahler surface  $M^2$ in $\,{\bf C}^{n+p}$ $(\,$i.e., $M$ is the  product submanifold of  a strongly minimal K\"ahler surface  $M^2$ in ${\bf C}^{p+2}$ and the identity map of the complex Euclidean $(n-2)$-space ${\bf C}^{n-2}\,)$. \end{theorem}
 
\subsection{Examples of strongly minimal K\"ahler submanifolds}   Every totally geodesic  K\"ahler  submanifold of a complex space form is trivially  strongly minimal.  There also exist nontrivial examples of strongly minimal  K\"ahler  submanifolds.

\begin{example} Consider the complex quadric $Q_2$ in $CP^3(4\epsilon)$ defined by 
\begin{align}Q_2=\Big\{(z_0,z_1,z_2,z_3)\in CP^3(4\epsilon): z_0^2+z_1^2+z_2^2+z_3^2=0\Big\},\end{align}
where $\{z_0,z_1,z_2\}$ is a homogeneous coordinate system of $CP^3(4\epsilon)$.

It is known that the scalar curvature $\tau$ of $Q_2$ is equal to $8\epsilon$ and $\inf K^r=0$. Thus, we obtain $\delta^r_4=8\epsilon$. Hence, $Q_2$ is a non-totally geodesic K\"ahler submanifold which satisfies \e{22.16} with $n=2$. Therefore, according to Theorem \ref{T:22.2}, $Q_2$ is a strongly minimal K\"ahler surface in $CP^3(4\epsilon)$.

On the other hand, it is also well-known that $Q_2$ is an Einstein-K\"ahler surface with Ricci tensor $S=4\epsilon g$, where $g$ is the metric tensor of $Q_2$. Thus, the equation of Gauss  yields
\begin{align}g(A_1^2X,Y)=cg(X,Y),\quad X,Y\in TQ_2.\end{align}
 Hence, with respect to a suitable choice of $e_1,e_2,$ $Je_1,Je_2,\xi_1,J\xi_1$, we have
\begin{align} A_1=\begin{pmatrix} A_1' &A_1^{''}\\ \\ A_1^{''} &-A_1^{'}\end{pmatrix},\quad 
 A_{1^*}=\begin{pmatrix} -A_1'' &A_1^{'}\\ \\ A_1^{'} &A_1^{''}\end{pmatrix},\end{align} where
\begin{align} & A_1^{'}=\begin{pmatrix} \sqrt{\epsilon} & 0
 \\  0  &-\sqrt{\epsilon}  \end{pmatrix},\quad A_1^{''}=\begin{pmatrix}  0 &0 \\  0  &0 \end{pmatrix} .\end{align}
This also shows that $Q_2$ is strongly minimal in $CP^3(4\epsilon)$.
\end{example}

Two non-trivial examples of strongly minimal  K\"ahler  surfaces in ${\bf C}^3$  are the following.

\begin{example} The complex surface:
\begin{align}\Big\{z\in {\bf C}^3\,:\,z^2_1+z_2^2+z_3^2=1\Big\}\end{align}
 is a strongly minimal complex surface in {\bf C}$^3$.
\end{example}

\begin{example}   \cite{su3} The  K\"ahler  surfaces:
\begin{align}N^2_k=\Big\{z\in {\bf C}^3\,:\,z_1+z_2+z_3^2=k\Big\},\quad  k\in {\bf C}, \end{align}
 are strongly minimal  K\"ahler  surfaces in ${\bf C}^3$.
\end{example} 

\subsection{Framed-Einsteinian  and strongly minimality}
 A Riemannian $n$-manifold $M$ is called {\it framed-Einstein} if there exist a function $\gamma$ and an orthonormal
frame $\{e_1,\ldots,e_n\}$ on $M$ such that the Ricci tensor $Ric$ of $M$ satisfies
$${\rm Ric}(e_i,e_i)=\gamma g(e_i,e_i)$$ for  $i=1,\ldots,n$. Clearly, every Einstein manifold is framed-Einstein, but not the converse in general.

The following result provides a simple relationship between strongly minimal surfaces and framed-Einsteinian.

\begin{proposition}  Let $M^2$ be a strongly minimal K\"ahler surface in a complex space form. Then $M^2$ is a framed-Einstein K\"ahler surface. \end{proposition}

  \subsection{Invariants $\delta^r_{\ell,k}$}
For a K\"ahler manifold $M$ of complex  dimension $n$, one may extend  the invariant $\delta^r_k$ to $\delta^r_{\ell,k}$ as
\begin{align}\delta^r_{\ell,k}(x)=\tau(x)-\frac{k}{\ell-1}\, \inf_{L^r_\ell}\, \tau(L^r_\ell),\quad x\in M,\end{align}
where $L^r_\ell$  runs over all totally real $\ell$-subspaces of $T_xM$. 

For each integer $\ell \in [2,n]$, the inequality  \eqref{22.12} was extended by B.  Suceav\u{a} \cite{su3} to the following inequality:
\begin{align}\label{su}\delta^r_{\ell,k}(x)\leq  \left\{2n^2+2n-\frac{k}{4}{\ell \choose 2}\right\}\epsilon \end{align}
for  K\"ahler  submanifolds in a complex space form $\tilde M^{m}(4\epsilon)$.   However, for $\ell \geq 3$, the equality sign of \eqref{su} occurs only for totally geodesic  K\"ahler  submanifolds.

\subsection{$\delta(2)$ and K\"ahler hypersurfaces}

For K\"ahler hypersurfaces in ${\bf C}^{n+1}$, Z.  Sent\"urk and L.  Verstraelen \cite{SV} proved the following.

\begin{proposition} \label{P:21.3} Let $M^n$ be 
a K\"ahler hypersurface $M^n$ in ${\bf C}^{n+1}$. The we have

\vskip.04in
{\rm (1)}   $\delta (2)$  satisfies  $\delta (2)\leq 2$;

\vskip.04in
{\rm (2)}  At any point of $M^n$,  $\delta (2) = 0$ holds if and only if (real) ${\rm rank} \,(A) \leq 2$ at that point. \end{proposition}

It follows from Proposition \ref{P:21.3} and Abe's complex version \cite{abe} of the Hartman-Nirenberg cylinder theorem that

\begin{theorem} {\rm \cite{SV}} The complex hypercylinders $\mathcal C^n$ in ${\bf C}^{n+1}$, i.e. the products of
any complex curve $C$ in a complex $2$-plane ${\bf C}^2$ in ${\bf C}^{n+1}$  with complex $(n-1)$-dimensional complex linear subspaces ${\bf C}^{n-1}$ of ${\bf C}^{n+1}$ which are perpendicular to the plane ${\bf C}^{2}$ of the curve $C$, are the complete K\"ahler hypersurfaces $M^n$ in ${\bf C}^{n+1}$ for which the  $\delta (2)$ curvature vanishes identically.
\end{theorem}

\section[$\delta^{\#}$-invariants]{Applications  to affine differential geometry (I) : $\delta^{\#}$-invariants}

\subsection{Basics of affine differential geometry} If $M$ is an $n$-dimensional manifold, let $f:M\to {\bf R}^{n+1}$ be a non-degenerate hypersurface of the affine $(n+1)$-space whose position vector field is nowhere tangent to $M$. Then $f$ can be regarded as a transversal field along itself. We call $\xi=-f$ the centroaffine normal. Following Nomizu, we call $f$ together with this normalization a centroaffine hypersurface.

The centroaffine structure equations are given by
\begin{align} \label{23.1} &D_Xf_*(Y)=f_*(\nabla_X Y)+h(X,Y)\xi,\\
\label{23.2} &D_X\xi=-f_*(X),\end{align} where $D$ denotes the canonical flat connection of ${\bf R}^{n+1}$, $\nabla$ is a torsion-free connection on $M$, called the induced centroaffine
connection, and $h$ is a non-degenerate symmetric $(0,2)$-tensor field, called the {\it centroaffine metric}.

  From now on we assume that the centroaffine hypersurface is definite, i.e., $h$ is definite.  In case that $h$ is negative definite, we shall replace $\xi=-f$ by $\xi=f$ for the affine normal. In this way, the second fundamental form $h$ is always positive definite. In both cases, \eqref{23.1}  holds. Equation \eqref{23.2}  change sign. In case $\xi=-f$, we call $M$ positive definite; in case $\xi=f$, we call $M$ negative definite.

Denote by $\hat\nabla$ the Levi-Civita connection of $h$ and by $\hat R$ and $\hat \kappa$ the curvature tensor and the normalized scalar curvature of $h$, respectively. The {\it difference tensor} $K$ is then defined by
\begin{align} \label{23.3} &K_XY=K(X,Y)=\nabla_X Y-\hat\nabla_X Y,\end{align}
which is a symmetric $(1,2)$-tensor field. The difference tensor $K$ and the cubic form $C$ are related by $$C(X,Y,Z)=-2h(K_XY,Z).$$ Thus, for each $X$, $K_X$ is self-adjoint with respect to $h$.

The Tchebychev form $T$ and the {\it Tchebychev vector field} $T^\#$ of $M$ are defined respectively by
\begin{align} \label{23.4} & T(X)=\frac{1}{n} {\rm trace}\,K_X, \\
\label{23.5} &h(T^\#,X)=T(X).\end{align}
If $T=0$ and if we consider $M$ as a hypersurface of the equiaffine space, then $M$ is a so-called  {\it proper affine hypersphere} centered at the origin.

 If the difference tensor $K$ vanishes, then $M$ is a quadric, centered at the origin, in particular an ellipsoid if $M$ is positive definite and a two-sheeted hyperboloid if $M$ is negative definite.

  An affine  hypersurface $\phi: M\to \bf R^{n+1}$ is called a {\it graph hypersurface} if the transversal vector field  $\xi$  is a constant  vector field.   A result of \cite{NP} states that a graph hypersurface $M$ is  locally affine equivalent to the graph immersion of a certain function $F$.  
Again in case that $h$ is nondegenerate, it defines a semi-Riemannian metric,  called the {\it Calabi metric} of the graph hypersurface. If $T=0$,  a graph hypersurface is a so-called  {\it improper affine hypersphere}.

Let $M_1$ and $M_2$ be two improper affine hyperspheres in ${\bf R}^{p+1}$ and ${\bf R}^{q+1}$ defined respectively by the equations: $$x_{p+1}=F_1(x_1,\dots,x_p),\quad y_{q+1}=F_2(y_1,\dots,y_q).$$ Then one can define a new improper affine hypersphere $M$ in ${\bf R}^{p+q+1}$ by
$$z=F_1(x_1,\dots,x_p) + F_2(y_1,\dots,y_q),$$
where $(x_1,\dots,x_p,y_1,\dots,y_q,z)$ are the coordinates on ${\bf R}^{p+q+1}$. The Calabi normal of $M$ is $(0,\dots,0,1)$. Obviously, the Calabi metric on $M$ is the product metric. Following \cite{D2} we call this composition the {\it Calabi
composition} of $M_1$ and $M_2$.

\subsection{Affine $\delta$-invariants and general fundamental inequalities} Analogous to the $\delta$-invariants $\delta(n_1,\ldots,n_k)$  we  defined in  \cite{cdv} affine $\delta$-invariants $\delta^\#(n_1,\ldots,n_k)$ as follows:
\begin{align} \label{affine delta} \begin{array}{l}
\delta^\#(n_1,\ldots,n_k)(p)=\hat \tau(p)-\sup\{ \hat \tau(L_1)+\cdots+\hat \tau (L_k)\},
\end{array} \end{align}
where $L_1,\ldots,L_k$ run over all $k$ mutually $h$-orthogonal
subspaces of $T_pM$ such that  $\dim L_j=n_j,\, j=1,\ldots,k$ and $\hat \tau (L)$ is the scalar curvature of $L$ {\it with respect to induced affine metric} $h$.

We have the following theorem from \cite{cdv}.

\begin{theorem} \label{T:23.1} Let $M$ be a definite centroaffine hypersurface in ${\bf R}^{n+1}$. Then,  for each $k$-tuple  $(n_1,\ldots,n_k)\in {\mathcal S}(n)$, we have
\begin{equation}\begin{aligned}\label{23.7}\delta^\#(&n_1,\ldots,n_k)\geq
\frac{1}{2}n(n-1)\varepsilon- \frac{1}{2}\sum_{j=1}^k
n_j(n_j-1)\varepsilon \\& \hskip.3in -\frac{n^2\(n+k-1-\sum_{j=1}^k n_j\)}{2\(n+k-\sum_{j=1}^k
n_j\)}h(T^\#,T^\#),\end{aligned}\end{equation}
where $\varepsilon=1$ or $-1$, according to $M$ is positive definite or
negative definite.

The equality case of inequality \eqref{23.7} holds at a point $p\in M$ if and only if $T^\#=0$ at $p$ and there exists an orthonormal basis $\{e_1,\ldots,e_n\}$ at $p$ such that with
respect this basis every linear map  $K_X,\, X\in T_pM,$   takes the following form:
\begin{align}\label{equality}\font\b=cmr8 scaled \magstep2 \def\bigzerol{\smash{\hbox{ 0}}}
\def\bigzerou{\smash{\lower.0ex\hbox{\b 0}}}
K_X= \begin{pmatrix} A^X_{1} &  & &\hskip-.1in \bigzerou \\   & \ddots& & \\
   & &A^X_k & \\ \bigzerou & & &\hskip-.1in  0   \end{pmatrix},\end{align}
where  $\{A^X_j\}_{j=1}^k$ are symmetric $n_j\times n_j$
submatrices satisfying
\begin{align}{\rm trace}\,(A^X_1)=\cdots={\rm trace}\,(A^X_k)=0.\end{align}
\end{theorem}

\subsection{Some immediate applications}
Two immediate consequences of Theorem \ref{T:23.1} are the following.

\begin{corollary} \label{C:23.1} {\rm \cite{cdv}} Let $M$ be a Riemannian $n$-manifold. If there is a $k$-tuple $(n_1,\ldots,n_k)$ in ${\mathcal S}(n)$ such that
\begin{align}\label{23.10}\delta^\#(n_1,\ldots,n_k)<\frac{1}{2}n(n-1)-\frac{1}{2}\sum_{j=1}^k
n_j(n_j-1)\end{align} at some point, then $M$ cannot be realized as an elliptic proper affine hypersphere.

  In particular, if there is a $k$-tuple $(n_1,\ldots,n_k)$ such that
$\delta^\#(n_1,\ldots,n_k)\leq 0$ at some point in $M$, then $M$ cannot be realized as an elliptic proper affine hypersphere in ${\bf R}^{n+1}$.
\end{corollary}

\begin{corollary} {\rm \cite{cdv}} Let $M$ be a Riemannian $n$-manifold. If there is a $k$-tuple $(n_1,\ldots,n_k)$ in ${\mathcal S}(n)$ such that
\begin{align}\label{23.11}\delta^\#(n_1,\ldots,n_k)<\frac{1}{2}\sum_{j=1}^k n_j(n_j-1)-\frac{1}{2}n(n-1)\end{align}
at some point, then $M$ cannot be realized as a hyperbolic proper affine hypersphere in ${\bf
R}^{n+1}$.\end{corollary}

Recall that a hyperovaloid in ${\bf R}^{n+1}$ is a compact strictly convex hypersurface embedded in ${\bf R}^{n+1}$. The inequalities  \eqref{23.7} give rise to some global
centroaffine curvature invariants for hyperovaloids. Moreover, they provide simple characterizations of hyperellipsoids in terms of these global invariants.

\begin{corollary} {\rm \cite{cdv}} Consider a centroaffine hyperovaloid $f:M\to {\bf R}^{n+1},n\geq 3$, with normalization as in \eqref{2.1} with $\xi=-f$. Then, for any $(n_1,\ldots,n_k)\in
{\mathcal S}(n)$,  we have the following global inequality:
\begin{equation}\begin{aligned}\label{23.12}  \int_M
\Bigg( \delta^\#(&n_1,\ldots,n_k)+\text {\small $\frac{n^2\big(n+k-1-\sum_{j=1}^k n_j\big)}{2\big(n+k-\sum_{j=1}^k n_j\big)} $} h(T^\#,T^\#)\Bigg)\omega(h) \\&\geq  \hskip.1in
\frac{1}{2}\Bigg(n(n-1)- \sum_{j=1}^k n_j(n_j-1) \Bigg){\rm vol}(M,h),
\end{aligned}\end{equation}
where $\omega(h)$ is the volume form and ${\rm vol}(M,h)$ is the volume of $(M,h)$.

The equality sign of \eqref{23.12} holds if and only if $M$ is a hyperellipsoid  centered at the origin. \end{corollary}

Moreover, we have the following theorems from \cite{cdv}.

\begin{theorem}\label{T:23.3} {\rm \cite{cdv}} If $\, M$ is a definite centroaffine hypersurface in ${\bf R}^{n+1},$ $n\geq 3,$ which satisfies the equality case of one of the inequalities in
\eqref{23.7}, then we have $\dim\, ({\rm Im}\, K)$ $<n$. \end{theorem} \vskip.01in

\begin{theorem} {\rm \cite{cdv}} \label{T:23.4} Let $M$ be a positive definite centroaffine hypersurface in ${\bf R}^{n+1}, \, n\geq 3,$  satisfying the equality case of \eqref{23.7} for some $k$-tuple $(n_1,\ldots,n_k)\in{\mathcal S}(n)$. If $\dim\, ({\rm Im}\, K)$ is constant on $M$, then $M$ is
foliated by $q$-dimensional ellipsoids centered at the origin of ${\bf R}^{n+1}$, where $q=n-\dim\, ({\rm Im}\, K)$.
\end{theorem}

\begin{theorem} {\rm \cite{cdv}} Let $M$ be a negative definite centroaffine hypersurface in ${\bf R}^{n+1}, \, n\geq 3,$  satisfying the equality case of \eqref{23.7} for a $k$-tuple $(n_1,\ldots,n_k)\in{\mathcal S}(n)$. If $\dim\, ({\rm Im}\, K)$ is constant, then $M$ is foliated by $q$-dimensional two-sheeted hyperboloids centered at the origin, where $q=n-\dim\, ({\rm Im}\, K)$.
\end{theorem}

Similarly, for graph hypersurfaces we have

\begin{theorem} \label{T:23.5}{\rm \cite{DV2}} Let $M$ be a definite graph hypersurface in ${\bf R}^{n+1}$. Then,  for each $k$-tuple  $(n_1,\ldots,n_k)\in {\mathcal S}(n)$, we have
\begin{equation}\begin{aligned}\label{23.13}\delta^\#(&n_1,\ldots,n_k)\geq
-\frac{n^2\(n+k-1-\sum_{j=1}^k n_j\)}{2\(n+k-\sum_{j=1}^k
n_j\)}h(T^\#,T^\#),\end{aligned}\end{equation}
where $\varepsilon=1$ or $-1$, according to $M$ is positive definite or
negative definite.

The equality case of inequality \eqref{23.13} holds at a point $p\in M$ if and only if $T^\#=0$ at $p$ and there exists an orthonormal basis $\{e_1,\ldots,e_n\}$ at $p$ such that with
respect this basis every linear map  $K_X,\, X\in T_pM,$   takes the following form:
\begin{align}\font\b=cmr8 scaled \magstep2 \def\bigzerol{\smash{\hbox{ 0}}}
\def\bigzerou{\smash{\lower.0ex\hbox{\b 0}}} K_X= \begin{pmatrix} A^X_{1} &  & &\hskip-.1in \bigzerou \\   & \ddots& & \\   & &A^X_k & \\ \bigzerou & & &\hskip-.1in  \bigzerou
   \end{pmatrix},\end{align}
where  $\{A^X_j\}_{j=1}^k$ are symmetric $n_j\times n_j$
submatrices satisfying
\begin{align}{\rm trace}\,(A^X_1)=\cdots={\rm trace}\,(A^X_k)=0.\end{align}
\end{theorem}

An immediate consequence of Theorem \ref{T:23.5} is the following.

\begin{corollary} {\rm \cite{DV2}}  If $(M, h)$ is a Riemannian manifold and for some $k$-tuple $(n_1,\ldots,n_k)\in {\mathcal S}(n)$ the $\delta$-invariant satisfies  $\delta^\#(n_1,\ldots,n_k)<0$ at some point, then $(M, h)$ cannot be
realized as improper affine sphere in some affine space.\end{corollary}

\subsection{Remarks and examples}
\begin{remark} For each $k$-tuple $(n_1,\ldots,n_k)$ in ${\mathcal S}(n)$, the inequality \eqref{23.7} is sharp. For instance,  for each $k$-tuple $(n_1,\ldots,n_k)$, any ellipsoid  in ${\bf R}^{n+1}$ centered at the origin satisfies $K=0$ identically. Hence, it satisfies the equality case of the inequality \eqref{23.7} trivially.
\end{remark}

\begin{remark} For any $k$-tuple $(n_1,n_2,\dots,n_k)\in \mathcal S(n)$, there also exist many definite centroaffine hypersurfaces  in ${\bf R}^{n+1}$ with nontrivial $K$ which satisfy the equality case of the inequality \eqref{23.7} identically for the $k$-tuple $(n_1,n_2,\dots,n_k)$.
\end{remark}

\begin{example} {\rm Let $\phi_1:M^{n_1}\to {\bf R}^{n_1+1}\times \{0\}$ and $\phi_2:M^{n_2}\to \{0\}\times {\bf R}^{n_2+1}$ be two elliptic affine hyperspheres centered at the origin.  We put $n=n_1+n_2+1$ and consider the immersion: $\phi^+:M=M_1\times M_2\times {\bf R}\to {\bf R}^{n+1}$ defined by
\begin{align}\label{23.16} \phi^+(u_1,u_2,t)=(\cos t)\phi_1(u_1)+ (\sin  t)\phi_2(u_2),\;\; \end{align} for $u_1\in M_1,\, u_2\in M_2,\, t\in {\bf R}.$ Then a straightforward long computation shows that $\phi^+$ is again an elliptic affine hypersphere, centered at the origin, and that the linear map $K_X, X\in TM$, of $\phi^+$ satisfies
\begin{align}\label{23.17}\font\b=cmr7 scaled \magstep2
\def\bigzerol{\smash{\hbox{ 0}}} \def\bigzerou{\smash{\lower.0ex\hbox{\b 0}}}
K_{X_1}= \begin{pmatrix}\; K^1_{X_{1}} & \bigzerou  & \bigzerou
\; \\ \\   \bigzerou &  \bigzerou &  \bigzerou\\\\  \bigzerou    & \bigzerou & 0
   \end{pmatrix},\;\;\; K_{X_2}= \begin{pmatrix}\;  \bigzerou &
\bigzerou  & \bigzerou\; \\
\\   \bigzerou &  K^2_{X_{2}} &  \bigzerou\\\\  \bigzerou
   & \bigzerou & 0  \end{pmatrix},\;\; K_{V}=0
\end{align}
for $X_1,X_2,V$ tangent to $M_1,M_2,{\bf R}$, respectively, where $K^1_{X_{1}}$ and $K^2_{X_{2}}$ denote the
corresponding linear maps of $\phi_1$ and $\phi_2$, respectively, which satisfy $\hbox{trace}\(K^1_{X_1}\)
   =\hbox{trace}\(K^2_{X_2}\)=0$. Therefore, the immersion $\phi^+$ gives rise to a positive definite centroaffine hypersurface  in
${\bf R}^{n+1}$  which satisfies the equality case of the inequality \eqref{23.7} for the 2-tuple $(n_1,n_2)\in
\mathcal S(n)$ identically.

In particular, if we choose $\phi_2$ to be an ellipsoid centered at the origin, then  $K^2=0$. Hence, the immersion $\phi^+$ also
satisfies the  inequality \eqref{23.7} for the 1-tuple $(n_1)\in \mathcal S(n)$ identically. 

More general, we can use the same procedure again many times and construct in this way a positive definite centroaffine
hypersurface in ${\bf R}^{n+1}$  which satisfies the equality case of the inequality \eqref{23.7} for any $k$-tuple $(n_1,n_2,\dots,n_k)$ in $ \mathcal S(n)$ identically, at least if $n-(n_1+n_2+\dots +n_k)\geq k-1$.} \end{example}

\begin{remark} Similarly, if we choose $\phi_1$ to be a hyperbolic affine hypersphere centered at the origin, and $\phi_2$ an elliptic affine hypersphere centered at the origin, then the corresponding immersion  defined by
\begin{align}\notag \phi^-(u_1,u_2,t)=(\cosh t)\phi_1(u_1)+ (\sinh  t)\phi_2(u_2),\;\;
u_1\in M_1,\; u_2\in M_2,\; t\in {\bf R}\end{align} gives rise to
a negative definite centroaffine hypersurface  in ${\bf R}^{n+1}$
which satisfies the equality case of the inequality \eqref{23.7} for the 2-tuple $(n_1,n_2)\in \mathcal S(n)$. 

Also here we can use the same procedure again many times and construct in this way a negative definite centroaffine
hypersurface in ${\bf R}^{n+1}$  which satisfies the equality case of the inequality \eqref{23.7} for any $k$-tuple
$(n_1,n_2,\dots,n_k)\in \mathcal S(n)$ identically, at least if $n-(n_1+n_2+\dots +n_k)\geq k-1$.
\end{remark}

\section[Warped products]{Applications  to affine differential geometry (II): warped products}

\subsection{A realization problem}  For a Riemannian $n$-manifold $(M,g)$ with Levi-Civita connection $\nabla$, \'E. Cartan and A. P. Norden studied nondegenerate affine immersions $f:(M,\nabla)\to {\bf R}^{n+1}$ with a transversal vector field $\xi$ and with $\nabla$ as the induced connection. 
  
  The well-known Cartan-Norden theorem states that if $f$ is  such an affine immersion, then  either $\nabla$ is flat and $f$ is a graph immersion or $\nabla$ is not flat and ${\bf R}^{n+1}$ admits a parallel Riemannian metric relative to which $f$ is an isometric immersion and $\xi$ is perpendicular to $f(M)$ (cf. for instance, \cite[p. 159]{N}) (see, also \cite{DNV}).
  
In \cite{c38,c44},  we study Riemannian manifolds in affine geometry from a view point different from Cartan-Norden. More precisely, we investigate  the following:  
  \vskip.1in
  
{\bf Realization Problem:}  {\it Which  Riemannian manifolds $(M,g)$ can be immersed  as affine
hypersurfaces in an affine space in such a way that the 
fundamental  form $h$ (e.g. induced by the centroaffine normalization or a constant 
transversal vector field) is the given Riemannian metric $g$}?
\vskip.1in

We say that a Riemannian manifold $(M,g)$ can be {\it realized as an affine hypersurface} if there exists a codimension one affine immersion from $M$ into some affine space in such a way that the induced affine metric $h$ is  exactly the Riemannian  metric $g$ of $M$ (notice that we do not put any assumption on the affine connection).  

In \cite{c44} we prove that Robertson-Walker space times can be be realized as centro-affine and graph hypersurfaces in some affine space.

\subsection{Existence results}
Also, we  show in \cite{c38} that  there exist many warped product Riemannian manifolds  which can be realized either as graph  or   centroaffine hypersurfaces. More precisely, we prove the following.

\begin{theorem}  \label{T:24.1}  Let  $f=f(s)$ be a positive function defined on an open interval $I$. Assume that ${\bf R}, S^n(a^2)$, $H^n(-a^2)$, and $\mathbb E^n$ are equipped with their canonical metrics. Then we have:
\vskip.04in

{\rm (a)} Every warped product surface $I\times_f {\bf R}$ can  be realized as a graph surface in the affine $3$-space ${\bf R}^3$.
\vskip.04in

{\rm (b)}  For each integer $n>2$, the warped product manifold  $I\times_f H^{n-1}(-a^2)$ can  be realized as a graph hypersurface
in ${\bf R}^{n+1}$. 
\vskip.04in

{\rm (c)} If $f'(s)\ne 0$ on $I$, then  the warped product manifold  $I\times_f {\mathbb E}^{n-1},n>2,$  can be realized as a graph hypersurface in ${\bf R}^{n+1}$. 
\vskip.04in

{\rm (d)} If $f'(s)^2>a^2$ on $I$ for some positive number $a$, then the warped product manifold  $I\times_f S^{n-1}(a^2),n>2,$ can  be realized as a graph hypersurface in ${\bf R}^{n+1}$. 
  \end{theorem}
  
  \begin{theorem}  \label{T:24.2} The following results hold.
\vskip.04in

{\rm (a)} If $n>2$ and  $f=f(s)$ is a positive function defined on an open interval $I$, then we have:
\vskip.04in

 {\rm (a.1)}  If $f'(s)^2>f^2(s)-a^2$ on $I$ for some positive number $a$, then $I\times_f H^{n-1}(-a^2)$ can  be realized as a centroaffine hypersurface
in ${\bf R}^{n+1}$. 
\vskip.04in

 {\rm (a.2)} If $f'(s)^2>f(s)^2$ on $I$, then   $I\times_f {\mathbb E}^{n-1}$  can be realized as a centroaffine hypersurface
in ${\bf R}^{n+1}$. 

\vskip.04in 
 {\rm (a.3)} If $f'(s)^2>f(s)^2+a^2$ on $I$ for some positive number $a$, then  $I\times_f S^{n-1}(a^2)$ can  be realized as a graph hypersurface in ${\bf R}^{n+1}$. 
\vskip.04in 

{\rm (b)} If $n=2$ and $f=f(s)$ is a  positive function defined on a closed interval $[\alpha,\beta]$, then the warped product surface $J\times _f {\bf R}, J=(\alpha,\beta)$, can always be realized as a centroaffine surface in ${\bf R}^3$.   
 \end{theorem}

\subsection{An inequality for graph hypersurfaces and its application}
We apply in \cite{c38} an affine $\delta$-invariant and prove  the following results.

\begin{theorem} \label{T:24.3} {\rm \cite{c38}} If a warped product manifold  $N_1\times_f N_2$  can be realized as a graph hypersurface
in  ${\bf R}^{n+1}$, then the warping function satisfies
\begin{align} \label{24.1}\frac{ \Delta f}{f}\geq -\frac{(n_1+n_2)^2}{4n_2}h(T^{\#},T^\#),\end{align} 
 where   $n=n_1+n_2$,  $n_1=\dim N_1$ and $n_2=\dim N_2$.\end{theorem}

The following result characterizes affine hypersurfaces which verify the equality case of inequality \e{24.1}.

\begin{theorem}  \label{T:24.4} {\rm \cite{c38}} Let  $\phi:N_1\times_f N_2\to {\bf R}^{n+1}$ be a realization of a warped product manifold as  a  graph  hypersurface. If the warping function satisfies the equality case of \eqref{24.1} identically,  then we have: 
\vskip.04in

{\rm (a)} The Tchebychev vector field $T^\#$ vanishes identically.
\vskip.04in

{\rm (b)} The warping function $f$ is a harmonic function.
\vskip.04in

{\rm (c)}  $N_1\times_f N_2$ is realized as an improper affine hypersphere.
  \end{theorem} 

An application of Theorem \ref{T:24.3} is the following.

 \begin{corollary} \label{C:24.1}  If $N_1$ is a compact Riemannian manifold, then every  warped product manifold  $N_1\times_f N_2$ cannot be realized
as an improper affine hypersphere in  ${\bf R}^{n+1}$. \end{corollary}

As an application of Theorems \ref{T:24.3} and \ref{T:24.4} we have the following. 

 \begin{theorem}  \label{T:24.5} {\rm \cite{c38}} If the Calabi metric of an improper affine hypersphere in an affine space  is the Riemannian product metric of  $k$ Riemannian manifolds,  then the improper affine hypersphere  is locally the Calabi composition of $k$ improper affine spheres.
 \end{theorem} 
 
 Theorem \ref{T:24.3} also implies the following.

\begin{corollary} \label{C:24.2}   If the warping function $f$ of a warped product manifold  $N_1\times_f N_2$ satisfies $\Delta f < 0$  at some point on $N_1$, then $N_1\times_f N_2$ cannot be realized
as an improper affine  hypersphere in  ${\bf R}^{n+1}$. \end{corollary}

\subsection{An inequality for centro-affine hypersurfaces and its application}
 Similarly, for centro-affine hypersurfaces we have the following \cite{c38}.

 \begin{theorem}  \label{T:24.6}  If a warped product manifold  $N_1\times_f N_2$ can be realized as a centroaffine hypersurface
in ${\bf R}^{n+1}$, then the warping function satisfies
\begin{align} \label{24.2}\frac{ \Delta f}{f}\geq  n_1\varepsilon- \frac{(n_1+n_2)^2}{4n_2}h(T^{\#},T^\#),\end{align} where $n=n_1+n_2$, $n_i=\dim N_i$, $i=1,2$,  $ \Delta$ is the Laplace operator of $N_1$, and $\varepsilon=1$ or $-1$ according to whether the centroaffine hypersurface is elliptic or hyperbolic. 
\end{theorem}

\begin{theorem}   \label{T:24.7} Let  $\phi:N_1\times_f N_2\to {\bf R}^{n+1}$ be a realization of  a warped product manifold $N_1\times_f N_2$ as  a centroaffine hypersurface. If the warping function satisfies the equality case of \eqref{24.2} identically,  then we have: 
\vskip.04in

{\rm (1)} The Tchebychev vector field $T^\#$ vanishes identically.
\vskip.04in

{\rm (2)} The warping function $f$ is an eigenfunction of the Laplacian $\Delta$ with eigenvalue $n_1\varepsilon$.
\vskip.04in

{\rm (3)}  $N_1\times_f N_2$ is realized as a proper affine hypersphere centered at the origin.  \end{theorem} 

Two immediate consequences of Theorem \ref{T:24.6} are the following.

\begin{corollary} \label{C:24.3}  If the warping function $f$ of a warped product manifold  $N_1\times_f N_2$ satisfies $\Delta f \leq 0$  at some point on $N_1$, then $N_1\times_f N_2$ cannot be realized
as an elliptic  proper affine hypersphere in ${\bf R}^{n+1}$. \end{corollary}
 
\begin{corollary} \label{C:24.4}  If the warping function $f$ of a warped product manifold  $N_1\times_f N_2$ satisfies  $(\Delta f)/ f< -\dim N_1$  at some point on $N_1$, then $N_1\times_f N_2$ cannot be realized
as a hyperbolic  proper affine hypersphere in   ${\bf R}^{n+1}$. \end{corollary}

Another interesting application of Theorem \ref{T:24.6} is the following.

 \begin{corollary} \label{C:24.5}  If $N_1$ is a compact Riemannian manifold, then every  warped product manifold  $N_1\times_f N_2$ with arbitrary warping function cannot be realized
as an elliptic  proper affine hypersphere in  ${\bf R}^{n+1}$. \end{corollary}

Another application of Theorem \ref{T:24.6} is the following.

 \begin{corollary} \label{C:24.6}  If $N_1$ is a compact Riemannian manifold, then every  warped product manifold  $N_1\times_f N_2$ cannot be realized
as an improper affine hypersphere in an affine space ${\bf R}^{n+1}$. \end{corollary}

\subsection{Remarks and examples}

The following examples  show that the results of this section are optimal.

\begin{example}\label{E:24.1} {\rm Let $M=N_1\times_{\cos s} N_2$  be the warped product of the open interval  $ N_1=(-\pi,\pi)$ and an open portion $N_2$ of the unit $(n-1)$-sphere $S^{n-1}(1)$ equipped with the warped product metric: 
\begin{equation}\begin{aligned}\label{24.3} &h={\rm d}s^2+\cos^2 s\Bigg({\rm d}u_2^2+\cos^2u_2{\rm d}u_3^2+\cdots+\prod_{j=2}^{n-1}\cos^2 u_j{\rm d}u_n^2\Bigg).
\end{aligned}\end{equation}
Consider the immersion of $M$ into the affine $(n+1)$-space  defined by
\begin{equation}\begin{aligned}\label{24.4} &\Bigg(\sin s,\sin u_2\cos s,\ldots,  \sin u_n \cos s\prod_{j=2}^{n-1}\cos^2 u_j ,\cos s\prod_{j=2}^{n}\cos u_j\Bigg). \end{aligned}\end{equation}
Then $M$ is a  centroaffine elliptic   hypersurface whose  centroaffine metric is the warped product metric  \eqref{24.3} and it satisfies $T^\#=0$. Moreover, the warping function $f=\cos s$ satisfies $$\frac{\Delta f}{f}=1=\varepsilon n_1.$$ Hence, this centroaffine hypersurface satisfies the equality case of \eqref{24.2} identically. Consequently, the estimate given in Theorem \ref{T:24.6} is optimal for centroaffine elliptic  hypersurfaces.}
\end{example}

\begin{example}\label{E:24.2} {\rm Let $M={\bf R}\times_{\cosh s} H^{n-1}(-1)$  be  the  warped product of the real line and the unit hyperbolic space $H^{n-1}(-1)$ equipped with warped product metric: 
\begin{equation}\begin{aligned}\label{24.5} &h={\rm d}s^2+\cosh^2 s \Bigg({\rm d}u_2^2+\cosh^2u_2{\rm d}u_3^2+\cdots+\prod_{j=2}^{n-1}\cosh^2 u_j{\rm d}
u_n^2\Bigg).
\end{aligned}\end{equation}
Consider the immersion of $M$ into the affine $(n+1)$-space  defined by
\begin{equation}\begin{aligned}\notag&\Bigg(\sinh s,\sinh u_2\cosh s,\ldots, \sinh u_n \cosh s\prod_{j=2}^{n-1}\cosh^2 u_j ,\cosh s\prod_{j=2}^{n}\cosh u_j\Bigg).\end{aligned}\end{equation}
Then $M$ is a centroaffine hyperbolic  hypersurface whose  centroaffine metric is the warped product metric  \eqref{24.5} and it
satisfies $T^\#=0$. Moreover, the warping function $f=\cosh s$ satisfies $$\frac{\Delta f}{f}=-1=\varepsilon n_1.$$ Therefore, this centroaffine hypersurface satisfies the equality case of \eqref{24.2} identically. Consequently, the estimate given in Theorem \ref{T:24.6} is optimal for centroaffine  hyperbolic hypersurfaces as well.}
\end{example}

\begin{example}\label{E:24.3} {\rm Let $M={\bf R}\times_{s} N_2$  be  the warped product of the real line and an open portion $N_2$ of  $S^{n-1}(1)$ equipped with the warped product metric: 
\begin{equation}\begin{aligned}\label{24.6} &h={\rm d}s^2+ s^2 \Bigg({\rm d}u_2^2+\cos^2u_2{\rm d}u_3^2+\cdots+\prod_{j=2}^{n-1}\cos^2 u_j{\rm d}
u_n^2\Bigg).
\end{aligned}\end{equation}
Consider the immersion of $M$ into the affine $(n+1)$-space  defined by
\begin{equation}\begin{aligned}\label{24.7} & s \Bigg( \sin u_2, \sin u_3\cos u_2,\ldots, \sin u_n \prod_{j=2}^{n-1}\cos^2 u_j , \prod_{j=2}^{n}\cos u_j, \frac{s}{2}\Bigg).
\end{aligned}\end{equation}
Then $M$ is a graph hypersurface with Calabi normal  given by $\xi=(0,\ldots,0,1)$ and it satisfies $T^\#=0$. Moreover, the  Calabi metric of this graph hypersurface is given by the warped product metric  \eqref{24.6}. Clearly, the warping function  is a harmonic function. Therefore, this warped product graph hypersurface satisfies the equality case of \eqref{24.1} identically. Consequently, the estimate given in Theorem \ref{T:24.3} is also optimal.}
\end{example}

\begin{remark} {\rm Example \ref{E:24.1} shows that the conditions $$\Delta f\leq 0$$ in Corollary \ref{C:24.3}  and  the ``harmonicity''   in Corollary  \ref{C:24.4}  are both necessary. }
\end{remark}

\begin{remark} {\rm Example \ref{E:24.1}  implies that the condition ``$N_1$ is a compact Riemannian manifold'' given   in Corollary \ref{C:24.6}  is  necessary.}
\end{remark}

\begin{remark} {\rm Example \ref{E:24.2} illustrates that the condition $$(\Delta f)/f< -\dim N_1$$   given in Corollary \ref{C:24.5}  is sharp.}
\end{remark}

\begin{remark} {\rm Example \ref{E:24.3} shows that the condition $$\Delta f< 0$$ in Corollary \ref{C:24.1}  is optimal as well.}
\end{remark}

\section[Eigenvalues]{Applications to affine differential geometry (III):  eigenvalues}

 For each integer $k\in [2,n]$, we introduce in \cite{c43}  the affine invariant
$\hat \theta_k$ on  $M$ by 
\begin{align} \label{25.1} \hat \theta_k(p) =\text{$\({1\over {k-1}}\)$} \sup_{L^k,X}
\hat S_{L^k}(X),\quad p\in T_pM,\end{align} where  $L^k$
runs over all linear $k$-subspaces in the tangent space $T_pM$ at $p$ and $X$ runs
 over all $h$-unit vectors in $L^k$. 

The {\it relative $K$-null space} ${\mathcal N}_p^K$ of $M$ in ${\bf R}^{n+1}$ is defined by
\begin{align} \label{25.2}{\mathcal N}_p^K=\Big\{X\in T_pM:K(X,Y)=0\;\hbox{\rm for all
}Y\in T_pM\Big\}.\end{align}

When $\dim {\mathcal N}^K_p$ is constant,  $ {\mathcal N}^K=\cup_{p\in M}  {\mathcal N}^K_p$ defines a subbundle of the tangent bundle, called  the {\it relative $K$-null subbundle.}
 
 For affine  hypersurfaces in ${\bf R}^{n+1}$ we have the following results. 

\begin{theorem}\label{T:25.1}{\rm \cite{c43}}  Let $f:M\to {\bf R}^{n+1}$ be a  locally strongly convex  centroaffine hypersurface in ${\bf R}^{n+1}$. Then, for any  integer $k \in [2,n]$, we have: 

\vskip.04in 

{\rm (1)} If $\hat\theta_k \ne \varepsilon$ at a point $p\in M$, then every eigenvalue of $\,K_{T^\#}$  at $p$ is greater than $\({{n-1}\over n}\)(\varepsilon -\r)$.
\vskip.04in 

{\rm (2)}  If $\r=\varepsilon$, every eigenvalue of $K_{T^\#}$ at $p$ is $\geq 0$.
\vskip.04in 

{\rm (3)}  A nonzero vector $X\in T_pM$ is an eigenvector of the operator $K_{T^{\#}}$  with eigenvalue 
$\({{n-1}\over n}\)(\varepsilon -\r )$ if and only if $\r=\varepsilon$ and $X$ lies in the relative $K$-null space ${\mathcal N}^K_p$ at $p$,
\vskip.04in

\noindent where $\varepsilon=1$ or $-1$ according to $M$ being of elliptic or hyperbolic type.
\end{theorem}

 \begin{theorem}\label{T:25.2}{\rm \cite{c43}}  Let $f:M\to {\bf R}^{n+1}$ be a graph hypersurface in ${\bf R}^{n+1}$ with positive definite Calabi metric. Then, for any integer $k\in [2,n]$, we have: 
\vskip.04in 

{\rm (1)}  If $\hat \theta_k\ne 0$ at a point $p\in M$, then every eigenvalue of $K_{T^\#}$  at $p$ is greater than $\({{1-n}\over n}\)\r$.
\vskip.04in 

{\rm (2)}  If $\hat \theta_k=0$ at  $p$, then every eigenvalue of $K_{T^\#}$ at $p$ is $\geq 0$.
\vskip.04in 

{\rm (3)}  A nonzero vector $X\in T_pM$ is an eigenvector of the operator $K_{T^{\#}}$  with eigenvalue 
$\({{1-n}\over n}\)\r $ if and only if we have $\hat \theta_k(p)=0$  and $X\in {\mathcal N}^K_p$.
\end{theorem}

Examples were given in \cite{c43} to show that the estimates of the eigenvalues given in Theorems \ref{25.1} and \ref{T:25.2} are both sharp.

As applications of Theorems \ref{25.1} and \ref{T:25.2}  we have:

\begin{corollary} \label{C:25.1}  Let $f:M\to {\bf R}^{n+1}$ be a locally strongly convex centroaffine   hypersurface in ${\bf R}^{n+1}$.  If  $\, \sup \hat K\ne \varepsilon$ at a point $p\in M$, then every eigenvalue of the  operator $K_{T^\#}$ at $p$ is greater than $\(\frac{n-1}{n}\)(\epsilon -\sup \hat K(p))$.
\end{corollary}

\begin{corollary} \label{C:25.2}  Let $f:M\to {\bf R}^{n+1}$ be a locally strongly convex  centroaffine  hypersurface  in ${\bf R}^{n+1}$.  If $\, \sup \hat S\ne \varepsilon$ at a point $p\in M$, then every eigenvalue of the  operator $K_{T^\#}$ at $p$ is greater than $\(\frac{n-1}{n}\)(\epsilon -\sup \hat S(p))$.
\end{corollary}

\begin{corollary} \label{C:25.3}  Let $f:M\to {\bf R}^{n+1}$ be a locally strongly convex  centroaffine hypersurface  in ${\bf R}^{n+1}$. 
If we have $\hat\theta_k<\varepsilon$ on $M$ for some integer $k\in[2,n]$, then   every eigenvalue of $K_{T^\#}$ is positive.\end{corollary}

\begin{corollary} \label{C:25.4}  An  elliptic centroaffine hypersurface $M$  in ${\bf R}^{n+1}$ is centroaffinely equivalent to an open portion of a hyperellipsoid if and only if we have $$n K_{T^\#}= (n-1)(1 -\hat\theta_k)I$$  on $M$ for some integer $k\in [2,n]$.\end{corollary}

\begin{corollary} \label{C:25.5} A hyperbolic centroaffine  hypersurface  $M$  in ${\bf R}^{n+1}$  is centroaffinely equivalent to an open portion of a two-sheeted hyperboloid  if and only if,   for some integer $k\in [2,n]$,  we have  $n K_{T^\#}= (1-n)(1 +\hat\theta_k)I$ identically on $M$.  \end{corollary}

\begin{corollary} \label{C:25.6}  Let $f:M\to {\bf R}^{n+1}$  be a graph hypersurface with positive definite Calabi metric. 
If we have either $\sup \hat K\ne 0$ or $\sup \hat S\ne 0$ at  a point $p\in M$, then every eigenvalue of the  operator $K_{T^\#}$ is
greater than $\(\frac{1-n}{n}\)\sup \hat K$  at $p$.
\end{corollary}

\begin{corollary} \label{C:25.7}  Let $f:M\to {\bf R}^{n+1}$  be a graph hypersurface with positive definite Calabi metric. 
If there exists an integer $k\in[2,n]$ such that $\hat\theta_k<0$ holds on $M$, then  every eigenvalue of $K_{T^\#}$ is positive. \end{corollary}

\begin{corollary} \label{C:25.8} Let $M$ be a Riemannian $n$-manifold. If  there exists an integer $k\in[2,n]$ such that $\hat\theta_k(p)<1$  at some point $p\in M$, then $M$ cannot be realized as an elliptic proper affine hypersphere in ${\bf R}^{n+1}$.\end{corollary}

\begin{corollary} \label{C:25.9} Let $M$ be a Riemannian $n$-manifold. If  there exists an integer $k\in[2,n]$ such that $\hat\theta_k(p)<-1$  at some point $p\in M$, then $M$ cannot be realized as a hyperbolic proper affine hypersphere in ${\bf R}^{n+1}$.\end{corollary}

\begin{corollary} \label{C:25.10} Let $M$ be a Riemannian $n$-manifold. If there exists an integer $k\in[2,n]$ such that $\hat\theta_k(p)<0$  at some point $p\in M$, then $M$ cannot be realized as an improper affine hypersphere in ${\bf R}^{n+1}$.\end{corollary}

\section{$\delta^{\#}(2)$ and affine hypersurfaces}

From \e{23.7} we see that the simplest  affine $\delta$-invariant different from the affine scalar curvature $\hat \tau$ is $\delta^{\#}(2)$, which is defined by
\begin{align}\label{26.1} \delta^\# (2)=\hat \tau- \sup \hat  K.\end{align}
 In this section we survey some results on affine hypersurfaces involving $\delta^{\#}(2)$.

\subsection{An inequality involving $\delta^\# (2)$} For  $\delta^{\#}(2)$ we have the following general optimal inequality (see Theorem \ref{T:23.1} for the general affine $\delta$-invariants).

\begin{theorem} \label{T:26.1} {\rm\cite{SS}} Let $M$ be an $n$-dimensional definite centroaffine hypersurface in an affine $(n+1)$-space ${\bf R}^{n+1}$. Then we have 
\begin{align}\label{26.2} \delta^{\#} (2)\geq  \frac{\epsilon}{2}(n+1)(n-2)-\frac{n^2(n-2)}{2(n-1)} h(T^{\#},T^{\#}),\end{align}
where $T^{\#}$ is the Tchebychev vector field, $h$ is the induced affine metric on $M$ and $\epsilon=1$ or $-1$ according to $M$ is positive-definite or negative-definite, respectively. \end{theorem}

\begin{theorem} \label{T:26.2} {\rm \cite{SS}} Let $M$ be an $n$-dimensional definite centroaffine hypersurface in an affine $(n+1)$-space ${\bf R}^{n+1}$.  If $M$ realizes the equality in \e{26.2} at a point, then the Tchebychev vector field $T^\#$ vanishes at $p$ and the equiaffine support  function is stationary at $p$. Moreover, if $M$ realizes the equality for every $p\in M$, then $M$ is a proper affine hypersphere centered at the origin. \end{theorem}

 When $n=3$, inequality \e{26.2} reduces to 
\begin{align}\label{26.3} \delta^{\#} (2)\geq 2\epsilon-\frac{9}{4} h(T^{\#},T^{\#}).\end{align}
Of course also for definite affine hyperspheres in equiaffine (Blaschke) differential geometry the above inequality (but with $T=0$) remains valid.

Note that $M$ is called an {\it affine hyperpshere} if and only if the affine shape operator $S$ satisfies
 $S=H I$, i.e.,  all affine normals are parallel or pass through a fixed point. 
 
We recall that as $M$ is definite, we can choose the Blaschke normal such that the affine metric is positive definite. This allows us to make still the following distinction between affine hyperpsheres:

\vskip.04in
{\rm (1)} elliptic affine hyperspheres, i.e. all affine normals pass through a fixed point and $H >0$,

\vskip.04in
{\rm (2)} hyperbolic affine hyperspheres, i.e. all affine normals pass through a fixed point and $H<0$,

\vskip.04in
{\rm (3)} parabolic affine hyperspheres, i.e. all the affine normals are parallel ($H=0$).
\vskip.04in

In \cite{KSV}, Kriele, Scharlach and Vrancken classified 3-dimensional elliptic affine spheres for which the equality in \e{26.3}  is assumed. This classification is achieved through reducing the problem to the problem of classifying those 3-dimensional minimal surfaces in the unit pseudo-sphere $S^5_3(1)$ of index 3 whose ellipses of curvature are circles. 

In \cite{KSV} they also investigated  2-dimensional minimal surfaces in $S^5_3(1)$  whose ellipses of curvature are circles.
Also, 3-dimensional hyperbolic  affine spheres  satisfying the equality in \e{26.3} are classified by Kriele and Vrancken in \cite{KV}. Here, the integrability conditions were solved directly. A crucial point, in both the elliptic and the hyperbolic case, is the existence of an adapted frame such that  the difference tensor $K$ has a special simple form.

\subsection{$\delta^\# (2)$  and affine hypersurfaces admitting a pointwise \\ symmetry}
It is known that if an affine 3-sphere satisfies the equality case of inequality \e{26.3}, then it admits a pointwise $S_3$-symmetry, but this family is bigger.

The idea to study pointwise symmetries comes from Bryant's article \cite{bryant} on special Lagrangian 3-manifolds. 
In \cite{v2}, Vrancken  studied affine hypersurfaces $M$ admitting a pointwise symmetry, i.e., there exists a subgroup $G$ of ${\rm Aut}(T_pM)$ for all $p\in M$ which preserves the affine metric $h$, the difference tensor $K$ and the affine shape operator $S$, i.e., for any $X,Y\in T_pM$ and $g\in G$, we have
\begin{equation}\begin{aligned}\notag &h(gX,gY)=h(X,Y),\;  \\&K(gX,gY)=g(K(X,Y)),\; \\&S(gX)=g(SX).\end{aligned}\end{equation}
 
 L. Vrancken introduces the following symmetric polynomial 
\begin{equation} \notag
f_p(x,y,z)=h(K(x e_1+ye_2 +z e_3,x e_1+y e_2+z e_3),x e_1+y e_2 +z e_3),
\end{equation}
where $\{e_1, e_2, e_3 \}$ is an orthonormal basis at the point $p$. The apolarity condition implies that the trace of  this polynomial  with respect to the metric vanishes. 

As far as such symmetric polynomials with vanishing trace  on a 3-dimensional real vector space are concerned, 
we have the following result by Bryant \cite{bryant}: 
\begin{theorem} Let $p \in M$ and assume that there exist an orientation preserving isometry which preserves $f_p$. Then there
exists an orthonormal basis of $T_pM$ such that either

\vskip.04in
{\rm (i)} $f_p=0$, in this case $f_p$ is preserved by every isometry,

\vskip.04in
{\rm (ii)}  $f_p = \lambda  (2 x^3-3 xy^2-3 xz^2)$, for some positive number $\lambda$ in which case $f_p$ is  preserved by
 a $1$-parameter group of rotations, 

\vskip.04in
{\rm (iii)}  $f_p=6 \lambda  xyz$ for some positive number $\lambda$, in which case $f_p$ is preserved by the discrete group
 $A_4$ of order $12$,

\vskip.04in
{\rm (iv)}  $f_p= \lambda (x^3-3 xy^2)$ for some positive number $\lambda$, in which case $f_p$ is preserved by the discrete
group $S_3$ of order $6$,

\vskip.04in
{\rm (v)}  $f_p =\lambda (2x^3-3x y^2 -3x z^2) + 6 \mu xyz$, for some $\lambda,\mu >0$, with $\lambda\neq \mu$, in which case $f_p$ is 
preserved by the group $\mathbb Z_2$ of order $2$,

\vskip.04in
{\rm (vi)}  $f_p = \lambda (2 x^3 -3 xy^2 -3 xz^2) +\mu (y^3-3 xy^2)$ for some $\lambda, \mu >0$, with $\mu \neq \sqrt{2} \lambda$,  in which 
case $f_p$ is preserved by the group $\mathbb Z_3$.
\end{theorem}   

Vrancken considers the special class of 3-dimensional affine hyperspheres, such that there exists an orientation preserving isometry which preserves $f_p$, at every point $p$. Thus he obtains six different types of expressions for the cubic form and he describes how to construct locally these affine hyperspheres of Type $k$, from (i) to (vi). In particular, Type (i) corresponds to quadrics, Type (iii) to flat affine hyperspheres and Type (iv) realizes the  equality of  \e{26.3} in terms of the affine $\delta$-invariant $\delta^{\#}(2)$.

In \cite{Lus}, Lu and Scharlach extend the above work and solve the corresponding algebraic problem for a general $S$. They determine the nontrivial stabilizers $G$ of the pair $(K, S)$ under the action of ${\rm SO}(3)$ on a Euclidean vector space $(V, h)$ and classify three-dimensional locally strongly convex affine hypersurfaces admitting a pointwise $G$-symmetry. Besides well-known examples, they obtain warped products of two-dimensional affine spheres and curves.

\subsection{An improved inequality involving $\delta^\# (2)$}  Recently, J. Bolton, F. Dillen,  J. Fastenakels and L. Vrancken \cite{BDFV} improve inequality \e{26.2} to the following.

\begin{theorem} \label{T:26.3} Let $M$ be an $n$-dimensional definite centroaffine hypersurface in an affine $(n+1)$-space ${\bf R}^{n+1}$. Then we have 
\begin{align}\label{26.4} \delta^{\#} (2)\geq  \frac{\epsilon}{2}(n+1)(n-2)-\frac{n^2(2n-3)}{2(2n+3)} h(T^{\#},T^{\#}),\end{align}
where  $\epsilon=1$ or $-1$ according to $M$ being positive or negative definite. \end{theorem}

It is proved in \cite{BDFV} that if $n\geq 4$, then the Tchebychev vector field $T$ vanishes. Moreover, they obtain a classification theorem for  definite centroaffine hypersurfaces in  ${\bf R}^{4}$ which satisfy the equality case of  \e{26.4} with $T\ne 0$. More precisely, they obtain the following.

\begin{theorem} \label{T:26.4} Let $M$ be a $3$-dimensional centroaffine hypersurface of ${\bf R}^{4}$  which attains equality at every point in \e{26.4}. Then  $M$ is locally  given by the one of the following immersions:
\vskip.05in 

{\rm (1)}  $$ f =\frac{(3\lambda-b_1) \rme^{-3t}}{\sqrt{b_1^2-9\lambda^2+\epsilon}}V+ \frac{\rme^{{-t}}}{\sqrt{b_1^2-9\lambda^2+\epsilon}}W,$$
where $V$ is a constant vector along the hypersurface, $W$ is a surface for which the Tchebychev form vanishes, and 
$b_1$ and  $\lambda$ are solutions of the following system of ordinary differential equations:
$$\frac{{\rm d}b_1}{{\rm d}t}=\frac{b_1^2-27\lambda^2+\epsilon}{3\lambda},\hskip.2in \frac{{\rm d}\lambda}{{\rm d}t}=-\frac{2}{3}b_1$$
and $b_1^{2}-9\lambda^2+\epsilon\ne 0$.

\vskip.05in

{\rm (2)} $f(t,u,v)=(tu,tv,tg(u,v)+\gamma_2(t),t),$ where $\gamma_{2}$ is a function satisfying 
$$t\gamma'''_{2}\gamma_{2}-t^{2}\gamma'''_{2}\gamma'_{2}-t^{2}(\gamma''_{2})^{2}+4\gamma''_{2}\gamma_{2}-4t\gamma'_{2}\gamma''_{2}=0$$ and  $(u, v) \mapsto (u, v, g(u, v))$ defines an improper affine sphere with affine normal $(0, 0, 1)$.
\end{theorem}

\section[Applications  to general relativity]{Applications of $\delta$-invariants to general relativity}

In 1916, Albert  Einstein and David Hilbert independently constructed the equations for the pure
gravitational field (cf. \cite{log})
$$G_{\lambda\mu}:={\rm Ric}_{\lambda\mu}-\frac{1}{2}g_{\lambda\mu} \rho=\kappa T_{\lambda\mu}$$
where Ric is the trace of the Ricci tensor, $T_{\lambda\mu}$  the energy-momentum tensor and $\kappa$  some
dimension-transposing parameter.

The remaining part of his scientific career Einstein searched for an unification between his description of gravity and the other known forces, in particular electromagnetism (cf. \cite{Goe}).  Several attempts have been made by using e.g., a non-symmetric metric, a connection with torsion, ..., etc.

 Shortly  after the publication of the theory of General Relativity, Theodor Kaluza noticed  in April 1919 that when he solved Einstein's equations for General Relativity using five dimensions, Maxwell's equations for electromagnetism emerged spontaneously. 
Kaluza wrote to Albert Einstein who encouraged him to publish. His theory was published in 1921 in \cite{Ka} with Einstein's support.  

In 1926, Oskar  Klein \cite{Kl} suggested that this fifth dimension would be compactified and unobservable on experimentally accessible energy scales.  However their work was neglected for many years as attention was directed towards quantum mechanics. 

The idea that fundamental forces can be explained by additional dimensions did not re-emerge until string theory was developed. This idea of compactifying the extra dimension has now dominated the search for a unified theory and lead to the 11D supergravity theory and more recently 10D superstring theory. Recently, this strategy of using higher dimensions to unify different forces is also an active area of research in particle physics   (see \cite{ow} for an overview).

Instead of compactifying the extra dimensions other approaches have also been developed during the last decade.  For example,  one particular variant of Kaluza--Klein  theory is space-time-matter  theory or induced matter theory, chiefly promulgated by Paul Wesson and other members of the so-called Space-Time-Matter Consortium (ses, http://astro.uwaterloo.ca/$\backsim$wesson/). 

In this version of the theory, it is noted that solutions to the equation $R_{AB} = 0$
with $R_{AB}$  the 5D Ricci curvature, may be re-expressed so that in four dimensions, these solutions satisfy Einstein's equation: \begin{align}G_{\mu\nu} = 8\pi T_{\mu\nu}\end{align}
with the precise form of the $T_{\mu\nu}$ following from the Ricci-flat condition on the 5D space. 

Since the energy-momentum tensor $T_{\mu\nu}$ is normally understood to be due to concentrations of matter in 4D space, the above result can be interpreted as saying that 4D matter is induced from geometry in a Ricci-flat 5D space.
In particular, the soliton solutions of \begin{align}R_{AB} = 0 \end{align}  can be shown to contain the Robertson-Walker metric in both matter-dominated (early universe) and radiation-dominated (present universe) forms. The general equations can be shown to be sufficiently consistent with classical tests of general relativity to be acceptable on physical principles, while still leaving considerable freedom to also provide interesting cosmological models (see \cite{W1,W2}).

There is another approach proposed by Lisa Randall and Raman Sundrum  in 1999. Randall-Sundrum models imagine  our Universe as a 5D anti de Sitter space and the elementary particles except for the graviton are localized on a (3 + 1)-D brane or branes.  

Their models attempt to address the hierarchy problem between the observed Planck and weak scales by embedding the 3-brane in a warped 5D metric; the warping of the extra dimension is analogous to the warping of space-time in the vicinity of a massive object, such as a black hole   (see  \cite{rs1,rs2} for details). The Randall--Sundrum scenario has gained a lot of support  recently (see \cite{c44}).

More recently, S. Haesen and L. Verstraelen \cite{HV} extend $\delta$-invariants to pseudo-Riemannian manifolds as follows:

 For a given set of mutually orthogonal plane sections $\{ L_j\}$  with dimensions $n_1,\ldots,n_k$ such that $n_1+\cdots+n_k\leq n$, the $\delta$-curvatures in the
semi-Riemannian case are given by
 \begin{align}\notag \Lambda(n_1,\ldots,n_k)=\tau -\inf \{ \tau(L_1)+\cdots+\tau(L_k): L_j\;  \text{\rm  a non-null plane section}\}\end{align}
 and
 \begin{align}\notag \hat \Lambda(n_1,\ldots,n_k)=\tau -\sup \{ \tau(L_1)+\cdots+\tau(L_k): L_j\;  \text{\rm  a non-null plane section}\}.\end{align}
Clearly, when the pseudo-Riemannian manifolds are Riemannian, these reduce to $\delta(n_1,\ldots,n_k)$ and $\hat \delta(n_1,\ldots,n_k)$, respectively.

In \cite{Hae,HV}, Haesen and Verstraelen prove the following.

\begin{theorem} Let $(M,g)$ be an $n$-dimensional Riemannian or Lorentzian manifold. Assume that $M$ is locally and isometrically embedded in a $(n+1)$-dimensional semi-Riemannian manifold
$(N, \tilde g)$ with diagonalizable Ricci tensor $\widetilde {{\rm Ric}},\, ($i.e., there exists an orthonormal basis $\{e_\alpha\}$
of $N$ such that $\widetilde {{\rm Ric}}=\sum_{a=1}^{n+1}\lambda_a e_a\otimes e_a)$. Then for any $k$-tuple $(n_1,\ldots,n_k)$ such that $n_1<m$ and $n_1+\cdots+n_k\leq n$, we have
\begin{align}\notag ||\overrightarrow H||^2_\perp \geq c(n_1,\ldots,n_k)\Lambda(n_1,\ldots,n_k)-\frac{1}{2}c(n_1,\ldots,n_k)\left\{\sum_{\alpha=1}^n \varepsilon_\alpha \lambda_\alpha-\lambda_{n+1}\right\},\end{align}
if {\rm sign}$(N)=(s_M+1,t_M)$, and 
\begin{align}\notag ||\overrightarrow H||^2_\perp \leq c(n_1,\ldots,n_k)\hat \Lambda(n_1,\ldots,n_k)-\frac{1}{2}c(n_1,\ldots,n_k)\left\{\sum_{\alpha=1}^n \varepsilon_\alpha \lambda_\alpha+\lambda_{n+1}\right\},\end{align}
if {\rm sign}$(N)=(s_M,t_M+1)$, where 
$$c(n_1,\ldots,n_k)=\frac{2(n+k-\sum_{j=1}^k n_j)}{n^2(n+k-1-\sum_{j=1}^k n_j)}.$$ 

Equality holds if and only if the second fundamental form
has the following form with respect to the eigenframe of the Ricci tensor $\widetilde {Ric}$:
$$    (\Omega_{\alpha\beta})=\font\b=cmr8 scaled \magstep2 \def\bigzerol{\smash{\hbox{ 0}}} \def\bigzerou{\smash{\lower.0ex\hbox{\b 0}}} A_{e_r}=\left( \begin{matrix} A^r_{1} & \hdots & 0 \\ \vdots  & \ddots& \vdots &\bigzerou \\ 0 &\hdots &A^r_k& \\ \\&\bigzerou & &\mu I_s \end{matrix} \right),
$$where $s=n-\sum_{j=1}^k n_j$ and $ A_{n_j}$ is a symmetric $n_j\times n_j$-matrix with ${\rm trace}\,(A_{n_j}) =\mu$. \end{theorem}

Haesen and  Verstraelen applied these $\delta$-invariants and inequalities  in \cite{HV} to investigate  the embedding problem of space-times from the view point of ideal embedding.  
 Among others, they found that ideally embedded hypersurfaces in a pseudo-Euclidean space contain both the de Sitter space-time and a Robertson-Walker model.
Embeddings of the de Sitter model and Robertson-Walker's  model were already considered by J. Ponce de Le\'on \cite{ponce} in 1988. It was later realized that his 5D embedding space was flat and this was used by S. Seahra and P. Wesson \cite{sw} to study the structure of the Big Bang.

Furthermore, S. Haesen and L. Verstraelen \cite{HV} also applied $\delta$-invariant to find a class of anisotropic perfect fluid models containing a timelike two-surface of constant curvature;  also been shown to be ideally embedded.  (For further results in this direction, see also \cite{Hae}.)

As we mentioned earlier in {\bf 24.1},  the author has recently constructed  explicit embeddings of Robertson-Walker's cosmological models as centroaffine and graph hypersurfaces in some flat affine space so that the induced affine metrics are exactly the space-time metrics on the Robertson-Walker space-times (see \cite{c44} for details). 

Since many physics quantities are represented by the metric and its curvatures, most physical  quantities of the Robertson-Walker space-times are preserved via the realizations constructed in \cite{c44}. 
These embeddings allow us to study Robertson-Walker space-times and their submanifolds using the method of affine differential geometry.  
Furthermore, in contrast to space-time-matter theory in which matter could be either space-like or  time-like depending on embeddings,  if  additional physical quantities  were induced from the extra dimension via our embeddings it would then be neutral (in the sense that it is neither space-like nor time-like).

\section{$k$-Ricci curvature and shape operator}

In this section we explain a sharp relationship between the $k$-Ricci curvature and the shape operator for arbitrary submanifolds in a real space form with arbitrary codimension (see \cite{c13} for details).  

Recall that for a submanifold $M$ in a Riemannian manifold, the {\it relative null space} of $M$ at a point $p\in M$ is defined by
$$N_p=\{X\in T_pM:h(X,Y)=0\;\mbox{ for all }Y\in T_pM\}.$$
Moreover, for each integer $k,\,2\leq k\leq n$, the invariant  $\theta_k$ on a Riemannian manifold is defined by
$$\theta_k(p)=\({1\over {k-1}}\)\inf_{L^k,X}{\rm Ric}_{L^k}(X),\quad p\in T_pM^n,$$ where $L^k$
runs over all $k$-plane sections at $p$ and $X$ runs over all unit vectors in $L^k$.

The following result  from \cite{c13}  provides a sharp relationship between the invariant $\theta_k$ and the shape operator $A_H$ for arbitrary submanifolds, regardless of codimension.
\vskip.1in
\begin{theorem} \label{T:28.1}  Let $x:M\to R^m(\epsilon)$ be an isometric immersion of a
Riemannian $n$-manifold $M$ into a real space form  $R^m(\epsilon)$ of constant sectional curvature $\epsilon$. Then, regardless of codimension, for any integer $k,\,2\leq k\leq n$, and any point $p\in M$ we have:  
\vskip.04in

$(1)$ If $\theta_k(p) \ne \epsilon$, then the shape operator at the mean curvature vector $H$ satisfies
 $$A_H>\frac {n-1} n (\theta_k(p) - \epsilon)\hbox{I}\quad \hbox{at }p,$$ where $\mbox{I}$ denotes the identity map of $T_pM$.
\vskip.04in

$(2)$ If $\theta_k(p) =\epsilon$, then $A_H\geq 0$ at $p$.
\vskip.04in

$(3)$ A unit vector $X\in T_pM$ satisfies $$A_HX={{n-1}\over n}(\theta_k(p) -
\epsilon)X$$ if and only if $\theta_k(p) =\epsilon$ and $X$ lies in the relative null space at $p$.
\vskip.04in

 $(4)$ We have
 $$A_H\equiv {{n-1}\over n}(\theta_k- \epsilon)\hbox{I}$$ at a point $p\in M$ if and only if $p$ is a  totally geodesic point, i.e., the second fundamental form vanishes identically at $p$.\end{theorem}
\vskip.1in

\begin{remark} Clearly, the estimate of $A_H$ given in statement (2) of Theorem \ref{T:28.1} is sharp. Here we provide an example to illustrate that statement (1) of Theorem \ref{T:28.1} is also sharp. 

Consider a hyperellipsoid in $\mathbb E^{n+1}$ defined by
$$ax_1^2+x_2^2+\ldots+x_{n+1}^2=1,$$
where $0<a< 1$. The principal curvatures of the hyperellipsoid are given by
$$a_1={a\over{(1+a(a-1)x_1^2)^{3/2}}},\quad $$
$$a_2=\ldots=a_{n}= {1\over{(1+a(a-1)x_1^2)^{1/2}}}.$$
Therefore, for any $k,\,2\leq k\leq n$, the $k$-Ricci curvatures at a point $p$ satisfies
$${\rm Ric}_{L^k}(X)\geq (k-1)\theta_k(p)  :={{(k-1)a}\over{(1+a(a-1)x_1^2)^2}}>0$$
for any $k$-plane section $L^k$ and any unit vector $X$ in $L^k$. Moreover, the eigenvalues $\kappa_1,\ldots,\kappa_n$ of the shape operator $A_H$ are given by
$$\kappa_1=\ldots=\kappa_{n-1}={{a+(n-1)(1+a(a-1)x_1^2)}
\over  {{n(1+a(a-1)x_1^2)^2}}}, $$ 
$$ \kappa_n={{a(a+(n-1)(1+a(a-1)x_1^2))}\over   {{n(1+a(a-1)x_1^2)^3}}}.
$$
From these it follows that $$A_H>{{n-1}\over n}{\theta_k(p) }  I_n$$ and
$$\kappa_1-{{n-1}\over n}\theta_k(p) =
{{a^2}\over   {{n(1+a(a-1)x_1^2)^3}}}\longrightarrow 0,\;\; {\rm as}\; a\rightarrow 0.$$
\vskip.1in

Hence, this example shows that our estimate of $A_H$ in statement (1) is also sharp. \end{remark}
\vskip.1in

One may also apply Theorem \ref{T:28.1}  to provide a lower bound for every eigenvalue of the shape operator $A_H$ for all isometric immersions of a given Riemannian $n$-manifold, regardless of codimension. As the simplest example, Theorem \ref{T:28.1} implies immediately the following.
\vskip.1in

\begin{corollary}\label{C:28.1} Let $M$ be a Riemannian $n$-manifold. If there is a point $p\in M$ such that every sectional curvature of $M$ at $p$ is equal to $1$, then, for any isometric immersion of $M$ in a Euclidean $m$-space with arbitrary codimension, each eigenvalue of the shape operator $A_H$ at $p$ is greater than ${(n-1)/n}$.\end{corollary}

\begin{remark} The estimate of the eigenvalues of $A_H$ given in Corollary \ref{C:28.1} is sharp. For instance, assume that $M$ is a surface in $\mathbb E^3$ whose two principal curvatures at $p$ are given respectively by $a$ and $1/a$ with $a\geq 1$. Then the smaller  eigenvalue of $A_H$ at $p$ is equal to $(a^2+1)/2a^2$ which approaches to $1/2$ when $a$ approaches $\infty$.  \end{remark}

  For an $n$-dimensional submanifold $M$ in $\mathbb E^m$, let $\mathbb E^{n+1}$ be the linear subspace of dimension $n+1$ spanned by the tangent space at a point $p\in M$ and the mean curvature vector $H(p)$ at $p$. 
  
  Geometrically, the shape operator $A_{n+1}$ of $M$ in $\mathbb E^m$ at $p$ is the shape operator of the orthogonal projection of $M^n$ into $\mathbb E^{n+1}$. Moreover, it is known that if the shape operator of a hypersurface in $\mathbb E^{n+1}$ is definite at a point $p$, then it is strictly convex at $p$.  For this reason a submanifold $M$ in $\mathbb  E^m$ is said to be $H$-{\it strictly convex\/} if the shape operator $A_H$ is positive-definite at each point  in $M$. 

Theorem \ref{T:28.1} implies immediately the following.
\vskip.1in

\begin{corollary}\label{C:28.2}  Let $M$ be a submanifold of a Euclidean space. If there is an integer $k,\,2\leq k\leq n,$ such that $k$-Ricci curvatures of $M$ are positive, then $M$ is $H$-strictly convex, regardless of codimension.\end{corollary}

By applying Theorem \ref{T:28.1}  and a classical result of W. S\"uss we have  the following.

\begin{corollary}\label{C:28.3}  If $M$ is a compact hypersurface of $\mathbb E^{n+1}$ with $\theta_k \geq 0$ (respectively, with $\theta_k>0$) for a fixed 
$k,\,2\leq k\leq n$, then $M$ is embedded as a convex (respectively, strictly convex)
hypersurface in $\mathbb E^{n+1}$. In particular, if $M$ has constant scalar curvature, then $M$ is a hypersphere of $\mathbb E^{n+1}$.\end{corollary}

\begin{corollary}\label{C:28.4}  If $M$ is a compact hypersurface of $\mathbb E^{n+1}$ with nonnegative Ricci curvature (respectively, with positive Ricci curvature),
then $M$ is embedded as a convex (respectively, strictly convex) hypersurface in $\mathbb E^{n+1}$. In particular, if $M$ has constant scalar curvature, then $M$ is a
hypersphere of $\mathbb E^{n+1}$.\end{corollary}

\section{General inequalities (I): $CR$-products}

Several general geometric optimal inequalities for submanifolds in complex space forms (or more generally in K\"ahler manifolds) have been discovered during the last decades. Although these inequalities do not  relate directly with the $\delta$-invariants, however due to their simplicity we present them as well.

\subsection{Segre embedding introduced in 1891} 

Let $(z_0^i,\ldots,z_{N_i}^i)$ $(1\leq i\leq s)$ denote the homogeneous coordinates of $CP^{N_i}$. Define a map:
$$S_{N_1\cdots N_s}:CP^{N_1}(4)\times\cdots\times CP^{N_s}(4)\to CP^N(4),\;\; N=\prod_{i=1}^s (N_i+1)-1,$$
which maps a point $$((z_0^1,\ldots,z_{N_1}^1),\ldots, (z_0^s,\ldots,z_{N_s}^s))$$ of the product K\"ahler manifold
$CP^{N_1}(4)\times\cdots\times CP^{N_s}(4)$ to the point $$(z^1_{i_1}\cdots z^s_{i_j})_{1\leq i_1\leq N_1,\ldots,1\leq i_s\leq N_s}\in CP^N(4).$$ 

The map $S_{N_1\cdots N_s}$, introduced by C. Segre \cite{segre} in 1981, is a K\"ahler embedding
which is known as the {\it Segre embedding}.  The Segre embedding plays some important role in algebraic geometry.

\subsection{``Converse'' of  Segre  embedding -- 1981}  In 1981, Chen \cite{c2}  and  Chen and W. E. Kuan \cite{ckuan} established the ``converse'' of the Segre embedding as  follows:

\begin{theorem} \label{T:29.1} Let $M_1,\ldots,M_s$ be K\"ahler manifolds of dimensions $N_1,\ldots,$ $N_s$, respectively. Then every K\"ahler immersion
$$\phi :M_1\times\cdots\times M_s\to CP^N(4),\quad N=\prod_{i=1}^s (N_i+1)-1,$$ of $M_1\times\cdots\times M_s$ into $CP^N(4)$   is locally  given by the Segre embedding, i.e., $M_1,$ $\ldots,M_s$ are open portions of $CP^{N_1}(4),\ldots, CP^{N_s}(4)$, respectively, and moreover, the K\"ahler immersion $\phi$ is congruent to the Segre embedding. \end{theorem}

\subsection{$CR$-submanifolds}

Let $(\tilde M,g,J)$ be a K\"ahler manifold with complex structure $J$ and  $M$  a Riemannian manifold isometrically immersed
in $\tilde M$. For each tangent vector $X$ of $M$, we put
\begin{align} JX=PX+FX,\end{align}
where $PX$ and $FX$ are the tangential and the normal components of $JX$. Then $P$ is an endomorphism of $TM$ and $F$ is a normal-bundle-valued 1-form.

For each point $x\in N$, we denote by $\mathcal D_x$ the maximal holomorphic subspace of the tangent space $T_xM$ defined by $$\mathcal D_x=T_xM\cap J(T_xM).$$ If the dimension of $\mathcal D_x$ is the same for all $x\in M$, $\mathcal D_x$'s define a  distribution $\mathcal D$ on $M$, which is called the {\it holomorphic distribution} of $M$.  

A submanifold $M$ in a K\"ahler manifold $\tilde M$ is called a $CR$-{\it submanifold} if there exists a  holomorphic distribution $\mathcal D$ on $M$ whose orthogonal complement $\mathcal D^\perp$ is a
{\it totally real distribution}, i.e., $J\mathcal D^\perp\subset T^\perp N$ (cf. \cite{bejancu}).
 If $\dim \mathcal D_x=0$ at each point $x\in M$, then the $CR$-submanifold $M$ is nothing but a totally real submanifold.  

 In \cite{BC}, D. E. Blair and the author proved that every $CR$-submanifold of a Hermitian manifold is a $CR$-manifold.

In 1978,  the author discovered that  the totally real distribution $\Cal D^\perp$ of a $CR$-submanifold of a K\"ahler manifold
is always completely integrable.  This integrability theorem implies that every proper $CR$-submanifold of a K\"ahler manifold is
foliated by totally real submanifolds. By applying this integrability theorem,  A. Bejancu \cite{b1979}  proved in 1979  that a $CR$-submanifold of a K\"ahler manifold is mixed totally geodesic if and only if each leaf of the totally real distribution is totally geodesic in the $CR$-submanifold.

Chen's integrability theorem  was later extended to $CR$-submanifolds of various families of Hermitian manifolds by various
geometers. For instance, this theorem was extended to $CR$-submanifolds of locally conformal symplectic manifolds by Blair and Chen in \cite{BC}. Furthermore, they  constructed in \cite{BC} $CR$-submanifolds in some Hermitian manifolds with non-integrable  totally real distributions.

For a $CR$-submanifold $M$ with Riemannian connection $\nabla$, let $e_1,\ldots,e_{2h}$ be an orthonormal frame field of the holomorphic distribution $\Cal D$. Put $\hat H=\,$trace $\hat \sigma$, where $\hat \sigma(X,Y)=(\nabla_XY)^\perp$ is the component of $\nabla_XY$ in the totally real distribution. The holomorphic distribution  $\Cal D$ is called a  minimal distribution if $\hat H=0$, identically. 

Although the holomorphic distribution is not necessarily integrable in general, the author proved in \cite{c2} that the holomorphic
distribution of a $CR$-submanifold is always a minimal distribution.
He also  discovered in \cite{c1981} a canonical cohomology class $c(M)\in H^{2h}(M;\hbox{\bf R})$ for every compact
$CR$-submanifold $M$ of a K\"ahler manifold. By applying this cohomology class, he proved the following: 

Let $M$ be a compact  $CR$-submanifold of a K\"ahler manifold. If the  cohomology group $H^{2k}(M;\hbox{\bf R})=\{0\}$ for some integer $k\leq h$, then either the holomorphic distribution $\Cal D$ is not integrable or the totally real distribution $\Cal D^\perp$ is not  minimal. 

The cohomology class $c(M)$ was applied  by S. Dragomir  in his study concerning the minimality of Levi distribution (cf. \cite{DO}).

A. Ros \cite{r1983} proved that if $M$ is an $n$-dimensional compact minimal $CR$-submani\-fold of $CP^m(4)$, then the
first nonzero eigenvalue of the Laplacian of $M$ satisfies $$\lambda_1\leq \frac{2}{n}(n^2+4h+p).$$

For further results on $CR$-submanifolds, see for instance \cite{bejancu,c2,c3,c1981,c18,YKon}.
 
 \subsection{$CR$-products}
 Recall that a submanifold of a K\"ahler manifold is called a $CR$-{\it product} if it is the Riemannian product $N_T\times N_\perp$ of a complex submanifold $N_T$ and a totally real submanifold $N_\perp$. 

Let $f:N_\perp\to CP^p(4)$ be a Lagrangian submanifold. Then the composition:
\begin{align} \begin{CD}CP^h(4)\times N_\perp @>\iota \times f >> CP^h(4)\times CP^p(4)  @>S_{hp}>> CP^{h+p+hp}(4)\end{CD}\end{align} 
is a $CR$-product in $CP^{h+p+hp}$, where $\iota :CP^h\to CP^h$ denotes the identity map and $S_{hp}$ is the Segre embedding.

  A $CR$-product $M=N_T\times N_\perp$ in  $CP^m(4)$ is called a {\it standard $CR$-product} if $m=h+p+hp$ and $N_T$ is a totally geodesic K\"ahler submanifold of $CP^m(4)$.
  
  It is known that (see  \cite{c2} for details):

 \begin{theorem} \label{T:29.4} The following statements hold:
\vskip.04in

{\rm (1)}  $CR$-products in complex hyperbolic spaces are either complex  or totally real.
\vskip.04in

{\rm (2)} A submanifold in a complex Euclidean space   is a $CR$-product if and only if it is the direct sum of a complex submanifold and a totally real submanifold of linear complex subspaces.
\vskip.04in

{\rm (3)} $CR$-products in the complex projective space $CP^{h+p+hp}(4)$ are obtained from the Segre embedding in a natural way.
\end{theorem}

\subsection{A general inequality for  submanifolds in complex space forms}

For an arbitrary submanifold $M$ in a complex space form $\tilde M^m(4\epsilon)$, we have the following general inequality.

\begin{theorem} \label{T:29.2} {\rm \cite{clm}} Let $M$ be a submanifold of  a complex space form $\tilde M^m(4\epsilon)$. Then we have
\begin{align} \label{29.2} ||\bar \nabla h||^2\geq 2 \epsilon^2 ||P||^2 ||F||^2,\end{align}
with  equality sign holding if and only if $M$ is a cyclic-parallel $CR$-submanifold. \end{theorem}

Here, by a cyclic-parallel submanifold we mean a submanifold whose second fundamental form $h$ satisfies
\begin{align} \label{29.3} (\bar \nabla_X h)(Y,Z)+(\bar \nabla_Y h)(Z,X)+(\bar \nabla_Z h)(X,Y)=0\end{align}
for any vectors $X,Y,Z$ tangent to $M$. 

It follows from polarization that a submanifold of a real space form is cyclic-parallel if and only if it is a parallel submanifold. For submanifolds in $CH^m(-4)$, the author, G. D.  Ludden and S. Montiel  \cite{clm} showed that a $CR$-submanifold $M$ of $CH^m(-4)$ is cyclic-parallel if and only if the preimage $\pi^{-1}( M)$ of $M$ (via the Hopf fibration $\pi: H_1^{2m+1}(-1)\to CH^m(-4)$) has parallel second fundamental form in  $H_1^{2m+1}(-1)$.  Similar result also
holds for cyclic-parallel $CR$-submanifolds in $CP^m(4)$ (see \cite{YKon}).
For the classification of cyclic-parallel $CR$-submanifolds of complex space forms, see \cite{clm,YKon}.

\subsection{Two inequalities for $CR$-products}
For $CR$-products in complex projective space $CP^m(4)$, we have the following general inequality  involving the squared norm of the second fundamental form.

\begin{theorem} \label{T:29.3} {\rm \cite{c2}} Let $M$ be a $CR$-product in $CP^m(4)$. Then we have
\begin{align} \label{29.4} ||h||^2\geq 4hp,\end{align}
where $h=\dim_{\bf C} \mathcal H_x$ and $p=\dim_{\bf R}\mathcal H^\perp$. 

 If the equality sign of \e{29.4} holds, then $N_T$ and $N_\perp$ are both totally geodesic in $CP^m(4)$. Moreover, the immersion is rigid and the $CR$-product is a standard one. \end{theorem}

For minimal $CR$-products in $CP^m(4)$, we also have the following.

 \begin{theorem}\label{T:29.5} {\rm \cite{c2}} If $M$ is a minimal $CR$-product in $CP^m(4)$, then the scalar curvature  of $M$ satisfies $$\tau\geq
2h^2+2h+\frac{1}{2}(p^2-p),$$ with the equality holding when and only when $||h||^2=4hp$ holds. \end{theorem}

\section{General inequalities (II): $CR$-warped products}

In \cite{c22}, the author proved that there do not exist any warped product submanifold of the form: $N_\perp\times_f N_T$ in any K\"ahler manifold $\tilde M$ such that $N_T$ is a complex submanifold and $N_\perp$ is a totally real submanifold of $\tilde M$. Moreover, he  combined in  \cite{c22} the notion of warped products with the notion of $CR$-submanifolds to introduce  the notion of {\it $CR$-warped products} as follows:

\begin{definition} A submanifold of a  K\"ahler manifold $M$ is called a {\it $CR$-warped product} if it is a warped product
$N_T\times_f N_\perp $ of a complex submanifold $N_T$ and a totally real submanifold $N_\perp$. \end{definition}

\subsection{Optimal inequalities involving $||h||^2$} 

For arbitrary $CR$-warped products in an arbitrary K\"ahler manifold, the author proved the following \cite{c22} .

\begin{theorem} \label{T:30.1}  Let $M=N_T\times_f N_\perp$ be a  $CR$-warped product in an arbitrary K\"ahler manifold $\tilde M$.  Then we have:

\vskip.04in 
{\rm (1)} The squared norm of the second fundamental form of $M$ satisfies
\begin{align} \label{30.1} ||h||^2\geq 2p\,||\nabla (\ln f)||^2,\end{align}
where $\nabla \ln f$ is the gradient of $\,\ln f$ and $p$ is the dimension of  $N_\perp$.

\vskip.04in
{\rm (2)} If the equality sign of \e{30.1} holds identically, then  $N_T$ is a totally geodesic submanifold and $N_\perp$ is a totally
umbilical submanifold of $\tilde M$. Moreover, $M$ is a minimal submanifold in $\tilde M$.

\vskip.04in
{\rm (3)}  If $M$ is anti-holomorphic and $p>1$, then the equality sign of \e{30.1} holds identically if and only if   $N_\perp$ is a totally umbilical submanifold of $\tilde M$.

\vskip.04in
{\rm (4)}  If $M$ is anti-holomorphic and $p=1$, then the equality sign of \e{30.1} holds identically if and only if  the characteristic vector field $J\xi$ of $M$ is a principal vector field with zero as its principal curvature. $($Notice that in this case, $M$ is a real hypersurface in $\tilde M.) $

Also, in this case, the equality sign of \e{30.1} holds identically if and only if $M$ is a minimal hypersurface in $\tilde M$.
 \end{theorem}

$CR$-warped products in complex space forms attaining the equality case of \e{30.1} were classified in \cite{c23}. In fact, we have the following results.

\begin{theorem} \label{T:30.2}  A non-trivial $CR$-warped product $\,N_T\times_{f} N_\perp$   in $CP^m(4)$ satisfies the basic equality case of inequality \e{30.1}, i.e., \begin{align}\label{30.2} ||h||^2= 2p||\nabla (\ln f)||^2\end{align}  if and only if we have
\vskip.04in 

$(1)$  $N_T$ is an open portion of  complex Euclidean $h$-space {\bf C}$^h$,
\vskip.04in 

$(2)$ $N_\perp$ is an open portion of  a unit $p$-sphere $S^p$, and
\vskip.04in 

$(3)$ up to rigid motions, the immersion {\bf x} of $N_T\times_{f} N_\perp$ into $CP^m(4)$ is  the composition $\pi\circ\breve{\hbox{\bf x}}$, where \begin{equation}\begin{aligned} \notag\breve{\hbox{\bf x}}&(z,w)=\Bigg( z_0+(w_0-1)a_0\sum_{j=0}^h a_j z_j,\,\cdots ,  z_h+(w_0-1)a_h\sum_{j=0}^h a_jz_j,\\&\hskip.8in  w_1\sum_{j=0}^h a_jz_j, \,\ldots,\; w_p\sum_{j=0}^ha_jz_j,0,\,\ldots,0\Bigg),
\end{aligned}\end{equation} where $\pi$ is the projection $\pi:\hbox{\bf C}^{m+1}_*\to CP^m(4)$, $z=(z_0,z_1,\ldots,z_h)
\in \hbox{\bf C}^{h+1}$ and $w=(w_0,\ldots,w_p)\in
S^p\subset \hbox{\bf E}^{p+1}$,  and  $a_0,\ldots,a_h$ are real numbers satisfying $$a_0^2+a_1^2+\ldots+a_h^2=1.$$  \end{theorem}

\begin{theorem} \label{T:30.3}  A $CR$-warped product $N_T\times_{f} N_\perp$  in $CH^m(-4)$ satisfies the equality case of \e{30.1} if and only if one of the following two cases occurs:
\vskip.04in 

$(1)$ $N_T$ is an open portion of  complex Euclidean $h$-space {\bf C}$^h$,  $N_\perp$ is an open portion of  a unit $p$-sphere $S^p$ and, up to rigid motions, the immersion  is  the composition $\pi\circ\breve{\hbox{\bf x}}$, where $\pi$ is the projection $\pi:\hbox{\bf C}^{m+1}_{*1}\to CH^m(-4)$ and 
 \begin{equation}\begin{aligned} \notag &\breve{\hbox{\bf x}}(z_,w)=\Bigg( z_0+a_0(1-w_0)\sum_{j=0}^h a_j z_j,z_1+a_1(w_0-1)\sum_{j=0}^h a_j z_j,\,\cdots, \\&\hskip.1 in  z_h+a_h(w_0-1)\sum_{j=0}^h a_jz_j,\;w_1\sum_{j=0}^h a_jz_j,\ldots,\; w_p\sum_{j=0}^h
a_jz_j,0,\ldots,0\Bigg), \end{aligned}\end{equation} 
 where  $z=(z_0,\ldots,z_h)\in\hbox{\bf C}^{h+1}_1$,  $w=(w_0,\ldots,w_p)\in S^p\subset\hbox{\bf E}^{p+1}$ and  $a_0,\ldots,a_{h}$ are real numbers  satisfying $a_0^2-a_1^2-\cdots-a_{h}^2=-1.$
\vskip.04in 

$(2)$ $p:=\dim N_\perp=1$,   $N_T$ is an open portion of  {\bf C}$^h$ and, up to rigid motions, the immersion  is  the composition $\pi\circ\breve{\hbox{\bf x}}$, where 
\begin{equation}\begin{aligned} \notag &\breve{\hbox{\bf x}}(z,t)=\Bigg( z_0+a_0(\cosh t-1)\sum_{j=0}^h a_j z_j,  z_1+a_1(1-\cosh
t)\sum_{j=0}^h a_jz_j,\\&\hskip.3in \ldots,z_h+a_h(1-\cosh t)\sum_{j=0}^h a_jz_j, \sinh t\sum_{j=0}^h a_jz_j,0,\ldots,0\Bigg)
\end{aligned}\end{equation} 
for some real numbers $a_0,a_1\ldots,a_{h}$ satisfying $a_0^2-a_1^2-\cdots-a_{h}^2=1.$
 \end{theorem}
 
 There exists another general optimal inequality involving $||h||^2$  for $CR$-warped products in complex space forms. In fact, we have the following results form \cite{c30}.
 
 \begin{theorem}\label{T:30.4} Let $\phi:N_T\times_f N_\perp\to \hbox{\bf C}^m$ be a  $CR$-warped product in complex Euclidean $m$-space {\bf C}$^m$. Then we have
\vskip.04in

{\rm (1)}  The squared norm of the second fundamental form of $\phi$ satisfies
\begin{equation}\label{30.3}||h||^2\geq 2p\big\{||\nabla(\ln f)||^2+\Delta(\ln f)\big\}.\end{equation}
\vskip.04in

{\rm (2)} If the $CR$-warped product satisfies the equality case of \eqref{30.3}, then we have
\vskip.04in

{\rm (2.a)} $N_T$ is an open portion of {\bf C}$^h_*$;
\vskip.04in

{\rm (2.b)} $N_\perp$ is an open portion of $S^p$;
\vskip.04in

{\rm (2.c)} There exists a natural number $\alpha\leq h$ and  a  complex coordinate system $\{z_1,\ldots,z_h\}$ on {\bf C}$^h_*$ such that the warping function $f$ is given by $$f=\sqrt{\sum_{j=1}^\alpha z_j\bar z_j};$$
\vskip.04in

{\rm (2.d)} Up to rigid motions of {\bf C}$^m$, the immersion $\phi$ is given by $\phi^{hp}_\alpha$ in a natural way; namely, we have
\begin{equation}\notag \phi(z,w)=\big(w_0 z_1,\ldots,w_p z_1,\ldots,w_0z_\alpha ,\ldots, w_p z_\alpha ,z_{\alpha +1},\ldots,z_h,0,\ldots,0\big)\end{equation}
 for $z=(z_1,\ldots,z_h)\in \hbox{\bf C}_*^h$ and $\,w=(w_0,\ldots,w_p)\in S^p\subset {\mathbb E}^{p+1}$.
\end{theorem}

\begin{theorem}\label{T:30.5} Let $\phi:N_T\times_{f} N_\perp\to CP^m(4)$  be a $CR$-warped product. Then 
\vskip.04in 

{\rm (1)} The squared norm of the second fundamental form of $\phi$ satisfies
\begin{equation}\label{30.4}||h||^2\geq 2p\big\{||\nabla (\ln f)||^2+\Delta(\ln f)\}+4hp.\end{equation}
\vskip.04in 

{\rm (2)} The $CR$-warped product  satisfies the equality case of \eqref{30.4} if and only if 
\vskip.04in 

{\rm (2.i)}  $N_T$ is an open portion of  complex projective $h$-space $CP^h(4)$;
\vskip.04in 

{\rm (2.ii)}  $N_\perp$ is an open portion of   unit $p$-sphere $S^p$; and
\vskip.04in 

{\rm (2.iii)} There exists a natural number $\alpha\leq h$ such that, up to rigid motions,  $\phi$ is the composition $\pi\circ\breve{\phi}$, where  \begin{equation}\begin{aligned}\notag \breve{\phi}(z,w)=&\Big(w_0 z_0,\ldots,w_p z_0,\ldots,w_0z_\alpha ,\ldots, w_p z_\alpha
, z_{\alpha +1}, \ldots,z_h,0\ldots,0\Big)\end{aligned}\end{equation} for $z=(z_0,\ldots,z_h)\in \hbox{\bf C}_*^{h+1}$ and $\,w=(w_0,\ldots,w_p)\in
S^p\subset {\mathbb E}^{p+1}$, and  $\pi$ being the projection $\pi:\hbox{\bf C}^{m+1}_*\to CP^m(4)$.
\end{theorem} 

\begin{theorem}\label{T:30.6} Let $\phi:N_T\times_{f} N_\perp\to CH^m(-4)$  be a $CR$-warped product. Then
\vskip.04in 

{\rm (1)} The squared norm of the second fundamental form of $\phi$ satisfies
\begin{equation}\label{30.5}||h||^2\geq 2p\big\{||\nabla (\ln f)||^2+\Delta(\ln f)\}-4hp.\end{equation}
\vskip.04in 

{\rm (2)} The $CR$-warped product satisfies the equality case of \eqref{30.5} if and only if 
\vskip.04in 

{\rm (2.a)}  $N_T$ is an open portion of  complex hyperbolic $h$-space $CH^h(-4)$;
\vskip.04in 

{\rm (2.b)}  $N_\perp$ is an open portion of  unit $p$-sphere $S^p$ $($or {\bf R},  when $p=1$$)$; and
\vskip.04in 

{\rm (2.c)} up to rigid motions, $\phi$ is the composition $\pi\circ\breve{\phi}$, where either $\breve\phi$ is 
 \begin{equation}\begin{aligned} \notag \breve{\phi}(z,w)=\Big(z_{0},&\ldots,z_{\beta},w_0 z_{\beta+1},\ldots,w_p z_{\beta+1},\ldots, w_0z_h ,\ldots, w_p z_h ,0\ldots,0\Big)\end{aligned}\end{equation} for  $0<\beta\leq h$,  $z=(z_0,\ldots,z_h)\in \hbox{\bf C}_{*1}^{h+1}$ and $\,w=(w_0,\ldots,w_p)\in S^p$, or $\breve \phi$ is 
\begin{align}\notag &\breve{\phi}(z,u)=\Big(z_0\cosh u,z_0\sinh u,z_1\cos u,z_1\sin u,\ldots,\ldots,\\&\hskip.9in z_\alpha\cos u,z_\alpha\sin u,z_{\alpha+1},\ldots,z_h,0\ldots,0\Big) \notag\end{align}
 for $z=(z_0,\ldots,z_h)\in \hbox{\bf C}_{*1}^{h+1}$ and $\,u\in \hbox{\bf R}$, and  $\pi$ being the projection $$\pi:\hbox{\bf C}^{m+1}_{*1}\to CH^m(-4).$$
\end{theorem}

\subsection{$CR$-warped products with compact holomorphic factor $N_T$}

When the holomorphic factor $N_T$ of a $CR$-warped product is a compact manifold, we have following additional  results from  \cite{c35} also involving the squared norm of the second fundamental form.

\begin{theorem}  \label{T:30.7} For any $CR$-warped product $N_T\times_f N_\perp$ in $CP^m(4)$ with compact $N_T$ and any $q\in N_\perp$, we have
\begin{align}\label{30.6}\int_{N_T\times\{q\}} ||h||^2{\rm d}V_{T}\geq 4hp\,\hbox{\rm vol}(N_T).\end{align} 

The equality sign of \eqref{30.6} holds identically if and only if  we have:
\vskip.04in 

{\rm (i)} The warping function $f$ is constant.
\vskip.04in 

{\rm (ii)}  $(N_T,g_{N_T})$ is holomorphically isometric to $CP^h(4)$ and it is isometrically immersed in $CP^m(4)$ as a totally geodesic complex submanifold.
\vskip.04in 

{\rm (iii)} $(N_\perp,f^2 g_{N_\perp})$ is isometric to an open portion of the real projective $p$-space $RP^p(1)$ of constant sectional curvature one  and it is isometrically immersed in $CP^m(4)$ as a totally geodesic totally real submanifold.
\vskip.04in 

{\rm (iv)} $N_T\times_f N_\perp$ is immersed linearly fully in a linear complex subspace $CP^{h+p+hp}(4)$ of $CP^m(4)$.

Moreover, the immersion is rigid.
\end{theorem}

\begin{theorem}  \label{T:30.8} If $\,N_T\times_f N_\perp\,$ is a $CR$-warped product in $\,CP^{h+p+hp}(4)\,$ with compact $N_T$, then $N_T$ is holomorphically isometric to $CP^h(4)$.
\end{theorem}

\begin{theorem}  \label{T:30.9} Let $N_T\times_f N_\perp$ be a $CR$-warped product with compact $N_T$ in $CP^m(4)$. If the warping function $f$ is non-constant, then, for each $q\in N_\perp$, we have
\begin{align}\label{30.7}\int_{N_T\times\{q\}} ||h||^2{\rm d}V_{T}\geq 2p\lambda_1 \int_{N_T}(\ln f)^2{\rm d}V_T+ 4hp\,\hbox{\rm vol}(N_T),\end{align}  where ${\rm d}V_T$, $\lambda_1$ and $\hbox{\rm vol}(N_T)$ are the volume element, the first positive eigenvalue of the Laplacian $\Delta$ and  the volume of $N_T$, respectively.

Moreover, the equality sign of \eqref{30.7} holds identically if and only if we have 
\vskip.04in 

{\rm (a)} $\Delta\ln f=\lambda_1 \ln f$.
\vskip.04in 

{\rm (b)} The $CR$-warped product is both $N_T$-totally geodesic and $N_\perp$-totally geodesic.
\end{theorem}

The following example shows that Theorems \ref{30.7}, \ref{30.8} and \ref{30.9}  are sharp.

\begin{example} Let $\iota_1$ be the identity map of $:CP^h(4)$ and let $$\iota_2:RP^p(1)\to CP^p(4)$$ be a totally geodesic Lagrangian embedding of $RP^p(1)$ into $CP^p(4)$. Denote by $$\iota=(\iota_1,\iota_2):CP^h(4)\times
RP^p(1)\to CP^h(4)\times CP^p(4)$$ the product embedding of $\iota_1$ and $\iota_2$. 
Moreover, let $S_{h,p}$ be the Segre embedding of  $CP^h(4)\times CP^p(4)$ into $CP^{hp+h+p}(4)$. Then the composition $\phi=S_{h,p}\circ\iota$:
\begin{align}  CP^h(4) &\times RP^p(1)\xrightarrow[\text{totally geodesic}]  {\text{$(\iota_1,\iota_2)$}}  CP^h(4)\times CP^p(4) \notag \\& \hskip.2in\xrightarrow[\text{Segre
embedding}] {\text{$S_{h,p}$}} CP^{hp+h+p}(4)\notag\end{align} is a $CR$-warped product in $CP^{h+p+hp}$ whose holomorphic factor $N_T=CP^h(4)$ is a compact manifold. 

 Since the second fundamental form of $\phi$ satisfies the equation: $||\sigma||^2=4hp$, we have the equality case of \eqref{30.6} identically. 
\end{example} 
 
The next example shows that the assumption of compactness in Theorems \ref{30.7}, \ref{30.8} and \ref{30.9} cannot be removed.

\begin{example} Let $\hbox{\bf C}^*=\hbox{\bf C}-\{0\}$ and $\hbox{\bf C}_*^{m+1}=\hbox{\bf C}^{m+1}-\{0\}$. Denote by $\{z_0,\ldots,z_h\}$   a natural complex coordinate system on $\hbox{\bf C}^{m+1}_*$.  Consider the action of
$\hbox{\bf C}^*$ on $\hbox{\bf C}_*^{m+1}$ defined by
$$\lambda\cdot (z_0,\ldots,z_m)=(\lambda z_0,\ldots,\lambda z_m)$$ for $\lambda\in {\bf C}_*$.  Let
$\pi(z)$ denote the equivalent class containing $z$ under this action. Then the  set of equivalent classes is the complex projective $m$-space $CP^m(4)$ with the  complex structure induced from the complex structure on {\bf C}$^{m+1}_*$. 

For any two natural numbers $h$ and $p$, we define a map: $$\breve \phi:\mathbb C^{h+1}_*\times S^p(1)\to\mathbb C^{h+p+1}_*$$ by
$$\breve \phi(z_0,\ldots,z_h;w_0,\ldots,w_p)=\big(w_0z_0,w_1z_0,\ldots,w_pz_0,z_1,\ldots,z_h\big)$$ 
for $(z_0,\ldots,z_h)$ in $\mathbb C^{h+1}_*$ and $(w_0,\ldots,w_p)$ in $S^{p}(1)$ with $\sum_{j=0}^p w_j^2=1$.

Since the image of $\breve \phi$ is invariant under the action of ${\bf C}_*$, the composition:
\begin{align}\pi\circ\breve \phi:{\bf C}^{h+1}_*\times S^p(1)\xrightarrow[\text{}] {\text{$\breve \phi$}}
{\bf C}^{h+p+1}_*\xrightarrow[\text{}] {\text{$\pi$}} CP^{h+p}(4)\notag\end{align} 
induces a $CR$-immersion of the product manifold $N_T\times S^p(1)$ into $CP^{h+p}(4)$, where $$N_T=\big\{(z_0,\ldots,z_h)\in CP^h(4):z_0\ne 0\big\}$$ is a proper open subset of $CP^h(4)$. Clearly, the induced metric on
$N_T\times S^p(1)$ is a warped product metric and the holomorphic factor $N_T$ is non-compact. 

Notice that the complex dimension of the ambient space is  $h+p$; far less than $h+p+hp$.
\end{example}

\subsection{Multiply $CR$-warped products} 

A multiply warped product $$N_T \times_{f_2} N_2 \times \cdots \times_{f_k} N_k$$ in a K\"ahler manifold is called a {\it multiply $CR$-warped product}  if $N_T$ is a holomorphic submanifold and $N_\perp:={}_{f_2} N_2 \times \cdots \times_{f_k} N_k$ is a totally real submanifold.

In \cite{cd}, Theorem \ref{30.1} was extended  to the following.

\begin{theorem} \label{T:30.10} Let $N=N_T \times_{f_2} N_2 \times \cdots \times_{f_k} N_k$ be a  multiply $CR$-warped product in an arbitrary Kaehler manifold $\tilde M$. Then the second fundamental form $h$ and the warping functions $f_2,\ldots, f_k$ satisfy
 \begin{align} \label{30.8} &||h||^2 \geq 2\sum_{i=2}^k n_i ||\nabla (\ln f_i)||^2. \end{align}
 
 The equality sign of \eqref{30.8} holds identically if and only if the following statements hold:
 \vskip.04in 

 {\rm (i)} $N_T$ is a  totally geodesic holomorphic submanifold of $\tilde M$;
\vskip.04in 
 
 {\rm (ii)} For each $i\in \{2,\ldots,k\}$,   $N_i$ is  a totally umbilical submanifold of $\tilde M$ with $-\nabla (\ln f_i)$ as its mean curvature vector; 
\vskip.04in 
   
 {\rm (iii)} ${}_{f_2} N_2 \times \cdots \times_{f_k} N_k$ is immersed as mixed totally geodesic submanifold in $\tilde M$;  and 
 \vskip.04in 

 {\rm (iv)} For each point $p\in N$, the first normal space ${\rm Im} \,h_p$ is a subspace of $J(T_pN_\perp)$, where $J$  is the almost complex structure of $\tilde M$.
 \end{theorem}

The following example shows that inequality \eqref{30.8} is sharp.

\begin{example} Assume that $h$ and $k$ are natural numbers with $h\geq k$. Let $$N_T={\bf  C}^{h}:=\{(z_1,\ldots,z_h): z_1,\ldots,z_h\in {\bf C}\}$$ and let $N_i=S^{n_i}$ denote the unit $n_i$-spheres for $i=2,\ldots,k$. 
 
 Consider the immersion $\psi$ of $N_T\times S^{n_2}\times \cdots\times S^{n_k}$ into ${\bf C}^{h+n_2+\cdots+n_k}$
 defined by 
\begin{align}\label{30.9} \psi= (z_1 w_{2,0},\ldots,z_1 w_{2,n_2}, \ldots, z_k w_{k,0},\ldots,z_k w_{k,n_k}, z_{k+1},\ldots,z_h),
\end{align}
where  
$(w_{i,0},\ldots,w_{i,n_i})\in {\bf R}^{n_i+1}$ satisfy $\sum_{\alpha=0}^{n_i} w^2_{i,\alpha}=1$ for $i=2,\ldots,k$.

It is easy to see that  the product manifold $${\bf C}^h\times S^{n_2}\times \cdots\times S^{n_k}$$ endowed with the induced metric via $\psi$ is the multiply warped product manifold $${\bf C}^h\times_{f_2} S^{n_2}\times \cdots\times_{f_k} S^{n_k}$$ with $f_i=|z_i|$.
Moreover, with respect to the canonical complex structure of ${\bf C}^{h+n_2+\cdots+n_k}$, the immersion $\psi$ is a 
 multiply $CR$-warped product submanifold.
 
A straightforward computation shows that this example of multiply $CR$-warped product submanifold satisfies the equality case of \eqref{30.8}. This examples shows that  inequality \eqref{30.8} is optimal.
\end{example}

By applying Theorem \ref{T:30.9} we have  the following.
 
\begin{corollary}\label{C:30.1} If $f_2,\ldots,f_k$ are  harmonic functions on $N_1$ or eigenfunctions of the Laplacian $\Delta$ on $N_1$ with positive eigenvalues, then the 
multiply warped product manifold $N_1 \times_{f_2} N_2 \times \cdots \times_{f_k} N_k$  cannot be isometrically immersed into any Riemannian manifold of negative sectional curvature as a minimal submanifold.
 \end{corollary}

\begin{example}\label{E:30.4}  Let  $M_1\times_{f_2} M_2\times \cdots\times_{f_k} M_k$ be a {\it multiply warped product representation\/} of  a Riemannian $m$-manifold $R^m(\epsilon)$ of constant curvature $\epsilon$. 
Assume that $\psi^i:N_i\to M_i,\, i=2,\ldots, k,$ are minimal immersions. Then the immersion: 
\begin{align} \psi: M_1\times_{f_2} N_2\times \cdots\times_{f_k} N_k\to M_1\times_{f_2} M_2\times \cdots\times_{f_k} M_k
\end{align} defined by $\psi=(id,\psi_2,\ldots,\psi_k)$ is a minimal isometric immersion of the multiply warped product manifold $M_1\times_{\rho_2} N_2\times \cdots\times_{\rho_k} N_k$ into $R^m(\epsilon)$.

On the other hand,  since $M_1\times_{f_2} M_2\times \cdots\times_{f_k} M_k$ is of constant curvature $\epsilon$,  the warping functions  $f_2,\ldots,f_k$ are  eigenfunctions of the Laplacian $\Delta$ of $M_1$ with eigenvalues given by $n_2\epsilon,\ldots,n_k\epsilon$, respectively. 
In particular, if $\epsilon=0$ the warping functions $f_2,\ldots,f_k$ are harmonic functions.

\medskip
Example \ref{E:30.4} illustrates that  the warping functions  $f_2,\ldots, f_k$ in Corollary \ref{C:30.1} cannot be replaced by eigenfunctions with negative eigenvalue. Moreover, the target space   in Corollary \ref{C:30.1} cannot be replaced either by Euclidean space or by spheres. Therefore, Corollary \ref{C:30.1} is sharp. 
\end{example}

\section[DDVV conjecture]{Inequality involving normal scalar curvature and DDVV conjecture}

\subsection{Inequalities of Chen and Guadelupe and Rodriguez}

The $\delta$-invariant $\delta(n_1,\ldots,n_k)$ reduces to the scalar curvature $\tau$ when the $k$-tuple $(n_1,\ldots,n_k)$ were chosen to be the empty set. Accordingly, Theorem \ref{T:5.1} reduces to the following \cite{c8}.

\begin{corollary}\label{C:31.1} Let $x:M \rightarrow R^m(\epsilon)$ be an isometric immersion of a Riemannian $n$-manifold $M$ with normalized scalar curvature $\rho$ into  an $m$-dimensional real space form $R^m(\epsilon)$  of constant sectional curvature $\bar c$. Then we have 
 \begin{equation}\label{31.1} H^2\geq \rho-\epsilon,\end{equation}
  equality holding at a point $p\in M$ if and only if $p$ is a totally umbilical point.\end{corollary}
  
  For a surface $M$ in a Riemannian manifold, the {\it ellipse of curvature} $E(x)$ at a point $x\in M$ is defined by
  $$E(x)=\{h(X,X)\,:\, X\in T_xM \;{\rm and} \; ||x||=1\},$$
where $h$ is the second fundamental form of $M$.

A  surface $M$ in $\mathbb E^4$ is called {\it superconformal}  if all of its ellipses of curvature of $M$ are circles.

On the other hand, Guadelupe and Rodriguez proved the following:

\begin{theorem} {\rm \cite{GR}}\label{T:31.1}  Let $M^2$ be a surface in a real space form $R^{2+m}(\epsilon)$. Denote by $K$  the Gaussian curvature of $M^2$ and by $K^\perp$ the normal scalar curvature. Then we have 
\begin{align}\label{31.0} K \leq   H^2 - K^\perp+\epsilon\end{align}
at every point $p\in M^2$, with equality if and only if the ellipse of curvature at $p$ is a circle.\end{theorem}

\subsection{DDVV conjecture}
Consider a submanifold $M^n$ of a real space form $R^{n+m}(\epsilon)$, the normalized normal scalar curvature $\rho^\perp$ is defined as
$$\rho^\perp=\frac{2}{n(n-1)}\sqrt{\sum_{1\leq i<j\leq n; 1\leq r<s\leq m}\hskip-.3in \<R^\perp(e_i,e_j)\xi_r,\xi_s\>^2}.$$

 Since the normalized scalar curvature $\rho$ is the higher-dimensional analogue of the Gaussian curvature $K$  and the normalized normal scalar curvature $\rho^\perp$  that of $K^\perp$,     De Smet, Dillen,  Verstraelen and Vrancken thus made the following conjecture in \cite{DDVV1}:
\vskip.1in

\noindent {\bf DDVV Conjecture.}  {\it Let $M^n$ be a submanifold of a real space form $R^{n+m}(\epsilon)$ of constant sectional curvature $\epsilon$. Denote by $\rho$ the normalized scalar curvature  and by $\rho^\perp$  the normalized normal scalar curvature. Then}
\begin{align}\label{31.2}  \rho \leq   H^2 -\rho ^\perp+\epsilon.\end{align}
\vskip.1in

DDVV Conjecture  is true for submanifolds of dimension $n$ and codimension $m$  if,  for every set $\{B_1,\ldots, B_m\}$ of symmetric $(n\times n)$-matrices with trace zero, the following inequality holds:
\begin{equation}\sum_{\alpha,\beta=1}^m ||[B_\alpha,B_\beta]||^2\leq \(\sum_{\alpha=1}^m ||B_\alpha||^2\)^2.\end{equation}

For normally flat submanifolds, in particular for hypersurfaces, inequality \e{31.2} is nothing but inequality \e{31.1}.  Moreover, the conjecture was proven for immersions with codimension 2 in \cite{DDVV2}; for immersions which are invariant with respect to the standard K\"ahlerian and Sasakian structures on $\mathbb E^{2k}$  and $S^{2k+1}(1)$ respectively in \cite{DFV1}; and for immersions which are totally real with respect to the nearly K\"ahler structure on $S^6(1)$ in \cite{DDVV1}. 

Let $N$ be a submanifold of a Riemannian manifold $M$. Then, according to \cite{HL},  $N$ is called  {\it austere\/} if for each normal vector $\xi$ the set of eigenvalues of $A_{\xi}$ is invariant under multiplication by $-1$; this is equivalent to the condition that all the invariants of odd order of the Weingarten map at each normal vector of $N$ vanish
identically.   Of course every austere submanifold is a minimal
submanifold. 

Recently,  Choi and Lu confirm  in \cite{CL} that DDVV conjecture is also true for 3-dimensional submanifolds with arbitrary codimension. Moreover, they prove the following result for 3-dimensional submanifolds.

\begin{theorem} If the equality in DDVV Conjecture  is valid at any point and M is minimal, then M is an austere $3$-fold.\end{theorem} 

Explicit construction of  Euclidean submanifolds of codimension two, free of minimal and umbilical points, that attain equality case of \e{31.2} are constructed in \cite{DT}. In particular, in the two-dimensional case, the construction yields all superconformal surfaces in $\mathbb E^4$. 

The notion of $H$-umbilical Lagrangian submanifolds introduced in \cite{c.4.1,c.4.2} is defined as follows:

A Lagrangian submanifold $M^n$ of a K\"ahler manifold $\tilde M^n$ is called {\it $H$-umbilical}  if it is a
 non-totally geodesic Lagrangian submanifold whose second fundamental form takes the following simple form:  \begin{equation}  \begin{aligned}  \label{31.3} & h(e_1,e_1)= \lambda Je_1,\quad h(e_2,e_2)=\cdots = h(e_n,e_n)=\mu Je_1,\\ & h(e_1,e_j)=\mu Je_j,\quad h(e_j,e_k)=0,\; \; j\not=k, \;\;\;\; j,k=2,\ldots,n\end{aligned}\end{equation}  for some suitable functions $\lambda$ and $\mu$ with respect to some suitable orthonormal local frame field $e_1,\ldots,e_n$. It is obvious  that  condition \e{31.3} is equivalent to 
\begin{equation}  \begin{aligned}  & h(X,Y)=\alpha\<JX,H\>\<JY,H\>H\\ & \quad+\beta\<H,H\>\{\<X,Y\>H+\<JX,H\>JY +\<JY,H\>JX\}\end{aligned}\end{equation}  for
any vectors $X,Y$ tangent to $M$, where
$$\alpha={{\lambda-3\mu}\over{\gamma^3}},\quad \beta={\mu\over{\gamma^3}},\quad\gamma={{\lambda+(n-1)\mu}\over n}$$ when $H\not=0$. 

According to \cite{c12}, an $n$-dimensional submanifold $M$  of a Riemannian $m$-manifold $N^m$ is called  {\it ultra-minimal\/} if, with respect to some suitable locally orthonormal frame fields, $e_1 \ldots,e_n,e_{n+1},\ldots,e_m$, the shape operators $A_r$ take the  form:
$$A_r=\left( \begin{array}{llllll} A^r_{1} & \cdots & 0&0 &\cdots & 0\\\vdots  & \ddots& \vdots&\vdots  & \ddots &\vdots \\ 0 &\cdots &A^r_k&0 &\cdots & 0 \\ 0&\cdots&0&0&\cdots&0\\\vdots  & \ddots & \vdots &\vdots &\ddots&\vdots \\0  &\cdots& 0&0 &\cdots & 0 \end{array}\right),\;\;\;\;\mbox{ trace}\; A^r_j=0 ,
$$ where $A^r_j, n+1\leq r\leq m,j=1,\ldots,k,$ are $n_j\times n_j$ symmetric submatrices. 

It is known that every ultra-minimal submanifold is always an ideal submanifold (see \cite{c12}). 

In \cite{DFV2}, F. Dillen, J. Fastenakels  and  J. Van der Veken  also show that the DDVV conjecture also holds for 4-dimensional ultra-minimal submanifolds in ${\bf C}^4$.

Also, Dillen, Fastenakels  and Van der Veken  prove in \cite{DFV1} the following.

\begin{proposition} Let $M^{2n+1}$ be a submanifold of $S^{2m+1}(1)$ which is invariant with respect to the standard Sasakian structure on the unit sphere. Then we have $\rho+\rho^\perp\leq 1$.
\end{proposition}

Recently, the DDVV conjecture was proved to be true in general by Zhiqin Lu in \cite{Luz}.

\subsection{Lorentzian version} 
Let $M$ be an $m$-dimensional, time-orientable, Lorentzian manifold and $\mathbb E^{m+2}_s$  an $(m+2)$-dimensional pseudo-Euclidean space of signature $(+,\cdots,+,-,\varepsilon_{m+1},\varepsilon_{m+2})$ with $\varepsilon_A=\pm 1, A = m+1,m+2$.

We consider $M$ to be locally and isometrically embedded in $\mathbb E^{m+2}_s$. Let $\Omega_{\alpha\beta}$ and $\Lambda_{\alpha\beta}$ be the components of the second fundamental form with respect to $\xi_{m+1},\xi_{m+2}$, respectively. Then the mean curvature vector $\overrightarrow H$ is defined as
$$\overrightarrow H=\frac{1}{m}(\varepsilon_{m+1}\Omega_\alpha^\alpha \xi_{m+1}+\varepsilon_{m+2}\Lambda_\alpha^\alpha \xi_{m+2}).$$

Let  $\{e_1,\ldots,e_{m-1},e_m\}$  be an orthonormal basis of $M$. Due to space-time applications
in mind we take $M$  time-orientable such that there exists a global, nowhere zero,
timelike vector field, denoted by $e_m$. 

Let $h$ be the second fundamental form of $M$ in $\mathbb E^{m+2}_s$. Put
$$h(X,Y)=\sum_{r=m+1}^{m+2} \varepsilon_r \eta(\tilde \nabla_X Y,\xi_r)\xi_r.$$ Then we have
$$\Omega_{\alpha m}=-\eta(e_m,\tilde\nabla_{e_\alpha}\xi_{m+1})=-\eta(e_\alpha,\tilde\nabla_{e_m}\xi_{m+1})$$
with $\alpha=1,\ldots,m-1$, and analogous relations hold for $\Lambda_{\alpha m}$. 

\begin{definition} An embedding $\phi: (M^m,g)\to  \mathbb E^{m+2}_s$ with $\varepsilon_{m+1}=\varepsilon_{m+2}=1$  is called {\it causal-type preserving} if and only if, with respect to some orthonormal basis $\{e_1,\ldots,e_m\}$, $\tilde\nabla_{e_\alpha} \xi_A$  is space-like, for $A = m + 1,m + 2$ and $\alpha = 1, . . . ,m - 1$. \end{definition}

\begin{definition} An embedding $\phi: (M^m,g)\to  \mathbb E^{m+2}_s$ with $\varepsilon_{m+1}=\varepsilon_{m+2}=-1$  is called {\it causal-type preserving} if and only if, with respect to some orthonormal basis $\{e_1,\ldots,e_m\}$, $\tilde\nabla_{e_m} \xi_A$  is time-like, for $A = m + 1,m + 2$. \end{definition}

Notice that causal-type preserving embeddings have $\Omega_{\alpha m}=\Lambda_{\alpha m}=0$ for each $\alpha=1,\ldots,m-1$.

If  $v = v_1\xi_{m+1} +v_2\xi_{m+2}$ is a vector in normal space, we define the norm:
 $$||v||^2_\perp= \varepsilon_{m+1}(v_1)^2+\varepsilon_{m+2} (v_2)^2.$$
Using these definitions we define the scalar normal curvature as
$$\rho^\perp=\frac{\sqrt{2}}{m(m-1)} ||[\Omega,\Lambda] ||.$$

In \cite{DHPV}, F. Dillen, S. Haesen, M.  Petrovi\'c and L. Verstraelen   extended  inequality \eqref{31.0} to  Lorentzian manifolds in a pseudo-Euclidean space with codimension 2. They obtain the following.

\begin{theorem}\label{T:31.3}  Let $\phi: (M^m,g)\to  \mathbb E^{m+2}_s$  be a causal-type preserving, local and isometric
embedding of a Lorentzian manifold $M^m$ in a pseudo-Euclidean space $ \mathbb E^{m+2}_s$. Then we have
\begin{align} ||\overrightarrow H||^2_\perp \geq \rho+\rho^\perp,\;\; \text{if $\varepsilon_{m+1}=\varepsilon_{m+2}=1$} \end{align}
 and
 \begin{align} ||\overrightarrow  H||^2_\perp \leq \rho+\rho^\perp,\;\; \text{if $\varepsilon_{m+1}=\varepsilon_{m+2}=-1$}.\end{align}
\end{theorem}

In \cite{DHPV},  they also determine the special form of the second fundamental form when the equality case  of Theorem \ref{T:31.3} occurs.  Furthermore, they showed that space-times which realize the equality case of Theorem \ref{T:31.3} are Petrov type D anisotropic fluid models with a time-like surface of constant curvature.

 \subsection{An optimal inequality involving normal $\delta$-invariant $\hat\delta^\perp(2)$ for K\"ahler hypersurfaces}
Recall that the $\delta$-invariant $\hat \delta(2)$ is defined as  $$\hat\delta(2) =\tau-\max K.$$ 
Similarly, Z. Sent\"urk and L.  Verstraelen \cite{SV} consider the scalar curvature function defined by \begin{align}\hat \kappa (2) =\tau   - \max H(e),\end{align} where $H(e)$ denotes the holomorphic sectional curvature with respect to a unit tangent vector $e$.

 Let $R^\perp$ denote the normal curvature tensor of a K\"ahler hypersurface $M^n$ in $CP^{n+1}(4\epsilon)$.
 we put 
\begin{align}&K^\perp_{\alpha\beta}=R^\perp(e_\alpha,e_\beta,\xi,J\xi),\\& \delta^\perp=\sum_{\alpha<\beta} |K^\perp_{\alpha\beta}|,\end{align}
 where $\xi$ is a unit normal vector and  $e_\alpha,e_\beta$ are orthogonal unit tangent vectors of $M^n$ in $CP^{n+1}(4\epsilon)$.

The scalar normal $\delta$-invariant $\hat\delta^\perp(2)$ is defined as 
 \begin{align}&\hat\delta^\perp(2) =\delta^\perp-\max |K^\perp|.\end{align}

It is proved by  Sent\"urk and  Verstraelen in \cite{SV}  that

\begin{theorem} For every K\"ahler hypersurface $M^n$ in the complex projective space $CP^{n+1}(4)$, the scalar valued curvatures  $\hat\kappa(s)$ and $\hat\delta^\perp(2)$ always satisfy the following inequality:
\begin{align}&\hat\kappa(2)\leq -\hat\delta^\perp(2) +2(n-1)(4n+3).\end{align}
And the only K\"ahler hypersurfaces $M^n$ in $CP^{n+1}(4)$ which are ideal in this
respect, i.e. for which the equality holds at all of their points, are  either open parts of
the totally geodesic complex space forms $CP^n(4)$ or open parts of the complex quadrics $Q_n$ in $CP^{n+1}(4)$.
 \end{theorem}

\section{Inequalities involving scalar curvature}

\subsection{Lagrangian submanifolds} 
There exists a general inequality for Lagrangian submanifolds in a complex space form involving the scalar curvature. In fact,  we have the following.

\begin{theorem} The scalar curvature $\rho$ and the squared mean curvature $H^2$ of a  Lagrangian submanifold $M$ in  complex space form $\tilde M^n(4\epsilon)$ satisfy the following general sharp inequality:
\begin{align}\label{La}H^2\geq {{2(n+2)}\over{n^2(n-1)}}\tau -\left( {{n+2}\over n}\right)\epsilon.\end{align}
The equality sign  holds if and only if, with
respect to an adapted Lagrangian frame field
$e_1,\ldots,e_n,e_{1^*},\ldots,e_{n^*}$ with $e_{1^*}$ parallel
to $JH$, the second fundamental form $h$ of $M$ in $\tilde M^n(4\epsilon)$
takes the following form: 
\begin{equation}\begin{aligned}& h(e_1,e_1)=3\lambda e_{1^*},\quad
h(e_2,e_2)=\cdots=h(e_n,e_n)=\lambda e_{1^*},\\  &h(e_1,e_j)=\lambda
e_{j^*}, \quad h(e_j,e_k)=0,
\quad 2\leq j\not= k\leq n.\end{aligned}\end{equation}\end{theorem}

Inequality \e{La} with $\epsilon=0$ and $n=2$ was proved in \cite{castro}. 
Their proof relies on complex analysis which is not applicable to $n\geq 3$. The general inequality was established in \cite{bcm} for
$\epsilon=0$ and arbitrary $n$; and in \cite{c10} for $\epsilon\ne 0$ and arbitrary $n$; (and independently in \cite{castro2}, for $\epsilon\ne 0$ with $n=2$, also using the method of complex analysis). 

If $\tilde M^n(4\epsilon)={\bf C}^n$, the equality of \e{La}  holds identically if and only if either the Lagrangian submanifold $M$ is an open portion of a Lagrangian $n$-plane or, up to dilations, $M$ is an open portion of the Whitney immersion \cite{bcm,ros}.

 Let $\cn(x,k)$ and $\dn(x,k)$ be the usual Jacobi's elliptic functions with modulus $k$. It is well-known that $\cn(x,k),\dn(x,k)$ are doubly periodic functions. We put $$\mu_a=\text{\small${{\sqrt{a^2-1}}\over \sqrt{2}}$}\,\cn\(ax,\text{\small${{\sqrt{a^2-1}}\over{\sqrt{2}a}}$}\),\quad
a>1,$$
$$\eta_a=\text{\small${{\sqrt{a^2+1}}\over\sqrt{2}}$}\,\dn\(\text{\small${{\sqrt{a^2+1}}\over \sqrt{2}}$}x, \text{\small$ {{\sqrt{2}a}\over{\sqrt{a^2+1}}}$}\),\quad
0<a<1,$$ 
$$\rho_a=\text{\small${{\sqrt{a^2+1}}\over\sqrt{2}}$}\,\cn\(ax, \text{\small${{\sqrt{a^2+1}}\over{\sqrt{2}a}}$}\),\quad a>1.$$

Let $S^{n}(\epsilon)$ and $H^{n}(-\epsilon)$ denote the $n$-sphere with constant sectional curvature $\epsilon$ and the real hyperbolic $n$-space with constant sectional curvature $-\epsilon$, respectively.
For $n\geq 3$, we denote by $P^n_a, D^n_{a}$ and $C^n_{a}$ the warped products $I\times_{\mu_a} S^{n-1}({{a^4-1}\over 4})$,
$\mathbb R\times_{\eta_a} H^{n-1}({{a^4-1}\over 4})$ and $I\times_{\rho_a}S^{n-1}({{a^4-1}\over 4})$ with warped functions
$\mu_a, \eta_a$ and $\rho_a$, respectively, and $I$  are the maximal open intervals containing 0 on which the corresponding warped functions are positive. 

Moreover, we denote by $F^n$ and $L^n$ the warped products $\mathbb R\times_{1/\sqrt{2}} H^{n-1}(-{1\over4})$ and $\mathbb
R\times_{\hbox{sech}(x)}\mathbb R^{n-1}$, respectively.
For $n=2$ we shall replace $S^{n-1}({{a^4-1}\over 4})$ or $ H^{n-1}({{a^4-1}\over 4})$ by the real line $\mathbb R$ to
 define $P^2_a,D^2_a,C^2_a,F^2,$ and $L^2$.

It is easy to see that  $F^2$ is a  flat surface and   $D^n_{a}$, $F^n$ and $L^n$ are  complete Riemannian $n$-manifolds, but $P^n_a$ and $C^n_a$ are not complete. Furthermore, these  Riemannian $n$-manifolds are conformally flat. Topologically,
$S^n$ is the two point compactification of both $P_a^n$ and  $C_a^n$

In \cite{c10}, the author proved that the one-parameter family of Riemannian $n$-manifolds, $P^n_a\, (a>1)$,  admit Lagrangian isometric immersions into $CP^n(4)$ satisfying the equality case of the inequality \e{La} for $\epsilon=1$;  the  two one-parameter families of Riemannian manifolds,  $\,C^n_a\, (a>1),D_a^n\,$ $ (0<a<1)$, and the  two exceptional $n$-spaces,  $F^n$ and $ L^n$, admit Lagrangian isometric immersion into
$ CH^n(-4)$ satisfying the equality case of the inequality for $\epsilon=-1$. 

It also proved in \cite{c10} that besides the totally geodesic ones, these are the only Lagrangian submanifolds in $  CP^n(4)$ and in $CH^n(-4)$ which satisfy the
 equality case of \e{La}  (for the case $n=2$, see also \cite{castro2}). 

The explicit expressions of those Lagrangian immersions of $P^n_a, C^n_a, D^n_a, F^n$ and $L^n$ satisfying the equality case of \e{La}  were completely determined  by B. Y. Chen and L. Vrancken in \cite{CV1}.

I. Castro and F. Urbano \cite{castro2} showed that a Lagrangian surface in $CP^2$ satisfies the equality case of \e{La} for $n=2$ and $c=1$ if and only if the Lagrangian surface has holomorphic twistor lift. 

\subsection{Slant submanifolds} 

Proper  slant submanifolds are  even-dimensional. \ Such submanifolds do exist extensively for any even dimension greater
than zero (cf. \cite{cbook3}).

A proper slant submanifold is called {\it K\"ahlerian slant } if the endomorphism $P$   is parallel with respect to the Riemannian connection, that is, $\nabla P = 0$. A K\"ahlerian slant submanifold is a K\"ahler manifold with respect to the induced metric and  the almost complex structure defined by ${\tilde J}= (\sec\theta)P$. 
K\"ahler submanifolds, totally real submanifolds and slant surfaces in a K\"aehler manifold  are examples of K\"ahlerian slant submanifolds. 

In general, let $M$ be a submanifold of a K\"ahler manifold $\tilde M$. Then $M$ satisfies $\nabla P=0$ if and only if $M$ is locally the Riemannian
product $M_1\times\cdots\times M_k$, where each $M_i$ is a K\"ahler submanifold, a totally real submanifold or a K\"ahlerian slant submanifold of $\tilde M$ (see \cite{cbook3}).

Slant submanifolds have the following topological properties:

\begin{theorem} \label{T:32.2} We have

{\rm (1)}  Let  $M$  be a compact $2k$-dimensional proper slant submanifold of a K\"ahler manifold, then $:$
$$\,H^{2i}(M;\hbox{\bf R})\ne \{0\}$$ for $i=1,\ldots,k$.
\vskip.04in

{\rm (2)}  Let $M$ be a slant submanifold in a complex Euclidean space. If $M$ is not totally real, then $M$ is non-compact.
\end{theorem}

Statements (1) and (2) of Theorem \ref{T:32.2} are due to \cite{cbook3} and \cite{CT2}, respectively.

An immediate consequence of statement (1) is the following.

\begin{corollary} If $M$ is a compact $2k$-manifold with $\,H^{2i}(M;\hbox{\bf R})= \{0\}$ for some $i\in \{1,\ldots,k\}$, then $M$ cannot be immersed in any K\"ahler manifold as a proper slant submanifold.
\end{corollary}

Although there  exist  no compact proper slant submanifolds in complex Euclidean spaces, there do exist compact proper slant submanifolds in complex flat tori.

A submanifold $N$ of a pseudo-Riemannian Sasakian manifold $(\tilde M,g,\phi,\xi)$ is called {\it contact $\theta$-slant} if the
structure vector field $\xi$ of $\tilde M$ is tangent to $N$ at each point of $N$ and, moreover, for each unit vector $X$ tangent to $N$ and orthogonal to $\xi$ at $p\in N$, the angle $\theta(X)$ between $\phi(X)$ and $T_pN$ is independent of the choice of $X$ and $p$.

Let $H^{2m+1}_1(-1)\subset C^{m+1}_1$ denote the anti-de Sitter space-time and $$\pi\colon\; H^{2m+1}_1(-1)\to  CH^m(-4)$$ the corresponding Hopf's  fibration. Then every $n$-dimensional proper $\theta$-slant submanifold $M$ in
$CH^{m}(-4)$ lifts to an $(n+1)$-dimensional proper contact $\theta$-slant submanifold $\pi^{-1}(M)$ in $H^{2m+1}_1(-1)$ via $\pi$.

Conversely, a proper contact $\theta$-slant submanifold of $H^{2m+1}_1(-1)$ projects to a proper $\theta$-slant submanifold of $CH^{m}(-4)$ via $\pi$. 

Similar correspondence also holds between proper $\theta$-slant submanifolds of $CP^m(4)$ and proper contact $\theta$-slant submanifolds of the Sasakian $S^{2m+1}(1)$ (see \cite{CT}).

For  K\"ahlerian $\theta$-slant submanifolds in complex space forms, we have the following inequality from \cite{c11,c21}.

\begin{theorem} Let $\phi:M\to \tilde M^n(4\epsilon),\,\epsilon\in\{-1,0,1\},$ be a K\"ahlerian $\theta$-slant submanifold of dimension $n$ in a complete simply-connected complex space form $\tilde M^n(4\epsilon)$. Then we have
\begin{align} \label{32.1} H^2\geq{{2(n+2)}\over{n^2(n-1)}}\tau-{{n+2}\over n}\left(1+{{3\cos^2\theta}\over {n-1}}
\right)\epsilon .\end{align}

The equality sign of \e{32.1} holds identically if and only if one of the following three cases occurs:
\vskip.04in 

{\rm (a)}   $\theta=0$ and $M$ is a totally geodesic complex submanifold of $\tilde M^n(4\epsilon)$, or

\vskip.04in 
{\rm (b)}  $\epsilon=0$ and $M$ is a totally geodesic $\theta$-slant submanifold in {\bf C}$^n$, or

\vskip.04in 
{\rm (c)}  $\epsilon=-1$, $n=2$, $\theta=\cos^{-1}\left({1\over 3}\right)$, and $M$ is a surface of constant curvature $-{2\over 3}$. Moreover,  up to rigid motions, the immersion $\phi$ is  the composition $\phi=\pi\circ z$, where $\pi$ is the hyperbolic Hopf fibration $\pi:H^5_1\to CH^2(-4)$ and 
$$z: \hbox{\bf R}^3\to H^5_1(-4)\subset {\bf C}^3_1$$ is the immersion defined by
\begin{equation}\begin{aligned} &z(u,v,t)=\rme^{\i t}\Bigg(\frac3 2\cosh av-\frac1 2+ \frac 1 6{u^2}\rme^{-av}-\frac \i 6 {{
\sqrt{6}u(1+\rme^{-av})}}, 
\\& \hskip.2in \frac1
3(1+2\rme^{-av})u+\i\sqrt{6}\(\frac {\rme^{av}}{4}+ \frac {\rme^{-av}} {12}+\frac {\rme^{-av}}{18}{u^2}-\frac 1 3 \),
    \\ &\hskip.3in  \frac{\sqrt{2}}6 (1-\rme^{-av})u+\i\sqrt{3}\(\frac 1 6+\frac {{\rme^{av}}} 4 +\frac  {\rme^{-av}} {18}{{u^2}}-{\frac
 {5\rme^{-av}}{12}}\)\hskip-.02in\Bigg)  \end{aligned}\end{equation} with $a=\sqrt{2/3}$.\end{theorem}

Inequality \e{32.1} has been extended by A. Oiag\u{a} \cite{oi} to the following.

\begin{theorem} Let $\phi:M\to \tilde M^n(4\epsilon),\,\epsilon\in\{-1,0,1\},$ be a purely real submanifold with $\nabla P=0$  in a  complex space form $\tilde M^n(4\epsilon)$. Then we have
\begin{align} \label{32.3} H^2\geq{{2(n+2)}\over{n^2(n-1)}}\tau-{{n+2}\over n}\left[1+{{3||P||^2}\over {n(n-1)}}\right]\epsilon .\end{align}
\end{theorem}

Here,  a purely real submanifold  is in the sense of \cite{1981},  i.e., a submanifold $M$ in $\tilde M^n(4\epsilon)$ whose holomorphic distribution is trivial, or equivalently, $\mathcal D_x=\{0\}$ for each $x\in M$.

\section{ Related articles}
 After  the invention of $\delta$-invariants, there are many  articles which study one of the topics treated in this survey.  Also,  there is a book by Adela Mihai \cite{Ma8} which also provides a  survey on this research area. 
  However, to enable further study in this very active field of research we divide those articles into several categories according to their main results and  applications.
  \vskip.1in

1.  {$\delta$-invariants:}
 \cite{c7,c12,c14,c16,c-nash,c34,c36,c37,c39,cdv,CG1,cmihai,HSV,Ma8,Op4,Sa3,Sa5,su2}.
  \vskip.1in

2.  Inequalities involving $\delta$-invariants:
 \cite{Ak,ACKY,AEMMO2,BDFV,ca1,ca2,c5,c7,c12,c14,c16,c24,c-nash,c36,c39,c40,c43,cdv,cdvv3,cmihai,ci,ci2,ci3,COi,DMV1,DMV2,GAK,HSV,HMT,KST1,KC,kk,kim,LHX,Ma3,Ma8,o,Op2,Op3,SV,su1,su2,su3,Tri,TKK2,vi,yoon1}.
  \vskip.1in

3. Related inequalities:
 \cite{bcm,ca1,ca2,ca3,c10,c11,c12,c15,c21,c22,c23,c24,c26,c27,c28,c29,c30,c31,c-nash,c33,c34,c35,c37,c38,c42,c43,cd,CJ,CT,cw,ci3,DDVV1,DDVV2,DFV1,DFV2,DHPV,GR,HMi,KST1,KTC2,kkd,Ma1,Ma2,Ma3,Ma4,Ma5,Ma6,Ma8,o,SV,SL,Win,yoon2,yoon3}.
  \vskip.1in

4. Ideal immersions-equality involving $\delta$-invariants:
 \cite
{ADV,blair,BV1,BV2,c5,c6,c12,c16,c26,c27,c33,c38,c39,c40,c41,c42,c43,cdfv,cdv,cdvv1,cdvv2,cdvv95,cdvv3,cdvv4,CG2,cmihai,CV2,CV3,d,DT,DMV1,Dillen1,DV1,DV2,DjV,DjV1,Fa,f,Hi,hong,KSV,KV,li1,li2,GWu,Lus,Ma8,Mi3,Mi5,OA,OT,MHi,Sa1,Sa2,Sa3,Sa5,SS,SV96,SV,su3,TKK3}.
 \vskip.1in

5. Equality case of related inequalities:
 \cite
{c11,c15,c19,c22,c23,c26,c27,c28,c30,c33,c38,c40,c41,c42,c43,CG2,CT,CV1,CL,DDVV1,DDVV2,DFV2,Ma8,Mi4,MT,PV,su3}.
 \vskip.1in

6. Applications to  immersions:
 \cite{ca1,ca2,c12,c16,c26,c29,c34,c37,c39,KY}.
\vskip.1in

7. General warped products:
\cite{BD,c26,c28,c29,c31,c33,c38,c42,cd,cw,f,kkd,MM,Ma4,Ma6,Ma7,Ma8,Mi2,Mi04,su2,yoon4,yoonC}. 
\vskip.1in

8. $CR$-products: \cite{bejancu,c2,c3,book1981,c18,c26,clm,CV2,Ma8,oi}

\vskip.1in
9. $CR$-warped products: 
\cite{Al, AEMMO1,BD,BM,c20,c22,c23,c28,c30,c35,cd,HMi,kkd,Ma8,Mi2,Mi04,Mu1,Mu2,Sa, SaG}. 
\vskip.1in

10. Ricci curvature:
\cite{AEMMO3,AEMMO4,c16,c17,c27,c31,FH,HT1,HT2,KDM,KTC2,kimP,Liu1,Liu2,LD,MMO2,Mi1,MiSS,oi2,Sa4,ST,WL,yoon2,yoon5}.  
  \vskip.1in

11. Shape operator: 
\cite{c8,c13,c16,ci3,FH,KST2,LS,LWS,Lu,MMO1,Ma2,Ma8,TKK1,YB}. 
 \vskip.1in

12. Contact  and Sasakian manifolds:
\cite{ACKY,cmihai,ci,ci2,ci3,DMV1,DMV2,FH,KC,kk,kim,LHX,Liu1,LS,LWS,MM,Ma8,Mi1,Mi2,Mi3,MT,Tri,TKK1,TKK2,TKK3,yoon3,yoon4,yoon5,yoonC}. 
   \vskip.1in

13. Affine differential geometry:
\cite{BDFV,c38,c43,c44,cdv,DV2,KSV,KV,Lus,SS,SV96}.
      \vskip.1in

14.  Lagrangian submanifolds of K\"ahler manifolds:
\cite{BDFV,BMV,bo1,bo2,BV1,bcm,ca3,c2,c3,c96,c16,c17,c19,c25,c41,cdvv1,cdvv3,CV1,CV3,Ma8,Op3}
\vskip.1in

15. Slant submanifolds:
\cite{ACKY,ca1,ca2,cbook3,c15,c21,COi,FH,GAK,kk,GWu,MM,MMO1,Ma5,Mi4,MT,Sa,SL,vi,yoon5}.
  \vskip.1in

16.  Other applications:
\cite{c24,c36,cdvv4,Hae,HSV,HV,SV}.

\vskip.4in
{\bf Acknowledgements:} The author would like to express his many thanks to Professors D. E. Blair, F. Dillen, I. Dimitric, J. Fastenakels, S. Haesen, I. Mihai, B. Suceav\u{a}, J. Van der Veken, L. Verstraelen, L. Vrancken and S. W. Wei for their valuable suggestions for the improvement of the presentation of this article.

\end{document}